\documentclass[preprint,sort&compress,12pt]{elsarticle}
\usepackage{amsmath}
\usepackage{mathrsfs}
\usepackage{stmaryrd}
\usepackage{mathrsfs}
\usepackage{bm}
\usepackage{amsmath}
\usepackage{amsfonts}
\usepackage{amssymb}
\usepackage{amsthm}

\journal{\quad}
\newcommand{\mf}{\mathbf}
\newcommand{\mm}{\mathrm}

\begin{document}
\begin{frontmatter}
\title{Nonlinear Rayleigh-Taylor Instability for Nonhomogeneous 
\\  Incompressible Viscous Magnetohydrodynamic Flows
}

\author[FJ]{Fei Jiang\corref{cor1}
}
\ead{jiangfei0591@qq.com}
 \cortext[cor1]{Corresponding.
 }
\author[sJ]{Song Jiang}
\ead{jiang@iapcm.ac.cn}
\author[FJ]{Weiwei Wang}
\ead{wei.wei.84@163.com}
\address[FJ]{College of Mathematics and Computer Science, Fuzhou University, Fuzhou, 350108, China.}
\address[sJ]{Institute of Applied Physics and Computational Mathematics, Beijing, 100088, China.}
\begin{abstract}
We investigate the nonlinear instability of a smooth
Rayleigh-Taylor steady-state solution (including the case of heavier density with increasing
height) to the three-dimensional incompressible nonhomogeneous
 magnetohydrodynamic (MHD) equations of zero resistivity in the presence of a uniform gravitational field.
We first analyze the linearized equations around the steady-state solution.
Then we construct solutions of the linearized problem that grow in time in the Sobolev space $H^k$,
thus leading to the linear instability. With the help of the constructed
unstable solutions of the linearized problem and a local well-posedness result of smooth solutions
to the original nonlinear problem, we establish the instability of the density, the horizontal and vertical velocities in the nonlinear problem.
Moreover, when the steady magnetic field is vertical and small, we prove the instability of the magnetic field. This verifies
the physical phenomenon: instability of the velocity leads to the instability of the magnetic field through the induction equation.
\end{abstract}

\begin{keyword} Incompressible MHD flows, steady solutions, Rayleigh-Taylor instability,
Navier-Stokes equations. \MSC[2000] 35Q35\sep  76D03.

\end{keyword}
\end{frontmatter}


\newtheorem{thm}{Theorem}[section]
\newtheorem{lem}{Lemma}[section]
\newtheorem{pro}{Proposition}[section]
\newtheorem{cor}{Corollary}[section]
\newproof{pf}{Proof}
\newdefinition{rem}{Remark}[section]
\newtheorem{definition}{Definition}[section]

\section{Introduction}
\label{Intro} \numberwithin{equation}{section}

This paper is concerned with nonlinear instability of a smooth Rayleigh-Taylor (RT) steady-state solution
to the following three-dimensional (3D) nonhomogeneous incompressible magnetohydrodynamic
(MHD) equations with zero resistivity (i.e. without magnetic diffusivity)
in the presence of a uniform gravitational field 
(see, for example, \cite{CHTA1,kulikovskiy1965magnetohydrodynamics,landau1984electrodynamics,LFZPGC}
on the derivation of the equations):
\begin{equation}\label{0101}\left\{\begin{array}{l}
 \rho_t+{\bf v}\cdot\nabla \rho=0,\\[1mm]
\rho\mathbf{v}_t+\rho {\bf v}\cdot\nabla {\bf v}+\nabla p=(\nabla\times \mf{M})\times \mf{M}+\mu\Delta
{\bf v}-\rho{g}\mf{e}_3,\\[1mm]
\mf{M}_t-\nabla\times (\mf{v}\times \mf{M})=\mf{0},\\[1mm]
\mathrm{div}\mathbf{v}=0,\quad \mathrm{div}\mathbf{M}=0.\end{array}\right.
\end{equation}
Here the unknowns $\rho:=\rho(t,\mf{x})$, $\mathbf{v}:=\mathbf{v}(t,\mf{x})$,
$\mf{M}:=\mf{M}(t,\mf{x})$ and $p:=p(t,\mf{x})$ denote the density,
velocity, magnetic field and pressure of the incompressible fluid, respectively;
$\mu>0$ stands for the coefficient of shear viscosity, $g>0$ for the gravitational
constant, $\mf{e}_3=(0,0,1)$ for the vertical unit vector, and $-\rho g\mf{e}_3$ for the gravitational force.

In the system \eqref{0101} the equation \eqref{0101}$_1$ is the continuity equation, while
\eqref{0101}$_2$ describes the balance law of momentum. It is well-known that the electromagnetic
field is governed by the Maxwell equations. In MHD, the displacement current can be neglected
\cite{kulikovskiy1965magnetohydrodynamics,landau1984electrodynamics}.
As a consequence, \eqref{0101}$_3$ is called the induction equation.
As for the constraint $\mm{div}\,\mf{M}=0$, it can be seen just as a
restriction on the initial value of $\mf{M}$ since
$(\mm{div}\,\mf{M})_t=0$. We remark that, the resistivity in
\eqref{0101}$_3$ is zero, which arises in the physics regime with
negligible electrical resistance, see
\cite{cowling1957magnetohydrodynamics}. In addition, if
$\mf{M}\equiv0$, the system \eqref{0101} reduces to the incompressible
Navier-Stokes equations in the presence of a uniform gravitational field.

In this paper we consider the problem of the RT instability in a horizontally periodic domain
$\Omega:=(2\pi L\mathbb{T})^2\times \mathbb{R}$, where $2\pi L \mathbb{T}$ stands for the 1D-torus
of length $2\pi L$. We assume that a smooth RT (steady-state) density profile
$\bar{\rho}:=\bar{\rho}({x}_3)\in L^\infty(\mathbb{R})$ exists and satisfies
\begin{eqnarray}\label{0102}
&&\bar{\rho}'\in C_0^{\infty}(\mathbb{R}),\quad \inf_{x_3\in \mathbb{R}}\bar{\rho}>0,\\[0em]
&& \label{0103}\bar{\rho}'(x_3^0)>0\;\;\mbox{ for some point }x_3^0\in \mathbb{R},
\end{eqnarray}
where $'=d/dx_3$.
We refer to \cite[Remark 1.1]{NJTSC2} for the construction of such $\bar{\rho}$.
Let $\bar{\mf{M}}\in \mathbb{R}^3$ be a constant magnetic field, then
the RT density profile $\bar{\rho}$ with $(\mathbf{v},\mathbf{M})(t,\mathbf{x})
\equiv(\mathbf{0},\bar{\mathbf{M}})$ defines a steady state solution
to \eqref{0101}, provided the steady pressure $\bar{p}$ is determined by
\begin{equation*}
\nabla \bar{p}=-\bar{\rho}g \mf{e}_3,\;\mbox{ i.e., }\;\frac{d\bar{p}}{dx_3}=-\bar{\rho}g.
\end{equation*}
We point out that by virtue of the condition \eqref{0103}, there is at least
a region in which the RT density profile has larger density
with increasing $x_3$ (height), thus leading to the classical
RT instability as shown in Theorem \ref{thm:0101} below.

 Now, we denote the perturbation to the RT steady state by
$$ \varrho=\rho -\bar{\rho},\quad \mathbf{u}=\mathbf{v}-\mathbf{0},\quad
\mathbf{N}=\mf{M}-\bar{\mf{M}},\quad q=p-\bar{p},$$
then, $(\varrho ,\mathbf{u},q)$ satisfies the perturbed equations
\begin{equation}\label{0105}\left\{\begin{array}{l}
\varrho_t+{\bf u}\cdot\nabla (\varrho+\bar{\rho})=0, \\[1mm]
(\varrho+\bar{\rho}){\bf u}_t+(\varrho+\bar{\rho}){\bf u}\cdot\nabla
{\bf u}+\nabla q+g \varrho \mf{e}_3=\mu \Delta \mathbf{u}+(\nabla\times \mf{N})\times (\mf{N}+\bar{\mf{M}}),\\[1mm]
\mf{N}_t=\nabla\times (\mf{u}\times (\mf{N}+\bar{\mf{M}})),\\[1mm]
\mathrm{div}\mathbf{u}=0,\ \mathrm{div}\mathbf{N}=0.\end{array}\right.  \end{equation}
For (\ref{0105}) we impose the initial and boundary conditions:
\begin{equation}\label{0106}
(\varrho,{\bf u},\mathbf{N} )|_{t=0}=(\varrho_0,{\bf u}_0,\mf{N}_0)\quad\mbox{in } \Omega
\end{equation}
and
\begin{equation}\label{0107}
\lim_{|{x}_3|\rightarrow +\infty}(\varrho,{\bf u},\mf{N})(t,\mathbf{x}',x_3)={\bf 0}\quad
\mbox{ for any }t>0,
\end{equation}where we have written $\mf{x}'=x_1\mf{e}_1+x_2\mf{e}_2$, $\mf{e}_1:=(1,0,0)$ and $\mf{e}_2:=(0,1,0)$.
Moreover, the initial data should satisfy the compatibility conditions
\begin{equation*}
\mathrm{div}\mathbf{u}_0=0\quad\mbox{ and }\quad\mathrm{div}\mathbf{N}_0=0.
\end{equation*}

If we linearize the equations \eqref{0105} around the steady state
$(\bar{\rho},\mathbf{0},\bar{\mf{M}})$, then the resulting linearized equations read as
\begin{equation}\label{0108}
\left\{\begin{array}{ll}
 \varrho_t+\bar{\rho}'{u}_3=0, \\[1mm]
  \bar{\rho}\mathbf{u}_t +\nabla q+g\varrho \mf{e}_3=\mu \Delta \mathbf{u}
  +(\nabla\times \mf{N})\times \bar{\mf{M}},\\[1mm]
  \mf{N}_t=\nabla\times (\mf{u}\times \bar{\mf{M}}),\\[1mm]
 \mathrm{div}\mathbf{u}=0,\quad \mathrm{div}\mathbf{N}=0.
\end{array}\right.\end{equation}

The RT instability is well-known as gravity-driven instability in fluid dynamics when heavy
fluid is on top of light one. The linear instability for an incompressible fluid was first introduced
by Rayleigh in 1883 \cite{RLAP}. The analogue of the RT instability
arises when fluids are electrically conducting and a magnetic
field is present, and the growth of the instability will be
influenced by the magnetic field due to the generated
electromagnetic induction and the Lorentz force. Some authors have
extended the partial results concerning the RT instability to the case of MHD fluids
by overcoming difficulties induced by presence of the magnetic field.
For example, Kruskal and Schwarzchild in 1954 first showed that a horizontal magnetic field
has no effect on the development of the linear RT
instability \cite{KMSMSP}. Then the influence of a vertical magnetic field was investigated
by Hide in \cite{HRWP} where the effect of finite viscosity and resistivity was included
and his analysis was encumbered with many parameters. By a variational approach,
Hwang in 2008 studied the nonlinear RT instability of \eqref{0105}--\eqref{0107}
for the inviscid case (i.e. $\mu=0$) in a 2D periodic domain \cite{HHVQ}.

 To our best knowledge, however, it is still open mathematically whether there exists an unstable
solution to the nonlinear RT problem \eqref{0105}--\eqref{0107} of 3D viscous MHD fluids.
The aim of this article is to rigorously verify the instability for the nonlinear RT
 problem \eqref{0105}--\eqref{0107}. Moreover, the impact of the magnetic filed on the instability will be analyzed,
 for example, we shall show that if the steady magnetic field is vertical and small, then the magnetic field is unstable,
 thus verifying the physical phenomenon: instability of the velocity leads to the instability of the magnetic field
 through the induction equation.
 The main result of this paper reads as follows.
\begin{thm}\label{thm:0101}
 Let
  \begin{equation}\begin{aligned}\label{0109}
M_\mathrm{c}:=\sqrt{\sup_{\psi\in H^1(\mathbb{R})\atop \psi \equiv\!\!\!\!/~0}
\frac{\int_{\mathbb{R}}g\bar{\rho}'|\psi|^2\mm{d}x}{\int_{\mathbb{R}}|\psi'|^2\mm{d}x}}>0.  \quad\qquad\mbox{(critical number)}
\end{aligned}\end{equation}
 Assume that the RT density profile $\bar{\rho}\in L^\infty(\Omega)$ satisfies \eqref{0102} and \eqref{0103}, and
 \begin{equation}\label{mageticfields}\mf{\bar{M}}=\left\{
                                \begin{array}{ll}
M\mf{e}_1, & M\neq 0, \mbox{ constant}, \\[1mm]
M\mf{e}_3, & |M|\in (0,M_\mm{c}), \mbox{ constant}.
              \end{array}
            \right.
\end{equation}
 Then, the steady state $(\bar{\rho},\mathbf{0},\mf{\bar{M}})$ of \eqref{0105}--\eqref{0107}
 is unstable, that is, there exist positive constants $\Lambda^*$, $\varepsilon$, $m_0$ and $\iota$,
 and a triple $(\bar{\varrho}_0,\bar{\mf{u}}_0,\bar{\mf{N}}_0)\in
H^\infty(\Omega):=\cap_{k=0}^\infty H^k(\Omega)$,
such that for any $\delta\in (0,\iota)$ and the initial data
 $(\varrho_0,\mathbf{u}_0,\mathbf{N}_0):=(\delta\bar{\varrho}_0,\delta\bar{\mf{u}}_0,\delta\bar{\mf{N}}_0)$,
there is a unique classical solution $({\varrho},\mathbf{u},\mf{{N}})$ of \eqref{0105}--\eqref{0107} on
$[0,T^{\max})$, but
\begin{equation}\label{0111}\|\varrho(T^\delta)\|_{L^2(\Omega)},\ \|({u}_1,u_2)(T^\delta)\|_{L^2(\Omega)},\
\|{u}_3(T^\delta)\|_{L^2(\Omega)}\geq {\varepsilon}\;
\end{equation}  for some escape time
$T^\delta:=\frac{1}{\Lambda^*}\mm{ln}\frac{2\varepsilon}{m_0\delta}\in
(0,T^{\max})$, where $\mathrm{div}\bar{\mathbf{u}}_0=0$, $\mathrm{div}\bar{\mathbf{N}}_0=0$, and $T^{\max}$ denotes
the maximal time of existence of the solution $(\varrho,\mathbf{u},\mf{N})$. Moreover,
\begin{equation*} \|N_3(T^\delta)\|_{L^2(\Omega)}\geq \varepsilon\quad
\mbox{ for the case }\;\bar{\mf{M}}=M\mf{e}_3.
\end{equation*}
\end{thm}
\begin{rem}
We should point out that since the above Sobolev spaces $H^k(\Omega)$ are
defined on the horizontally periodic domain $\Omega=(2\pi L\mathbb{T})^2\times \mathbb{R}$,
the solution $({\varrho},\mathbf{u},\mf{{N}})$ constructed in Theorem \ref{thm:0101}
is horizontally periodic. One should keep this in mind in what follows.

Throughout the rest of this article we shall repeatedly use the abbreviations:
\begin{align} &
I_T:=(0,T),\ \bar{I}_T:=[0,T],\ W^{m,p}:=W^{m,p}(\Omega),\;\; H^m:=H^m(\Omega),\
H^\infty:=\cap_{m=1}^\infty H^{m}(\Omega),\nonumber\\
&L^p:=L^p(\Omega),\ \|\cdot\|_{W^{m,p}}:=\|\cdot\|_{W^{m,p}(\Omega)},\;\;
\|\cdot\|_{H^m}:=\|\cdot\|_{H^m(\Omega)},\;\;
\|\cdot\|_{L^p}:=\|\cdot\|_{L^p(\Omega)},\quad \mbox{etc.}\nonumber
\end{align}
\end{rem}
\begin{rem}The conclusion \eqref{0111} in Theorem \ref{thm:0101} also holds for the general
horizontal magnetic field $\bar{\mf{M}}=({M}_1,{M}_2,0)$. In fact,
rotating the $o$-$xy$ coordinates properly so that
$\bar{\mf{M}}=({M},0,0)$ with ${M}=\sqrt{{M}_1^2+{M}_1^2}$,
we have the same case as in Theorem \ref{thm:0101} under the rotated coordinates,
since the $L^2$-norms of  the density, horizontal velocity, and vertical velocity are
invariant under the horizontal rotation.
We should point out here that for the both cases of the horizontal and
vertical magnetic fields, the ordinary differential equation (ODE)
\eqref{0207} associated with the normal mode solutions enjoys good variational structure.
For a general constant magnetic field $\bar{\mf{M}}=({M}_1,{M}_2,M_3)$, similarly to that in the derivation of \eqref{0207}, we can also
deduce a ODE corresponding to $\bar{\mf{M}}=({M}_1,{M}_2,M_3)$. The resulting
ODE, however, possesses the terms $\mm{i}\psi'$ and $\mm{i}\psi'''$ which destroy the good variational structure, and consequently,
one could not directly construct the growing mode solutions to the linearized problem.
In addition, Theorem \ref{thm:0101}
also holds when $\bar{\mf{M}}=\mf{N}=\mf{0}$, hence the weak RT instability for incompressible viscous fluids given
 in \cite{NJTSC2} can be further shown to be strongly RT unstable as in \eqref{0111} in the case of the horizontally periodic domain.
\end{rem}

\begin{rem}
In \cite{HHJGY} Hwang and Guo proved the nonlinear RT
 instability for 2D nonhomogeneous incompressible inviscid flows
 (i.e. $\mu=0$ and $\bar{\mf{M}}=\mf{N}=\mf{0}$ in \eqref{0105}, \eqref{0106}) with boundary condition
 $\mathbf{u}\cdot \mathbf{n}|_{\partial\Omega'}=0$ where
$\Omega'=\{(x_1,x_2)\in \mathbb{R}^2~|~-l<x_2<m\}$ and $\mathbf{n}$ denotes the
 outer normal vector to $\partial\Omega'$. Later, Hwang \cite{HHVQ} further investigated the nonhomogeneous incompressible inviscid MHD
 fluid on a periodic domain, and get the instability of the norm $\|(\varrho,\mathbf{u},\mathbf{N})\|_{L^2(\Omega')}$. Our result is
 more precise than those in \cite{HHJGY,HHVQ} in the sense that Theorem \ref{thm:0101} reveals that the vertical velocity induces the instability of
  the density and horizontal velocity; and moreover, we can show the instability of the magnetic field for the case $\bar{\mf{M}}=M\mf{e}_3$.
This mathematically verifies the physical phenomenon: instability of the velocity further leads to the instability of the magnetic field
through the induction equation.
\end{rem}
\begin{rem}\label{criticalmget}
The number $M_\mm{c}$ in \eqref{0109} is infinite for $\int_{\mathbb{R}}\bar{\rho}'\mm{d}x> 0$
(i.e., $\bar{\rho}(+\infty)>\bar{\rho}(-\infty)$) and is finite for $\int_{\mathbb{R}}\bar{\rho}'\mm{d}x< 0$
(i.e., $\bar{\rho}(+\infty)<\bar{\rho}(-\infty)$), see Proposition \ref{pro:n0201} for the detailed proof.
This means that any vertical steady magnetic field cannot restrain growth of the nonlinear RT instability for
the case $\int_{\mathbb{R}}\bar{\rho}'\mm{d}x> 0$. Now it is still not
clear to us that whether a sufficiently large vertical steady magnetic field has impact on growth of the nonlinear RT instability
 for the case $\int_{\mathbb{R}}\bar{\rho}'\mm{d}x< 0$ due to some technical difficulties.
 However, if we consider $\Omega =(2\pi L\mathbb{T})^2\times (-l,m)$,
then, similarly to the derivation of (3.71) in \cite{WYC}, we can show the stability of the velocity for any classical solution
of the linearized problem satisfying boundary conditions $\mf{u}|_{x_3=-l}=\mf{u}|_{x_3=m}=0$, provided
the vertical steady magnetic field is sufficiently large. This result does not contradict that in \cite{HHVQ} where for any vertical steady
magnetic field, the nonlinear RT instability for a nonhomogeneous incompressible inviscid MHD flow in a periodic domain (the vertical direction
is also periodic) is shown. These mathematical results reveal that the domain and boundary conditions of the velocity do have impact on
the instability of MHD flows.
\end{rem}

The proof of Theorem \ref{thm:0101} is divided into four steps given in Sections \ref{sec:02}--\ref{sec:05}:
(i) First, we make an ansatz to seek the ``normal mode'' solutions of the linearized
equations \eqref{0108}, which grow exponentially in time by the factor $e^{\lambda(\xi)t}$ with $\xi\in\mathbb{R}^2$ being the
horizontal spatial frequency and $\lambda (\xi)>0$. This reduces
the equations to a system of ODEs defined on $\mathbb{R}$ with
$\lambda(\xi)>0$ for some $\xi$ (see \eqref{0205}, \eqref{0206}). All such points $\xi$ constitute a solvable domain.
Because of presence of the magnetic field, the solvable domain is not a ball for the horizontal case as in \cite{GYTI2},
resulting in some difficulties in constructing a solvable domain. In order to circumvent such difficulties,
similarly to \cite{WYC,JFJSWWWOA}, we introduce the notions of the critical frequency function $S(\xi)$
and critical frequency constant $|\xi|_{\mathrm{vc}}^M$ to define the solvable domain $\mathbb{A}^{\mm{g}}$
(cf. \eqref{0217}, \eqref{0216}). Thus, by careful constructing the solutions and
adapting the modified variational method in \cite{GYTI2},
we can also solve the ODEs for any given $\xi\in \mathbb{A}^{\mm{g}}$ and obtain a normal mode with $\lambda(\xi)>0$,
thus leading to a mechanism for the global linear RT instability.
Using the normal modes, we can further construct real-valued solutions of the linearized problem \eqref{0106}--\eqref{0108} that
grow in time, when measured in $H^k(\Omega)$ for any $k\geq 0$.
In particular, the density, horizontal velocity and vertical velocity in the solutions of the linearized problem are not zero, this fact
will play a key role in the nonlinear instability in \eqref{0111}.
(ii) In Section \ref{sec:04}, we state a local well-posedness result of the perturbed problem \eqref{0105}--\eqref{0107},
which will be proved in Section \ref{appendix}. Then we derive the integrand form of Gronwall's inequality of high-order energy
estimate $\mathcal{E}(\varrho,\mathbf{u},\mathbf{N})$ for the perturbed problem, and this
makes the escape time occurring before break-down of the classical solutions.
Since the equilibrium state of the magnetic field $\bar{\mf{M}}$ is a no-zero vector,
we shall introduce a simple technique to deal with the terms including $\bar{\mf{M}}$
in the energy estimates, see Subsection \ref{sec:0404}. (iii) Finally, in Section \ref{sec:05},
with the help of the results established in Sections \ref{sec:02}--\ref{sec:04}, we adapt a careful bootstrap argument as in \cite{GYSWIC}
 to establish the instability of the nonlinear problem. We mention that although
the approach in \cite{GYSWIC} has been also used to treat the instability of other problems
(see \cite{GYSWIC,GYHCSDDC,JJTIIA,wang2011viscous} for examples), but Duhamel's principle in the standard bootstrap argument
can not be directly applied to our problem,
 since the nonlinear terms in \eqref{0105} do not satisfy the compatibility condition of divergence-free.
 To circumvent this obstacle, we shall use some specific energy estimates to replace Duhamel's principle. Moreover,
we can also find in the proof that $\Lambda$ is indeed a sharp exponential growth rate for
general solutions to the linearized problem (see Remark \ref{rem0401n}).

We end this section by briefly reviewing some of the
previous results on the nonlinear RT instability for two layer incompressible
fluids separated by a free interface (stratified fluids), where the RT steady-state solution is
a denser fluid lying above a lighter one separated by a free interface.
When the densities of two layer fluids are two constants,
Wang and Tice \cite{wang2011viscous} proved the (local)
existence of nonlinear unstable solutions in a horizontally periodic domain
$\mathbb{T}^2\times (-b,1)$, where the instability term is described
by the sum of $L^{2}$-norm of the velocity and the moving internal interface.
Pr\"uss and Simonett used the $C^0$-semigroup theory and the Henry instability theorem to show the
existence of nonlinear unstable solutions in the Sobolev-Slobodeckii spaces in $\mathbb{R}^3$ \cite{PJSGOI5},
where the instability term is described by the sum of $\|\mathbf{u}\|_{W_p^{2-2/p}(\mathbb{R}^3)}$ and
$\|\mathbf{h}\|_{W_p^{3-2/p}(\mathbb{R}^3)}$ (see \cite[Theorem 1.2]{PJSGOI5} for
details).  When densities of two layer fluids are variable, to our best knowledge, the (local) existence
of solutions to the nonlinear problem  \eqref{0105}--\eqref{0107} still remains open, consequently,
the strong nonlinear instability is still open. For compressible fluids, there are very few results on
the nonlinear RT instability. Guo and Tice proved the instability
of immiscible compressible inviscid fluids in the frame of Lagrangian coordinates
under the existence assumption of solutions \cite{GYTI1}. This is in some sense a compressible analogue to the local
ill-posedness of the RT problem for incompressible fluids given in \cite{EDGTC1111}.

Finally we mention related results on the instability for stratified MHD fluids.
Wang \cite{WYC} introduced the critical number of stratified MHD fluids
(denoted by $M_\mathrm{c}^{\mathrm{s}}$) to investigate the linear
RT instability of stratified MHD fluids in a infinite slab domain $\mathbb{R}^2\times (-l,l)$.
Later, Jiang et al. \cite{JFJSWWWOA} further
showed the weak nonlinear RT instability of stratified MHD fluids and found that
$M_\mathrm{c}^{\mathrm{s}}=\sqrt{gl(\varrho_+-\varrho_-)/2}$, where $\varrho_+> \varrho_-$,
and $\varrho_+ $ resp. $\varrho_-$ denotes the density of the upper- resp. lower-layer fluid.
Obviously, $M_\mathrm{c}^{\mathrm{s}}\to\infty$ as the height $2l\rightarrow \infty$. Hence,
the critical number of stratified MHD fluids in $\mathbb{R}^3$ is infinite. This
means that the vertical magnetic field may have no impact on growth of the RT instability in $\mathbb{R}^3$.

\section{Construction of solutions to the linearized problem}\label{sec:02}
We wish to construct a solution to the linearized equations
\eqref{0108} that has growing $H^k$-norm for any $k$. We will
construct such solutions via some synthesis as in \cite{GYTI2} by first constructing
a growing mode for any but fixed spatial frequency. Moreover, we shall introduce
the techniques of the critical frequency function and critical frequency constant
 to carefully analyze the terms involving with the magnetic field.

\subsection{Linear growing modes}
To start with, we make a growing mode ansatz of solutions, i.e., for some $\lambda>0$,
\begin{equation*}
{\varrho}(t,\mathbf{x})=\tilde{\rho} (\mathbf{x})e^{\lambda t},\;
\mathbf{u}(t,\mathbf{x})=\tilde{\mathbf{v}}(\mathbf{x})e^{\lambda t},\;
{q}(t,\mathbf{x})=\tilde{p}(\mathbf{x})e^{\lambda t},\;
\mf{N}(t,\mf{x})=\tilde{\mf{M}}(\mathbf{x})e^{\lambda t}.
\end{equation*}
Substituting this ansatz into \eqref{0108}, and then eliminating
$\tilde{\rho}$ and $\tilde{\mf{M}}$ by using the first and third equations, we arrive at the
time-independent system for
$\tilde{\mathbf{v}}=(\tilde{v}_1,\tilde{v}_2,\tilde{v}_3)$ and $\tilde{p}$:
\begin{equation}\label{0201}
\left\{
                              \begin{array}{ll}
\lambda^2\bar{\rho}\tilde{\mathbf{v}}+\lambda\nabla \tilde{{p}}
= (\nabla \times (\nabla\times (\tilde{\mf{v}}\times \bar{\mf{M}})))\times \bar{\mf{M}}+\lambda\mu \Delta \tilde{\mathbf{v}}+g\bar{\rho}'\tilde{v}_3\mf{e}_3,\\[1mm]
\mathrm{div}\, \tilde{\mathbf{v}}=0
\end{array}
                            \right.
\end{equation}
with
\begin{equation}\label{0202}
\lim_{|x_3|\rightarrow
+\infty}(\tilde{\mathbf{v}},\tilde{\mathbf{M}})(\mathbf{x'},x_3)=\mathbf{0},
\end{equation}where $\bar{\mf{M}}$ is given by \eqref{mageticfields}.
After a straightforward calculation, we find that
\begin{equation}\label{n0203}\begin{aligned}
\nabla \times (\nabla\times (\tilde{\mf{v}}\times \bar{\mf{M}}))\times \bar{\mf{M}}
&=\nabla \times (\bar{\mf{M}}\cdot\nabla\tilde{\mf{v}} )\times \bar{\mf{M}}\\
&=M^2\left\{
                                 \begin{array}{ll}
(0,\partial_{1}^2 \tilde{v}_2-\partial_{21}^2\tilde{v}_1,\partial_{11}^2\tilde{v}_3
-\partial_{31}^2\tilde{v}_1),
& \hbox{for }\bar{\mf{M}}=M\mf{e}_1; \\[1mm]
(\partial_{3}^2 \tilde{v}_1-\partial_{13}^2\tilde{v}_3,\partial_{33}^2\tilde{v}_2
-\partial_{23}^2\tilde{v}_3,0),
& \hbox{ for }\bar{\mf{M}}=M\mf{e}_3.
                                 \end{array}
                                  \right.
\end{aligned}\end{equation}

 We fix a spatial frequency $\xi=(\xi_1,\xi_2)\in \mathbb{R}^2$, and
define the new unknowns
\begin{equation*}
\tilde{v}_1(\mf{x})=-\mathrm{i}\varphi(x_3)e^{\mm{i}\mf{x}'\cdot\xi}, \;\;
\tilde{v}_2(\mf{x})=-\mathrm{i}\theta(x_3)e^{\mm{i}\mf{x}'\cdot\xi},\;\;\tilde{v}_3(\mf{x})=\psi(x_3)e^{\mm{i}\mf{x}'\cdot\xi},\;\;
{{\tilde{p}}}(\mf{x})=\pi(x_3)e^{\mm{i}\mf{x}'\cdot\xi}.
\end{equation*}Then, in view of \eqref{0201}--\eqref{n0203}, we see that $\varphi$, $\theta$, $\psi$ and
$\lambda$ satisfy the following system of ODEs:
\begin{equation}\label{0205} \left\{
                              \begin{array}{ll}
\lambda^2\bar{\rho} \varphi-\lambda\xi_1\pi+\lambda\mu (|\xi|^2\varphi-\varphi'')=M^2B_1,\\[1mm]
\lambda^2\bar{\rho} \theta-\lambda\xi_2\pi+\lambda\mu (|\xi|^2\theta-\theta'')=M^2B_2,\\[1mm]
\lambda^2\bar{\rho} \psi+\lambda\pi'+\lambda\mu (|\xi|^2\psi-\psi'')-g\bar{\rho}'\psi=M^2B_3,\\[1mm]
\xi_1\varphi+\xi_2\theta+\psi'=0
\end{array}
                            \right.
\end{equation}
 with
\begin{eqnarray}\label{0206}
\varphi(-\infty)=\theta(-\infty)=\psi(-\infty)=\varphi(+\infty)=\theta(+\infty)=\psi(+\infty)=0,
\end{eqnarray}
where
\begin{equation*}  \left( \begin{array}{c}
B_1\\[1mm]
B_2\\[1mm]
B_3
\end{array}\right)=\left\{\begin{array}{ll}
\left(  \begin{array}{c} 0    \\[1mm]
          (\xi_1\xi_2\varphi-\xi_1^2\theta)\\[1mm]
 -(|\xi_1|^2\psi+\xi_1\varphi') \end{array} \right)
& \mbox{ for }\bar{\mf{M}}=M\mf{e}_1,
 \\[10mm]
\left(
                              \begin{array}{c}
(\varphi''+\xi_1\psi')\\[1mm]
(\theta''+\xi_2\psi')\\[1mm]
0
\end{array}
                            \right)&\mbox{ for }\bar{\mf{M}}=M\mf{e}_3.\end{array}\right.
\end{equation*}

Eliminating $\pi$ from the third equation in \eqref{0205}, we obtain
the following ODE for $\psi$
\begin{equation}\label{0207}\begin{aligned}
&-\lambda^2[\bar{\rho}|\xi|^2\psi-(\bar{\rho}\psi')']\\
=&\lambda\mu
(|\xi|^4\psi-2|\xi|^2\psi''+\psi'''') -g|\xi|^2\bar{\rho}'\psi + M^2\left\{
                                                                \begin{array}{ll}
\xi_1^2(|\xi|^2\psi-\psi''), & \hbox{ for }\bar{\mf{M}}=M\mf{e}_1, \\[1mm]
(\psi''''-|\xi|^2\psi''), & \hbox{ for }\bar{\mf{M}}=M \mf{e}_3,
                                                                \end{array}
                                                              \right.
\end{aligned}
\end{equation}
with
\begin{eqnarray}\label{0208}
&&\psi(-\infty)=\psi'(-\infty)=\psi(+\infty)=\psi'(+\infty)=0.
\end{eqnarray}

Next we use the modified variational method to construct a solution of \eqref{0207}, \eqref{0208}. This idea
can be found in Guo and Tice's paper on compressible viscous stratified
flows \cite{GYTI2}, and has been adapted by other authors to study
the instability for other fluid models \cite{JFJSWWWO,JJTIIA,WYC,DRJFJSRS}.
To this end, we now fix a non-zero vector $\xi\in \mathbb{R}^2$ and $s>0$. From
\eqref{0207} and \eqref{0208} we get a family of the modified problems
\begin{equation} \label{0209}\begin{aligned}
&-\lambda^2[\bar{\rho}|\xi|^2\psi-(\bar{\rho}\psi')']\\
& = s\mu (|\xi|^4\psi-2|\xi|^2\psi''+\psi'''')-g|\xi|^2\bar{\rho}'\psi+M^2\left\{
                                                                \begin{array}{ll}
 \xi_1^2(|\xi|^2\psi-\psi''), & \hbox{ for }\bar{\mf{M}}=M\mf{e}_1, \\[1mm]
(\psi''''-|\xi|^2\psi''), & \hbox{ for }\bar{\mf{M}}=M\mf{e}_3,
                                                                \end{array}
                                                              \right.
\end{aligned}\end{equation}
coupled with the condition \eqref{0208}. We define the energy functional of \eqref{0209} by
\begin{equation}\label{0210}E(\psi)=|\xi|^2E_0(\psi)+s{E}_1(\psi)\end{equation} with an associated
admissible set
\begin{equation}\label{0211}
\mathcal{A}=\left\{\psi\in
H^2(\mathbb{R})~\bigg|~J(\psi):=\int_{\mathbb{R}}\bar{\rho}(|\xi|^2|\psi|^2+|\psi'|^2)\mathrm{d}x_3=1\right\},
\end{equation}
where
\begin{eqnarray}\label{0212}
&&E_0(\psi)=\left\{
               \begin{array}{ll}
\int_{\mathbb{R}}M^2\xi_1^2\left(|\psi|^2
+\frac{|\psi'|^2}{|\xi|^2}\right)
-g\bar{\rho}'\psi^2\mm{d}x_3, & \hbox{ for }\bar{\mf{M}}=M \mf{e}_1, \\[3mm]
\int_{\mathbb{R}}M^2\left(|\psi'|^2+\frac{|\psi''|^2}{|\xi|^2}\right)
-g\bar{\rho}'\psi^2\mm{d}x_3, & \hbox{ for }\bar{\mf{M}}=M \mf{e}_3,
               \end{array}
             \right.\\[2mm]
&&\label{n03010212} E_1(\psi)= \mu \int_{\mathbb{R}}(4|\xi|^2|\psi'|^2
+||\xi|^2\psi+\psi''|^2)\mm{d}x_3.
\end{eqnarray}
 Thus we can find a $-\lambda^2$ (depending on $\xi$) by minimizing
\begin{equation}\label{0214}
-\lambda^2(\xi)=\alpha(\xi):=\inf_{\psi\in
\mathcal{A}}E(\psi).\end{equation}
In order to emphasize the dependence on $s\in (0,\infty)$ we will sometimes write
\begin{equation*}
E(\psi,s):=E(\psi)\mbox{ and } \alpha(s):=\inf_{\psi\in
\mathcal{A}}E(\psi,s)<+\infty.\end{equation*}

Before constructing the growth solutions, we shall introduce some preliminary results,
which will be used in Subsections \ref{sec:0202}--\ref{sec:0204}.
Let the critical frequency function $S(\xi)$ for the horizontal case
and the critical frequency constant $|\xi|_{\mathrm{vc}}^M$ for the vertical case be given by the following variational forms
\begin{equation}\label{0217}
S(\xi):=\sqrt{\sup_{\psi\in H^1({\mathbb{R}}),\,\psi\not\equiv 0}
\frac{g|\xi|^2\int_{\mathbb{R}}\bar{\rho}'\psi^2\mm{d}x_3/(M\xi_1)^2-
\int_{{\mathbb{R}}}|\psi'|^2\mathrm{d}x_3}{\int_{\mathbb{R}}|\psi|^2\mathrm{d}x_3}}
>0\end{equation}
for $0<|M\xi_1|/|\xi|<M_\mathrm{c}$, and
\begin{equation}\begin{aligned}\label{0216}
|\xi|_{\mathrm{vc}}^M:=\sqrt{\inf_{\psi\in \mathcal{M}_{\mathrm{vc}}
}\frac{M^2\int_{\mathbb{R}}|\psi''|^2\mm{d}x_3}{
\int_{\mathbb{R}}(g\bar{\rho}'\psi^2-M^2|\psi'|^2)\mm{d}y}}\quad\mbox{ for }|M|\in (0,M_\mathrm{c}),
\end{aligned}\end{equation}
respectively,
where $M_\mathrm{c}$ is given by \eqref{0109}, and
$$\mathcal{M}_{\mathrm{vc}}:=\Big\{\psi\in H^2(\mathbb{R})~\Bigg|~\int_{\mathbb{R}}
(g\bar{\rho}'\psi^2-M^2|\psi'|^2)\mm{d}x_3>0\Big\}. $$

We remark that it depends on the choice of $\bar{\rho}$ whether $M_\mathrm{c}$ becomes infinite or finite. More precisely,
we have the following conclusions:
\begin{pro}\label{pro:n0201}
Assume that $\bar{\rho}'\in C_0^0(\mathbb{R})$,
\begin{enumerate}
  \item[(1)] if $\int_{\mathbb{R}}\bar{\rho}'\mm{d}x>0$, then
$M_\mathrm{c}$   is infinite,
  \item[(2)]  if $\int_{\mathbb{R}}\bar{\rho}'\mm{d}x<0$, then $M_\mathrm{c}$ is finite.
\end{enumerate}
\end{pro}
\begin{pf}
(1) We first show the first assertion. Let $\mm{sup}\bar{\rho}'\subset (-l,l)$ with $l>0$, and
\begin{equation*}
\psi_n(x):=\left\{
                       \begin{array}{ll}
1,& x\in [-l,l],\\
1+(x+l)/n, & x\in (-l-n,-l),\\
1-(x-l)/{n}, & x\in (l,l+n),\\
0& x\in\!\!\!\!\!/~ (-l-n,l+n).
\end{array}
                            \right.
\end{equation*}
Then $\psi_n(x)\in H^1(\mathbb{R})$ and
  \begin{equation*}\begin{aligned}
\frac{\int_{\mathbb{R}}g\bar{\rho}'|\psi_n|^2\mm{d}y}{\int_{\mathbb{R}}|\psi'_n|^2\mm{d}y}
\geq g n\int_{\mathbb{R}}\bar{\rho}'\mm{d}x \rightarrow \infty\mbox{ as }n\rightarrow\infty,
\end{aligned}\end{equation*}
which implies $M_\mathrm{c}=\infty$. We mention that, in such a case, we can infer that the
critical frequency constant $|\xi|_{\mathrm{vc}}^M$ defined by \eqref{0216} is equal to zero.

(2) We turn to show the second assertion by contradiction. Assume that $M_\mathrm{c}$ is infinite, i.e., there exists
a sequence of functions $\{\psi_n\}_{n=1}^\infty$
such that
$$\psi_n\in H^1(\mathbb{R}),
\mbox{ and }0<\frac{\int_{\mathbb{R}}g\bar{\rho}'|\psi_n|^2\mm{d}x}
{\int_{\mathbb{R}}|\psi'_n|^2\mm{d}x}
\rightarrow \infty\mbox{ as }n\rightarrow \infty.
$$

Denote
$$\varphi_n = \frac{\psi_n}
{\sqrt{\int_{\mathbb{R}}g\bar{\rho}'|\psi_n|^2\mm{d}x}}\in H^1(\mathbb{R}),$$
then
$$0<\frac{1} {\int_{\mathbb{R}}|\varphi'_n|^2\mm{d}x}
\rightarrow \infty\mbox{ as }n\rightarrow \infty,  $$
which yields that
$$ {\int_{\mathbb{R}}|\varphi'_n|^2\mm{d}x} \rightarrow 0\mbox{ as }n\rightarrow \infty.  $$
Now, using  Newton-Leibniz's formula, we immediately get that, for any given $\varepsilon\in
(0,1)$, there exists a $N_\varepsilon>0$ (may depend on $\varepsilon$), such that,
for any $n>N_\varepsilon$ and for any $x\in (-l,l)$,
\begin{equation}\label{A010216}
|\varphi_n(x)-\varphi_n(0)|=\left|\int_{0}^x \varphi'_n(y)\mm{d}y\right|\leq \left(\int_{0}^x |\varphi'_n(y)|^2\mm{d}y\right)^{\frac{1}{2}}x^{\frac{1}{2}}\leq \sqrt{l}\varepsilon.
\end{equation}
Thus, recalling the conditions $\int_{\mathbb{R}}\bar{\rho}'\mm{d}x< 0$ and $\bar{\rho}\in C_0^0(\mathbb{R})$, and using the following relation
\begin{equation}\label{A02}
\begin{aligned}
1=&\int_{\mathbb{R}}g\bar{\rho}'|\varphi_n|^2\mm{d}x
=\int_{\mathbb{R}}g\bar{\rho}'(\varphi_n-
\varphi_n(0)+\varphi_n(0))^2
\mm{d}x\\
=&\varphi_n^2(0)\int_{\mathbb{R}}g\bar{\rho}'\mm{d}x
+2\varphi_n(0)\int_{\mathbb{R}}g\bar{\rho}'
(\varphi_n(x)-\varphi_n(0))\mm{d}x\\
&+ \int_{\mathbb{R}}g\bar{\rho}'(\varphi_n-
\varphi_n(0))^2\mm{d}x:=R(n),  \end{aligned}
\end{equation}
we infer that there exists a constant $c$, independent of $n$ and $\varepsilon$,  such that
\begin{equation}\label{A030217}|\varphi_n(0)|\leq c\mbox{ for any }n>N.\end{equation}
In fact, if this is not true, then there exists a subsequence, denoted
by $\{\varphi_{n_m}(0)\}_{m=1}^\infty$, satisfying $|\varphi_{n_m}(0)|\geq m$,
which implies that $1=R(n_m)\rightarrow\infty$. This is a contradiction.

Finally, making use of \eqref{A010216}--\eqref{A030217} and $\int_{\mathbb{R}}\bar{\rho}'\mm{d}x< 0$, we have that for any $n>N_\varepsilon$,
\begin{equation*}
R(n) \leq (2l^{\frac{3}{2}}c + l^2)g\|\bar{\rho}'\|_{L^\infty}\varepsilon.
\end{equation*}
Consequently, letting $\varepsilon\rightarrow 0$,
we immediately see that $R(n)$ can be sufficiently small
for some $n$, and this contradicts with $R(n)\equiv 1$. Hence the second assertion holds.
\hfill$\Box$
\end{pf}

Now we define a solvable domain for a growing solution
\begin{equation}\label{0215}
\mathbb{A}^{\mm{g}}=\left\{
             \begin{array}{l}
      \big\{\xi\in \mathbb{R}^2~|~
|\xi_1M|/|\xi|\in (0,M_\mm{c}),\ |\xi|<S(\xi)\big\}\cup \{\xi_1=0\}
          \setminus\{\mf{0}\}\ \hbox{ for }\bar{\mf{M}}=M\mf{e}_1,
         \\[2mm]
 \big\{\xi\in \mathbb{R}^2~|~
        |\xi|_{\mathrm{vc}}^M<|\xi|\big\} \hbox{ for }\bar{\mf{M}}={M}\mf{e}_3
        \mbox{ with }|M|\in (0,M_\mathrm{c}).
             \end{array}
           \right.
\end{equation}
 By virtue of the definition of $M_\mathrm{c}$, it is easy to verify that the above two definitions
\eqref{0217}, \eqref{0216} make sense. In addition, it is easy to see that
 $|\xi|<S(\xi)$ in \eqref{0215} for any $|\xi|>0$ with sufficiently small $\xi_1$ (the smallness of $\xi_1$ depends on $|\xi|$).
Hence, the set $ \big\{\xi\in \mathbb{R}^2~|~|\xi_1M|/|\xi|\in (0,M_\mm{c}),\ |\xi|<S(\xi)\big\}$ is not empty.
Moreover, $|\xi|_{\mathrm{vc}}^M$ and $S(\xi )$ are finite. Next we introduce some properties concerning
with the critical frequency and the solvable domain $\mathbb{A}^{\mm{g}}$.

\begin{pro}
The supremum in \eqref{0217} is achieved for each $\xi$ and $M$ with
$\xi_1 M\neq 0$. Moreover, $S(\xi)$ is continuous on $\mathbb{D}:
=\{\xi\in\mathbb{R}^2~|~0<|M\xi_1|/|\xi|<M_\mm{c}\}$ for any given $M\neq 0$, and
 \begin{equation}\label{behavior}
 S( \xi)\rightarrow +\infty\;\;\mbox{ as }\;\frac{\xi_1}{|\xi|}\to 0\;\mbox{ and }\;\xi\in \mathbb{D}.
\end{equation}
\end{pro}
\begin{pf}
 We rewrite \eqref{0217} as
 \begin{equation*}\begin{aligned}
{S^2( \xi)}=\sup_{\psi\in \mathcal{A}_{L^2}}Q(\psi),
\end{aligned}\end{equation*}
where
 \begin{equation} \label{0219}
 Q(\psi):=\frac{g|\xi|^2}{(M\xi_1)^2}\int_{\mathbb{R}}\bar{\rho}'\psi^2\mm{d}x_3-
\int_{{\mathbb{R}}}|\psi'|^2\mathrm{d}x_3,\quad
 \mathcal{A}_{L^2}:=\left\{ H^1(\mathbb{R})~\bigg|~
\|\psi\|^2_{L^2(\mathbb{R})}=1\right\}. \end{equation}
Then, it is easy to see that ${S}^2(\xi )>0$, since $0<|M\xi_1|/|\xi|<M_\mathrm{c}$.
Let $\psi_n\in  {\mathcal{A}}_{L^2}$ be a minimizing sequence of $S^2(\xi)$,
i.e., $\limsup_{n\to\infty} Q(\psi_n)=\sup_{\psi\in\mathcal{A}_{L^2}}Q(\psi)$, we have from \eqref{0219}
that $\psi_n$ is uniformly bounded in $H^1(\mathbb{R})$ on $n$.
Hence there exists a $\psi_0\in H^1(\mathbb{R})$,
such that $\psi_n\to\psi_0$ weakly in $H^1(\mathbb{R})$
and strongly in $L^2_{\mathrm{loc}}(\mathbb{R})$. Hence,
 \begin{equation*}\begin{aligned}
0<\sup_{\psi\in\mathcal{A}_{L^2}}Q(\psi)=\limsup_{n\rightarrow \infty} Q(\psi_n)\leq Q(\psi_0)\mbox{ and }0<\|\psi_0
\|_{L^2(\mathbb{R})}\leq \liminf_{n\rightarrow \infty} \|\psi_n\|_{L^2(\mathbb{R})}=1.
\end{aligned}\end{equation*}
We proceed to verify that $\|\psi_0\|_{L^2(\mathbb{R})}=1$.
Suppose by contradiction that $\|\psi_0\|_{L^2(\mathbb{R})}<1$, then we
may scale up $\psi_0$ by $\alpha>1$ so that $\alpha \psi_0\in{\mathcal{A}}_{L^2}$.  From
this we deduce that
\begin{eqnarray*}
\sup\limits_{\psi\in{\mathcal{A}}_{L^2}}Q(\psi) \geq Q(\alpha \psi_0)=\alpha^2Q(\psi_0)\geq
\alpha^2\sup\limits_{\psi\in{\mathcal{A}}_{L^2}}Q(\psi)>\sup\limits_{\psi\in{\mathcal{A}}_{L^2}}Q(\psi),
\end{eqnarray*}
which is a contradiction. Hence
$\|\psi\|_{L^2(\mathbb{R})}=1$, which shows that $Q(\psi)$
achieves its infinimum on $\mathcal{A}_{L^2}$.

Next, we prove the continuity of ${S}^2(\xi_0)$ for each given $\xi_0$ in $\mathbb{D}$.
Letting $\xi\in \mathbb{D}\to\xi_0$, we have ${|\xi|^2}/{\xi_{1}^2}\to {|\xi_0|^2}/{\xi_{01}^2}$,
where $\xi_{01}$ represents the first component of $\xi_0$.
 Without loss of generality, we assume that ${|\xi|^2}/{\xi_{1}^2}$ and
 ${|\xi_0|^2}/{\xi_{01}^2}$ belong to a finite interval $(a,b)$ with $a>0$.
  Denote $\delta:={|\xi|^2}/{\xi_{1}^2}-{|\xi_0|^2}/{\xi_{01}^2}$,
  then $\delta\rightarrow 0$ as $\xi\rightarrow \xi_0$.

On the other hand, $Q(\psi)$
achieves its infinimum on $\mathcal{A}_{L^2}$, i.e. for any $\xi\in \mathbb{D}$, there is
$\psi_{\xi}\in \mathcal{A}_{L^2}$, such that
 \begin{equation}\begin{aligned}\label{0221}
 S^2(\xi):=\frac{g|\xi|^2}{(M\xi_1)^2}\int_{\mathbb{R}}\bar{\rho}'\psi^2_\xi\mm{d}x_3-
\int_{{\mathbb{R}}}|\psi'_\xi|^2\mathrm{d}x_3.
\end{aligned}\end{equation}
Substituting ${|\xi_0|^2}/{\xi_{01}^2}={|\xi|^2}/{\xi_{1}^2}-\delta$ into \eqref{0221}, one has
 \begin{equation}\begin{aligned}\label{0222}
{S}^2(\xi)=&\frac{g|\xi_0|^2}{(M\xi_{01})^2}\int_{\mathbb{R}}\bar{\rho}'\psi^2_{\xi}\mm{d}x_3-
\int_{{\mathbb{R}}}|\psi'_{\xi}|^2\mathrm{d}x_3+ \frac{g\delta}{M^2}
\int_{\mathbb{R}}\bar{\rho}'\psi^2_{\xi}\mm{d}x_3\\
\leq &{S}(\xi_0)+ \frac{g\delta}{M^2}\int_{\mathbb{R}}\bar{\rho}'\psi^2_{\xi}\mm{d}x_3.
\end{aligned}\end{equation}
Similarly to \eqref{0222}, we obtain
 \begin{equation*}\begin{aligned}
{S}^2(\xi_0)=&\frac{g|\xi|^2}{(M\xi_{
1})^2}\int_{\mathbb{R}}\bar{\rho}'\psi^2_{\xi_0}\mm{d}x_3-
\int_{{\mathbb{R}}}|\psi'_{\xi_0}|^2\mathrm{d}x_3- \frac{g\delta}{M^2}
\int_{\mathbb{R}}\bar{\rho}'\psi^2_{\xi_0}\mm{d}x_3\\
\leq &{S}(\xi)- \frac{g\delta}{M^2}\int_{\mathbb{R}}\bar{\rho}'\psi^2_{\xi_0}\mm{d}x_3,
\end{aligned}\end{equation*}
which, together with \eqref{0222}, yields that
 \begin{equation*}\begin{aligned}
\frac{g\delta}{M^2}
\int_{\mathbb{R}}\bar{\rho}'\psi^2_{\xi_0}\mm{d}x_3 \leq{S}^2(\xi)-{S}^2(\xi_0)
\leq \frac{g\delta}{M^2}\int_{\mathbb{R}}\bar{\rho}'\psi^2_{\xi}\mm{d}x_3.
\end{aligned}\end{equation*}
Hence,
\begin{equation*}\begin{aligned}
  |{S}^2(\xi)-{S}^2(\xi_0)|\leq \frac{g\delta}{M^2}\|\bar{\rho}'\|_{L^\infty(\mathbb{R})}
  \to 0\mbox{ as }\delta\rightarrow 0,
  \end{aligned}\end{equation*}
and ${S}^2(\xi)$ is continuous at each given point $\xi_0\in \mathbb{D}$. Finally, \eqref{behavior} obviously holds by the definition
\eqref{0217}.
The completes the proof.\hfill$\Box$
\end{pf}
\begin{pro}\label{pro:0202} Let $\mathbb{A}^{\mm{g}}$ be defined by \eqref{0215}, then
\begin{enumerate}[\quad \ (1)]
 \item the set $\mathbb{A}^{\mm{g}}$ is symmetric on $x$-axis and $y$-axis
in $\mathbb{R}^2$, respectively;
 \item there exist countably infinite    lattice points
of $(L^{-1}\mathbb{Z})^2$ belongs to $\mathbb{A}^{\mm{g}}$
;
 \item the set $\mathbb{A}^{\mm{g}}$ is a nonempty open set in $\mathbb{R}^2$.
 \end{enumerate}
\end{pro}
\begin{pf}
The first two assertions obviously hold by virtue of the definition of $\mathbb{A}^{\mm{g}}$. It suffices to show
the last assertion. Here we only give the proof of the horizontal case for the reader's convenience.

Let $0\neq \xi_0\in \mathbb{A}^{\mm{g}}$, then $\xi_0\in \{\xi_1=0\}$ or
$\xi_0\in \big\{\xi\in\mathbb{R}^2~|~ |\xi_1M|/|\xi|\in (0,M_\mm{c}),\ |\xi|<S(\xi)\big\}$.
 For the first case, noting that $|\xi_1M|/|\xi|\in (0,M_\mm{c})$ and $|\xi|<S(\xi)$
hold for any $|\xi|>0$ with sufficiently small $\xi_1$              
there exists a sufficiently small disk $B_{\xi_0}^\delta:=\{\xi\in\mathbb{R}^2~|~|\xi-\xi_0|<\delta\}
\subset \mathbb{A}^{\mm{g}}$. For the latter case, noting that
$0<|\xi_{01}M|/|\xi_0<M_\mm{c}$ and $|\xi_0|<S(\xi_0)$, we use
the continuity of $|\xi_1 M|/|\xi|$, $|\xi|$ and $S(\xi)$ as $\xi\rightarrow \xi_0\neq 0$ to deduce that
there also exists a sufficiently small disk
$B_{\xi_0}^\alpha:=\{\xi\in\mathbb{R}^2~|~|\xi-\xi_0|<\alpha\}\subset \big\{\xi\in \mathbb{R}^2~|~
|\xi_1M|/|\xi|\in (0,M_\mm{c}),\ |\xi|<S(\xi)\big\}$.
Summing up the previous discussions, we immediately conclude that
the set $\mathbb{A}^{\mm{g}}$ is a nonempty open set in $\mathbb{R}^2$ by the definition of open sets.
\hfill$\Box$
\end{pf}

\begin{pro}\label{pro:0203}Let $\xi\neq\mf{0}$. Then,
\begin{itemize}
   \item if $\xi\in \mathbb{A}^{\mm{g}}$, there exists a
\begin{equation*}\psi_0\in \left\{
             \begin{array}{ll}
   H^1(\mathbb{R})&\;\hbox{ for the horizontal case,}
               \\
   \mathcal{M}_{\mathrm{vc}}&\;\hbox{ for the vertical case,}      \end{array}
           \right.
\end{equation*}
such that $E_0(\psi_0)<0$;
   \item else,
$|\xi|^2E_0(\psi)$ $\geq 0$ for any
\begin{equation*}\psi\in \left\{
             \begin{array}{ll}
      H^1(\mathbb{R})&\;\hbox{ for the horizontal case; }
               \\[1mm]
 H^2(\mathbb{R})&\;\hbox{ for the vertical case.}    \end{array}
           \right.
\end{equation*}
 \end{itemize}
\end{pro}
\begin{pf}
The above assertions in fact follow from the definitions \eqref{0109}, \eqref{0212}, \eqref{0217}, \eqref{0216} and \eqref{0215}, here we only give the proof of the horizontal case
for the reader's convenience.

Let $\xi\in \mathbb{A}^{\mm{g}}$. If $\xi_1=0$, then obviously, there exists a $\psi_0$, such that
$$E_0= \int_{\mathbb{R}}\Big[(M\xi_1)^2\Big(|\psi_0|^2+\frac{|\psi'_0|^2}{|\xi|^2}\Big)
-g\bar{\rho}'\psi^2_0\Big]\mm{d}x_3= \int_{\mathbb{R}} -g\bar{\rho}'\psi^2_0\mm{d}x_3<0.$$
If $\xi$ satisfies $|M\xi_1|/|\xi|\in (0,M_\mm{c})$ and $|\xi|<S(\xi)$,
then by virtue of the definition of \eqref{0217}, there also exists a $\psi_0$, such that
$$E_0= \int_{\mathbb{R}}\Big[(M\xi_1)^2\Big(|\psi_0|^2+\frac{|\psi'_0|^2}{|\xi|^2}\Big)
-g\bar{\rho}'\psi^2_0\Big] \mm{d}x_3<0.$$
Summing up the above discussions, we see that there exists a $\psi_0$,
such that $E_0(\psi_0)<0$ for $\xi\in \mathbb{A}^{\mm{g}}$.

Let $\xi\not\in\mathbb{A}^{\mm{g}}\cup\{\mf{0}\}$, then $\xi$ can be divided by two cases:
 (i) $|M\xi_1|/|\xi|\in (0, M_\mm{c})$ with $|\xi|\geq S(\xi)$, and
(ii) $|M\xi_1|/|\xi|\geq M_\mm{c}$ if $M_\mm{c}<\infty$
(noting that this case will not appear if $M_\mm{c}=\infty$).
For the first case, in view of the definition of \eqref{0217}, we find that
$$|\xi|^2E_0(\psi)= |\xi|^2  \int_{\mathbb{R}}\Big[(M\xi_1)^2\Big(|\psi|^2+\frac{|\psi'|^2}
{|\xi|^2}\Big)-g\bar{\rho}'\psi^2\Big]\mm{d}x_3\geq 0\;\;\mbox{ for any }\psi\in H^1(\mathbb{R}). $$
Finally, for the second case, by the definition of \eqref{0109}, we have
$$  \int_{\mathbb{R}}\Big(\frac{(M\xi_1)^2}{|\xi|^2}|\psi'|^2 -
g\bar{\rho}'\psi^2\Big)\mm{d}x_3\geq 0\;\;\mbox{ for any }\psi\in H^1(\mathbb{R}),$$
which yields
$$|\xi|^2E_0= |\xi|^2  \int_{\mathbb{R}}\Big[(M\xi_1)^2\Big(|\psi|^2+\frac{|\psi'|^2}{|\xi|^2}\Big)
-g\bar{\rho}'\psi^2\Big]\mm{d}x_3\geq 0\;\;\mbox{ for any }\psi \in H^1(\mathbb{R}).$$
Summarizing the above discussions, we see that $|\xi|^2E_0(\psi)$ $\geq 0$ for any $\psi\in H^1(\mathbb{R})$
if $\xi \in\!\!\! \! \! /~ \mathbb{A}^{\mm{g}}\cup\{\mf{0}\}$.
This completes the proof for the horizontal case. \hfill $\Box$
\end{pf}
\subsection{Solutions to the variational problem}\label{sec:0202}
In this seusection we show that a minimizer
of \eqref{0214} exists for the case of
$\inf_{\psi\in \mathcal{A}}E(\psi,s)<0$ which will be shown to be true for sufficiently small $s$
in Proposition \ref{pro:0205} below, and that the corresponding
Euler-Lagrange equations are equivalent to \eqref{0208}, \eqref{0209}.

\begin{pro}\label{pro:0204}
 For any fixed $s>0$ and $\xi$ with $|\xi|\neq 0$, we have
 \begin{enumerate}
   \item[(1)] $\inf_{\psi\in\mathcal{A}}E(\psi,s)>-\infty$.
   \item[(2)] if there exists a
$\bar{\psi}\in \mathcal{A}$, such that $E(\bar{\psi})<0$, then $E({\psi})$
achieves its infinimum on $\mathcal{A}$.
   \item[(3)] let $\tilde{\psi}$ be
a minimizer and $-\lambda^2:=E(\tilde{\psi})$, then the pair $(\tilde{\psi}$,
$\lambda^2)$ satisfies \eqref{0208}, \eqref{0209}. Moreover,
$\tilde{\psi}\in H^\infty(\mathbb{R}):=\cap_{k=0}^\infty H^k(\mathbb{R})$.
 \end{enumerate}
\end{pro}
\begin{pf} We only show the proposition for the horizontal case, the vertical case
can be dealt with in the same manner.

 (1) Noticing that for any $\psi\in \mathcal{A}$,
\begin{equation}\label{0226}E(\psi)\geq -{g|\xi|^2}\int_{\mathbb{R}}
\bar{\rho}'\psi^2\mathrm{d}x_3\geq
-{g}\left\|\frac{\bar{\rho}'}{\bar{\rho}}\right\|_{L^\infty(\mathbb{R})}\int_{\mathbb{R}}
\bar{\rho}|\xi|^2\psi^2\mathrm{d}x_3\geq
-{g}\left\|\frac{\bar{\rho}'}{\bar{\rho}}\right\|_{L^\infty(\mathbb{R})},
\end{equation}
we see that $E$ is bounded from below on $\mathcal{A}$ by virtue of \eqref{0102}.
This proves (1).

(2) We proceed to show (2). Let $\psi_n\in\mathcal{A}$ be a minimizing sequence,
then $E(\psi _n)$ is bounded.
This together with \eqref{0210} and \eqref{0226} implies that
$\psi_n$ is bounded in $H^2(\mathbb{R})$. So, there exists a $\psi_0\in
H^2(\mathbb{R})$, such that $\psi_n\rightarrow \psi_0$ weakly in $H^2(\mathbb{R})$
and strongly in $H^1_{\mathrm{loc}}(\mathbb{R})$.
Moreover, by the lower semi-continuity, \eqref{0102} and the assumption that $E(\bar{\psi})<0$
for some $\bar{\psi}\in \mathcal{A}$, we deduce that
\begin{equation*}E(\psi_0)\leq \liminf_{n\rightarrow
\infty}E(\psi_n)=\inf_{\psi\in \mathcal{A}}E(\psi)<0,\quad\mbox{and }\; 0<J(\psi_0)\leq 1.
\end{equation*}

Suppose by contradiction that $J(\psi_0)<1$. By the homogeneity of $J$
we may find an $\alpha>1$ so that $J(\alpha \psi_0)=1$, i.e., we may scale up
$\psi_0$ so that $\alpha \psi_0\in{\mathcal{A}}$. From this we deduce that
\begin{eqnarray*}
E(\alpha\psi_0)=\alpha^2E(\psi_0)\leq\alpha^2\inf\limits_{\mathcal{A}}
E<\inf\limits_{\mathcal{A}} E<0,
\end{eqnarray*}
which is a contradiction since $\alpha\psi_0\in{\mathcal{A}}$. Hence
$J(\psi_0)=1$ and $\psi_0\in{ \mathcal{A}}$. This shows that $E(\psi)$
achieves its infinimum on $\mathcal{A}$.

(3) Finally we prove (3). Notice that since $E(\psi)$ and $J(\psi)$ are homogeneous
of degree $2$, \eqref{0214} is equivalent to
  \begin{equation}\label{0227}\alpha(s)=\inf_{\psi\in
  H^2(\mathbb{R}),\psi \equiv\!\!\!\!/~0}\frac{E(\psi)}{J(\psi)}.
\end{equation}
For any $\tau\in \mathbb{R}$ and $\psi\in H^2(\mathbb{R})$
we take $\psi(\tau):=\tilde{\psi}+\tau\psi$, then \eqref{0227} implies
  \begin{equation*}E(\psi(\tau))+\lambda^2J(\psi(\tau))\geq 0.
\end{equation*}
If we set $I(\tau)=E(\psi(\tau))+\lambda^2J(\psi(\tau))$, then
$I(\tau)\geq 0$ for all $\tau\in \mathbb{R}$ and $I(0)=0$. This
implies $I'(0)=0$. By virtue of \eqref{0210} and \eqref{0211}, a
direct computation leads to
  \begin{equation}\label{0228}\begin{aligned}
  & s\mu \int_{\mathbb{R}}\Big(4|\xi|^2\tilde{\psi}'\psi'+(|\xi|^2\tilde{\psi}+\tilde{\psi}'')
  (|\xi|^2\psi +\psi'') +(M\xi_1)^2(|\xi|^2\tilde{\psi}\psi+\tilde{\psi}'\psi')\Big)\mathrm{d}x_3\\
 & =   g|\xi|^2\int_\mathbb{R}\bar{\rho}'\tilde{\psi}\psi\mathrm{d}x_3
  -\lambda^2\int_{\mathbb{R}}\bar{\rho}(|\xi|^2\tilde{\psi}\psi+\tilde{\psi}'\psi')\mathrm{d}x_3,
\end{aligned}\end{equation}
where we have used the upper boundedness of $\bar{\rho}$.

 By further assuming that $\psi$ is compactly supported in $\mathbb{R}$, we
find that $\tilde{\psi}$ satisfies the equation \eqref{0209} in the weak sense
on $\mathbb{R}$ for the horizontal case. In order to improve the regularity of $\tilde{\psi}$,
we rewrite \eqref{0228} as
  \begin{equation}\label{0229}\begin{aligned}
\int_\mathbb{R}
\tilde{\psi}''\psi''\mathrm{d}x_3
=&\frac{1}{s\mu}\int_{\mathbb{R}}\left(g|\xi|^2\bar{\rho}'\tilde{\psi}
-\lambda^2(|\xi|^2\bar{\rho}\tilde{\psi}-(\bar{\rho}\tilde{\psi}')')\right.\\
&  \left. +s\mu(2|\xi|^2\tilde{\psi}''
-|\xi|^4\tilde{\psi})+(M\xi_1)^2(\tilde{\psi}''-|\xi|^2\tilde{\psi} )\right)\psi\mathrm{d}x_3\\
:=&\int_{\mathbb{R}}f \psi\mathrm{d}x_3.
\end{aligned}\end{equation}

For any $n\geq 1$, let $\psi_{1,n}$, $\psi_2\in C_0^\infty(\mathbb{R})$
satisfy $\psi_{1,n}(x_3)\equiv 1$ for $|x_3|\leq n$. If one takes
$\psi=\psi_{1,n}\int_{-\infty}^{x_3}\psi_2\mathrm{d}y$ in \eqref{0229}, then one has
  \begin{equation*}\begin{aligned}\int_\mathbb{R}
(\psi_{1,n}\tilde{\psi}'')\psi_2'\mathrm{d}x_3=&
\int_{\mathbb{R}}\left(f\psi_{1,n}\int_{-\infty}^{x_3}\psi_2\mathrm{d}y
-\psi_{1,n}''\psi''\int_{-\infty}^{x_3}\psi_2\mathrm{d}y-2\psi'_{1,n}\tilde{\psi}''
\psi_2\right)\mathrm{d}x_3 \\
= & \int_{\mathbb{R}}\left(\int_{x_3}^{+\infty}(f\psi_{1,n}-\psi_{1,n}''\tilde{\psi}'')
\mathrm{d}y-2\psi'_{1,n}\tilde{\psi}'' \right)\psi_2\mathrm{d}x_3,
\end{aligned}\end{equation*}
which, recalling $\tilde{\psi}\in H^2(\mathbb{R})$, implies $\tilde{\psi}''\in H^1_{\mathrm{loc}}(\mathbb{R})$ and
\begin{equation*}\begin{aligned}
\tilde{\psi}'''=(\psi_{1,n}\psi'')'=\int_{x_3}^{+\infty}(f\psi_{1,n}-\psi_{1,n}''\tilde{\psi}'')\mathrm{d}y\quad
 \mbox{ for any }x_3 \mbox{ with }|x_3|\leq n.
\end{aligned}\end{equation*}
Integrating by parts, we can rewrite \eqref{0229} as
  \begin{equation*}\begin{aligned}
 -\int_\mathbb{R}
\tilde{\psi}'''\psi'\mathrm{d}x_3
=&\frac{1}{s\mu }\int_{\mathbb{R}}\left(g|\xi|^2\bar{\rho}'\tilde{\psi} -
\lambda^2(|\xi|^2\bar{\rho}\tilde{\psi}-(\bar{\rho}\tilde{\psi}')')\right.\\ &\left.\qquad
 \qquad  +s\mu(2|\xi|^2\tilde{\psi}''
-|\xi|^4\tilde{\psi})+(M\xi_1)^2(\tilde{\psi}''-|\xi|^2\tilde{\psi} ) \right)\psi\mathrm{d}x_3,
\end{aligned}\end{equation*}
which, keeping in mind that $\tilde{\psi}\in H^2(\mathbb{R})$, yields
$\tilde{\psi}''''\in L^2(\mathbb{R})$. Hence $\tilde{\psi}\in H^4_{\mathrm{loc}}(\mathbb{R})\cap
C^{3,1/2}_{\mathrm{loc}}(\mathbb{R})$, and
$\tilde{\psi}^{'}(\infty)=\tilde{\psi}^{''}(\infty)=\tilde{\psi}^{'''}(\infty)=0$. Using
these facts, H\"{o}lder's inequality, and integration by
parts, we conclude that
  \begin{equation*}\begin{aligned}
\|\tilde{\psi}'''\|_{L^2(\mathbb{R})}^2=\int_{\mathbb{R}}|\tilde{\psi}'''|^2\mathrm{d}x_3
=-\int_{\mathbb{R}}\tilde{\psi}''\tilde{\psi}''''\mathrm{d}x_3\leq\|\tilde{\psi}''\|_{L^2(\mathbb{R})}\|
\tilde{\psi}''''\|_{L^2(\mathbb{R})},
\end{aligned}\end{equation*}
i.e., $\tilde{\psi}'''\in L^2(\mathbb{R})$. Consequently, $\tilde{\psi}\in H^4(\mathbb{R})$
solves \eqref{0208}, \eqref{0209}. This immediately gives
$\tilde{\psi}\in H^\infty(\mathbb{R})$. \hfill $\Box$
\end{pf}

Next, we want to show that there is a fixed point such that
$\lambda=s$. To this end, we first give some properties of
$\alpha(s)$ as a function of $s> 0$.
\begin{pro}\label{pro:0205}
The function $\alpha(s)$ defined on $(0,\infty)$ enjoys the following properties:
\begin{enumerate}[\quad \ (1)]
 \item $\alpha(s)\in C_{\mathrm{loc}}^{0,1}(0,\infty)$ is nondecreasing.
    \item  for any $\xi\in \mathbb{A}^{\mm{g}}$,
 there exist constants $c_1$, $c_2>0$ depending on  $g$, $M$, $\bar{\rho}$,
  $\mu$, and $\xi$, such that
  \begin{equation}\label{0231}\alpha(s)\leq -c_1+sc_2
  .\end{equation}
  \end{enumerate}
\end{pro}
\begin{pf} We still give the proof for the horizontal case only, and
the vertical case can be dealt with in the same way.

(1) Let  $\{\psi^n_{s_2}\}\subset\mathcal{A}$ be  a minimizing sequence of
$\inf_{\psi\in \mathcal{A}}E(\psi,s_2)$.   From \eqref{0210} and \eqref{n03010212} it follows that
\begin{equation*}
\alpha(s_1)\leq \limsup_{n\rightarrow\infty}E(\psi_{s_2}^n,s_1)
\leq \limsup_{n\rightarrow\infty}E(\psi_{s_2}^n,s_2)=\alpha(s_2)\;
\mbox{ for any }0<s_1<s_2<\infty.
\end{equation*}
Hence $\alpha(s)$ is nondecreasing on $(0,\infty)$. Next we shall use this fact to
show the continuity of $\alpha(s)$.

Let $I:=[a,b]\subset\mathbb{R}^+$ be a bounded interval.
In view of \eqref{0226} and the monotonicity of $\alpha(s)$, we know that
\begin{equation} \label{02321} |\alpha(s)|\leq \max\left\{|\alpha(b)|,{g}
 \left\| {\bar{\rho}'}/{\bar{\rho}}
\right\|_{L^\infty(\mathbb{R})}\right\}<\infty\;\;\mbox{ for any }s\in I.
\end{equation}
On the other hand, for any $s\in I$, there exists a minimizing sequence
$\{\psi^n_{s}\}\subset \mathcal{A}$ of $\inf_{\psi\in \mathcal{A}}E(\psi,s)$, such that
$$ |\alpha(s)-E(\psi_{s}^n,s)|<1. $$
  Making use of \eqref{0210}--\eqref{n03010212}
and  \eqref{02321}, we infer that
\begin{equation*}\begin{aligned} 0\leq E_1(\psi_s^n,s)=&\frac{E(\psi_s^n,s)}{s}
+ \frac{g|\xi|^2}{s}\int_{\mathbb{R}} \bar{\rho}'\psi^2\mathrm{d}x_3-\frac{(M\xi_1)^2}{s}
\int_{\mathbb{R}}\left({|\xi|^2}|\psi|^2+ {|\psi'|^2}\right)\mathrm{d}x_3  \\
\leq & \frac{1+\max\{|a(b)|,{g}\left\|{\bar{\rho}'}/{\bar{\rho}}
\right\|_{L^\infty(\mathbb{R})}\}}{a}+\frac{g}{a}\left\|\frac{\bar{\rho}'}{\bar{\rho}}
\right\|_{L^\infty(\mathbb{R})}:=K.
\end{aligned}\end{equation*}

For $s_i\in I$ ($i=1,2$), we further find that
\begin{equation}\begin{aligned}\label{0235}\alpha(s_1)\leq \limsup_{n\rightarrow
\infty}E(\varphi_{s_2}^n,s_1)\leq & \limsup_{n\rightarrow
\infty}E(\psi_{s_2}^n,s_2)+|s_1-s_2|\limsup_{n\rightarrow
\infty}E_1(\psi_{s_2}^n)\\
\leq & \alpha(s_2)+K|s_1-s_2|.
\end{aligned}\end{equation}
Reversing the role of the indices $1$ and $2$ in the
derivation of the inequality \eqref{0235}, we obtain the same boundedness with
the indices switched. Therefore, we have
\begin{equation*}\begin{aligned}|\alpha(s_1)-\alpha(s_2)|\leq K|s_1-s_2|,
\end{aligned}\end{equation*}
which yields $\alpha(s)\in C_{\mathrm{loc}}^{0,1}(0,\infty)$.

(2) Since $\xi\in \mathbb{A}^{\mm{g}}$, by virtue of Proposition \ref{pro:0203}, there exists
a $\bar{\psi}\in H^1(\mathbb{R})$, such that
\begin{equation}\label{0236}
E_0(\psi)=     \int_{\mathbb{R}}(M\xi_1)^2\left(|\bar{\psi}|^2
+\frac{|\bar{\psi}'|^2}{|\xi|^2}\right)-g\bar{\rho}'\bar{\psi}^2\mm{d}x_3<0.
\end{equation}
On the other hand, $H^2(\mathbb{R})$ is dense in $H^1(\mathbb{R})$,
thus there is a function sequence $\bar{\psi}_n\in H^2(\mathbb{R})$, so that
\begin{equation}\label{0237}
\bar{\psi}_n\rightarrow\bar{\psi}\;\;\mbox{ stronly in }H^1(\mathbb{R}).
\end{equation}

Putting \eqref{0236} and \eqref{0237} together, we see that there is a subsequence
$\bar{\psi}_{n_0}\in H^2(\mathbb{R})$, such that $\bar{\psi}_{n_0}\not\equiv 0$ and
$$ E_0(\bar{\psi}_{n_0})=\int_{\mathbb{R}}\Big[(M\xi_1)^2\Big( |\bar{\psi}_{n_0}|^2
+\frac{|\bar{\psi}_{n_0}'|^2}{|\xi|^2}\Big) -g\bar{\rho}'\bar{\psi}_{n_0}^2\Big]\mm{d}x_3<0.
$$
Thus, we have
$$ \alpha(s)= \inf_{\psi\in \mathcal{A}}E(\psi)=\inf_{\psi\in
 H^2(\mathbb{R})}\frac{E(\psi)}{J(\psi)} \leq
\frac{E(\bar{\psi}_{n_0})}{J(\bar{\psi}_{n_0})}=|\xi|^2\frac{E_0(\bar{\psi}_{n_0})}{J(\bar{\psi}_{n_0})}
+s\frac{E_1(\bar{\psi}_{n_0})}{J(\bar{\psi}_{n_0})}:= -c_1+sc_2
$$
for two positive constants $c_1: =c_1(g,M,\bar{\rho},|\xi|)$ and $c_2:=c_2(g,M,\mu,\bar{\rho},|\xi|)$.
This completes the proof for the vertical case.  \hfill $\Box$
\end{pf}

Given $\xi\in \mathbb{A}^g$, by virtue of
\eqref{0231}, there exists a $s_0>0$ depending on the quantities $g$, $M$,
$\mu$, $\bar{\rho}$ and  $\xi$, such that for any $s\in (0, s_0]$, $\alpha(s)<0$. Let
\begin{equation}\label{0240}\mathfrak{S}_{\xi}:=\sup\{s~|~\alpha(\tau)<0\mbox{ for any }\tau\in
(0,s)\},\end{equation}  then $\mathfrak{S}_{\xi}>0$. This allows
us to define $\lambda(s)=\sqrt{-\alpha(s)}>0$ for any $s\in
\mathcal{S}_{\xi}:=(0,\mathfrak{S}_{\xi})$. Therefore, as a
result of Proposition \ref{pro:0204}, we have the following existence for the modified problem
\eqref{0208}, \eqref{0209}.
\begin{pro}\label{pro:0206}
For each $\xi\in \mathbb{A}^g$ and $s\in \mathcal{S}_{\xi}$ there is
a solution $\psi=\psi(\xi,x_3)\not\equiv 0$ with
$\lambda=\lambda(\xi,s)>0$ to the problem \eqref{0208}, \eqref{0209}.
Moreover, $\psi\in H^\infty(\mathbb{R})\cap \mathcal{A}$.
\end{pro}

 Now, we can use Proposition \ref{pro:0205}, \eqref{0226} and \eqref{0240}
to find that $\lambda(s)\in C_{\mathrm{loc}}^{0,1}(\mathcal{S}_{\xi})$ is nonincreasing,
$\lambda(s)\leq \sqrt{g}\|\sqrt{{\bar{\rho}'}/{\bar{\rho}}}\|_{L^\infty(\mathbb{R})}$, $\lim_{s\rightarrow 0}\lambda(s)>0$
and $\lim_{s\rightarrow \mathfrak{S}_{\xi}}\lambda(s)=0$
if $\mathfrak{S}_{|\xi|}<+\infty$. Hence, we can employ a fixed-point argument to find
$s\in\mathcal{S}_{\xi}$ so that $s=\lambda(\xi, s)$, and thus obtain a solution to
the original problem \eqref{0207}, \eqref{0208}.
\begin{pro}\label{pro:0207} Let $\xi\in \mathbb{A}^g$,
then there exists a unique $s\in \mathcal{S}_{\xi}$, such that
$\lambda(\xi,s)=\sqrt{-\alpha(s)}>0$ and $s=\lambda(\xi,s)$.
\end{pro}
\begin{pf}
We refer to \cite[Theorem 3.8]{GYTI2} (or \cite[Lemma 3.7]{WYC}) for
a proof.\hfill $\Box$
\end{pf}
Moreover, in view of Propositions \ref{pro:0206} and \ref{pro:0207},
 we conclude the following existence for the problem \eqref{0207}--\eqref{0208}.
\begin{thm}\label{thm:0201}
For each $\xi\in \mathbb{A}^g$, there exist
$\psi=\psi(\xi,x_3)\equiv\!\!\!\!\!\!/\ 0$ and $\lambda(\xi)>0$
satisfying \eqref{0207} and \eqref{0208}. Moreover, $\psi\in H^\infty(\mathbb{R})\cap \mathcal{A}$.
\end{thm}

We end this subsection by giving  additional properties of the solutions
established in Theorem \ref{thm:0201} in terms of $\lambda(\xi)$,
which show that $\lambda$ is a bounded, continuous function of $\xi$.
\begin{pro}\label{pro:0208}
The positive function
$\lambda :\mathbb{A}^\mm{g}\rightarrow \mathbb{R}^+$ is continuous and satisfies
\begin{equation}\label{0241}
\sup_{\xi\in \mathbb{A}^{\mm{g}}}\lambda(\xi)\leq
\sqrt{g\left\|{\bar{\rho}'}/{\bar{\rho}}\right
\|_{L^\infty(\mathbb{R})}}.\end{equation}
\end{pro}
\begin{pf} The boundedness of $\lambda$ in \eqref{0241} follows from \eqref{0226}.
As for the proof of the continuity of $\lambda$, we still give the proof
for the horizontal case only.

First, let $\xi_0\in \mathcal{A}$ be arbitrary but fixed, and $\xi\rightarrow \xi_0$.
Without loss of generality, assume $\xi\in \mathbb{A}^\mm{g}$, since $\mathbb{A}^\mm{g}$
is an open set by Proposition \ref{pro:0202}. Then there exists an interval $[a,b]\subset \mathbb{R}^+$
so that $|\xi|$ and $|\xi_0|\in (a,b)$. Let
 $\delta=|\xi|^2-|\xi_0|^2$, then
$\delta\rightarrow 0 $ as $|\xi|\rightarrow |\xi_0|$.

(i) We begin with the proof of  the following conclusion:
\begin{equation}\label{0242}
\lim_{\xi\rightarrow \xi_0}
\alpha(\xi,s)=\alpha(\xi_0,s)\;\;\mbox{ for any }s\in
\mathcal{S}_{\xi}.\end{equation}

By virtue of Proposition \ref{pro:0206}, for any $\xi\in \mathbb{A}^{\mm{g}}$,
there exists a function $\psi_{\xi}\in \mathcal{A}$, such that
\begin{equation}\label{0243}
\alpha(\xi)=\int_{\mathbb{R}}\Big[s\mu(4|\xi|^2|\psi'_{\xi}|^2
+||\xi|^2\psi_{\xi}+\psi''_{\xi}|^2) +(M\xi_1)^2\left(|\xi|^2|\psi_{\xi}|^2
+{|\psi'_{\xi}|^2}\right) - g|\xi|^2\bar{\rho}'\psi^2_{\xi}\Big] \mathrm{d}x_3<0,
\end{equation}
which, together with \eqref{0211}, yields
\begin{equation}\label{0244}\|\psi_{\xi}\|_{H^2(\mathbb{R})}\leq  c_3,\end{equation}
where $c_3$ depends on $g$, $M$, $\mu$, $\bar{\rho}$, $a$, $b$ and $s$.

To deal with the term involved with $|\xi_1|$, we  denote $\delta_1:=|\xi_1|^2-|\xi_{01}|^2$.
 Substitution of $|\xi|^2=|\xi_0|^2+\delta$ and $|\xi_1|^2=|\xi_{01}|^2+\delta_1$ into \eqref{0243} results in
\begin{equation}\label{0245}\begin{aligned}
\alpha(\xi)=&\int_{\mathbb{R}}\Big[ s\mu(4|\xi_0|^2|\psi'_{\xi}|^2
+||\xi_0|^2\psi_{\xi}+\psi''_{\xi}|^2) +(M\xi_{01})^2\left(|\xi_0|^2|\psi_{\xi}|^2
+{|\psi'_{\xi}|^2}\right) - g|\xi_0|^2\bar{\rho}'\psi^2_{\xi}\Big] \mathrm{d}x_3\\
& +\delta f(\delta,\psi_{\xi}) +\delta_1 h(\psi_{\xi})\geq\alpha(\xi_0)
+\delta f(\delta,\psi_{\xi}) +\delta_1 h(\psi_{\xi}),\end{aligned}\end{equation}
where
$$ f(\delta,\psi_{\xi})=\int_{\mathbb{R}}\Big[ s\mu(4|\psi'_{\xi}|^2 +2\delta
\psi_{\xi}(|\xi_0|^2\psi_{\xi}+\psi''_{\xi})+\delta\psi^2_{\xi}+(M\xi_{01})^2|\psi_{\xi}|^2-
g\bar{\rho}'\psi^2_{\xi}\Big] \mathrm{d}x_3 $$
and
$$ h(\delta,\psi_{\xi})=
\int_{\mathbb{R}}M^2\left((|\xi_0|^2+\delta)|\psi_{\xi}|^2+{|\psi'_{\xi}|^2}\right)\mathrm{d}x_3. $$
By H\"{o}lder's inequality and \eqref{0244}, we can bound
\begin{equation}\label{0246}\begin{aligned}|
f(\delta,\psi_{\xi})|+|h(\delta,\psi_{\xi})|\leq c_4\;\;\mbox{ for some constant }c_4.
\end{aligned}\end{equation}

Similarly to \eqref{0245} and \eqref{0246}, we also have
\begin{equation}\label{0247}\begin{aligned}\alpha(\xi_0)\geq\alpha(\xi)-\delta
f(-\delta,\psi_{\xi_0})-\delta_1
h(\psi_{\xi_0})\end{aligned}\end{equation}
{ and }\begin{equation}\label{0248}\begin{aligned}|f(-\delta,\psi_{\xi_0})|+|
h(-\delta,\psi_{\xi_0})|\leq
c_5.\end{aligned}\end{equation}

Combining \eqref{0245} with \eqref{0247}, we get
$$ \delta f(-\delta,\psi_{\xi_0}) +\delta_1 h(\psi_{\xi_0})
\geq \alpha(\xi)-\alpha(\xi_0)\geq \delta f(\delta,\psi_{\xi})+\delta_1 h(\psi_{\xi}), $$
which, together with \eqref{0246} and \eqref{0248}, implies  \eqref{0242}. Hence,
\begin{equation}\label{0249nn}
\lim_{\xi\rightarrow \xi_0}
\lambda(\xi,s)=\lambda(\xi_0,s)\;\;\mbox{ for any }s\in
\mathcal{S}_{\xi},\end{equation}
because of $\lambda(\xi,s)=\sqrt{-\alpha(\xi,s)}$.

(ii) In view of \eqref{0249nn} and Proposition \ref{pro:0207}, we see that
for any $\varepsilon>0$, there exists a $\delta>0$ such that
 $|\lambda(\xi,s_{\xi_0})-
 \lambda(\xi_0,s_{\xi_0})|<\varepsilon$  and
 $s_{\xi_0}=\lambda({\xi_0},s_{\xi_0})=\sqrt{-\alpha({\xi_0,s_{\xi_0}})}$ for any
$\xi\in B_{\xi_0}^\delta:=\{\mf{y}\in \mathbb{R}^2~|~|\mf{y}-\xi_0|<\delta\}\subset\mathbb{A}^\mm{g}$.
 On the other hand, for each $|\xi|>0 $, $\lambda(s)$ is nonincreasing and continuous
 on $\mathcal{S}_{\xi}$, and there exists a unique
 $s_{\xi}\in \mathcal{S}_{\xi}$ satisfying
$\lambda(\xi,s_{\xi})=s_{\xi}>0$ by Proposition \ref{pro:0207}.
Consequently, we immediately infer that
$|\lambda(\xi,s_{\xi})-\lambda(\xi_0,s_{\xi_0})|<\varepsilon$
with  $s_{\xi}=\lambda(\xi,s_{\xi})$. Therefore, $\lambda(\xi)$
is continuous.  This completes the proof of the proposition.  \hfill $\Box$
\end{pf}

\subsection{Construction of a solution to  the system  \eqref{0205}, \eqref{0206}}
A solution to \eqref{0207}, \eqref{0208} gives rise to a solution of
the ODEs \eqref{0205}, \eqref{0206} for the growing mode velocity $\mathbf{u}$ as well.
\begin{thm}\label{thm:0302}
For each $\xi\in \mathbb{A}^\mm{g}$, there exists a
solution $({\varphi},{\theta},{\psi},{\pi}):= ({\varphi},{\theta},{\psi},{\pi})(\xi,x_3)$
with $\lambda =\lambda(\xi)>0$ to \eqref{0205}, \eqref{0206}, and
the solution belongs to $H^\infty(\mathbb{R})$. Moreover $\psi\in \mathcal{A}$.
\end{thm}
\begin{pf}
We still give the proof for the horizonal case only. First, in view of Theorem \ref{thm:0201},
we have a solution $(\psi ,\lambda ):=(\psi (\xi,x_3),\lambda (\xi))$
satisfying \eqref{0207}, \eqref{0208}. Moreover, $\lambda >0$ and
$\psi\in\mathcal{A}\cap H^\infty(\mathbb{R}^2)$.
Then, multiplying \eqref{0205}$_1$ and \eqref{0205}$_2$ by $\xi_1$
and $\xi_2$ respectively, adding the resulting equations, and
utilizing \eqref{0205}$_4$, we find that $\pi$ can be expressed by $\psi$, i.e.,
\begin{eqnarray}
&&\pi =\pi(\xi,x_3)\nonumber\\
&&\quad =[-{\lambda\mu (\xi_1\varphi+\xi_2\theta)''+(\lambda^2\bar{\rho}+\lambda\mu
|\xi|^2+M^2\xi_1^2)(\xi_1\varphi+\xi_2\theta)}-M^2\xi_1|\xi|^2\varphi]/({\lambda|\xi|^{2}})\nonumber\\
&&\label{0249} \quad =[{\lambda\mu \psi'''-(\lambda^2\bar{\rho}+\lambda\mu
|\xi|^2+M^2\xi_1^2)\psi'}-M^2\xi_1|\xi|^2\varphi]/({\lambda|\xi|^{2}}).
 \end{eqnarray}
Thus \eqref{0205}$_1$ can be rewritten as
\begin{equation}\label{0250}-\varphi''+\sigma \varphi=\omega \end{equation}
with boundary conditions
\begin{equation}\label{0251}
 \varphi(-\infty)=\varphi(+\infty)=0,
\end{equation}
where $\sigma=(\lambda^2\bar{\rho}+\lambda\mu |\xi|^2+M^2\xi_1^2)/{(\lambda\mu)}>0$
 and  $\omega=\xi_1[{\lambda\mu \psi'''-(\lambda^2\bar{\rho}+\lambda\mu
|\xi|^2+M^2\xi_1^2)\psi'}]/({\lambda\mu|\xi|^{2}})$ $\in H^\infty(\mathbb{R})$.
By the ODE theory on a bounded interval and the domain expansion technique,
we can obtain a unique solution $\varphi:=\varphi(\xi,x_3)\in H^\infty(\mathbb{R})$ to
\eqref{0250}, \eqref{0251}. In view of \eqref{0249}, thus $\pi$ can be uniquely determined by  the known functions $\psi$ and $\varphi$. Employing arguments similar to those in the construction of $\varphi$,
we can obtain a unique solution $\theta:=\theta(\xi,x_3)\in H^\infty(\mathbb{R})$ to
\eqref{0205}$_3$ with $\theta(-\infty)=\theta(+\infty)=0$, where the solution $\theta$
depends on the known functions $\pi$ and $\varphi$.

Finally, by a simple computation, we can check that $({\varphi},{\theta},{\psi},{\pi})$
with $\lambda$ constructed as above indeed is a solution to \eqref{0205}, \eqref{0206}.
 \hfill
$\Box$
\end{pf}
\begin{rem}\label{rem:0201}
For each $x_3$, it is easy to see that the solution
$({\varphi}(\xi,\cdot),{\theta}(\xi,\cdot),{\psi}(\xi,\cdot),{\pi}(\xi,\cdot),\lambda(\xi))$
constructed in Theorem \ref{thm:0302} possesses the following properties:
\begin{enumerate}[\quad \ (1)]
  \item $\lambda(\xi)$, ${\psi}(\xi,\cdot)$ and ${\pi}(\xi,\cdot)$
  are even on $\xi_1$ or $\xi_2$, when another variable is fixed;
  \item ${\varphi}(\xi,\cdot)$ is odd on $\xi_1$, but even on $\xi_2$,
  when another variable is fixed;
  \item ${\theta}(\xi,\cdot)$ is even on $\xi_1$, but odd on $\xi_2$,
  when another variable is fixed.
  \end{enumerate}
\end{rem}
\begin{rem}
We mention that the system  \eqref{0205}, \eqref{0206} for the vertical case enjoys
rotational structure.  Hence, we can also
use the rotation method as in \cite{NJTSC2} to construct a solution
of \eqref{0205}, \eqref{0206}, which is simper than the above construction process.
\end{rem}

 The next lemma provides an estimate for the $H^k(\mathbb{R})$-norm of
the solution $({\varphi},{\theta},{\psi},{\pi})$ with $\xi$ varying,
which will be useful in the next subsection. To emphasize the dependence on
$\xi$, we write these solutions as
$(\varphi ,\theta ,\psi ,\pi )(\xi):=(\varphi ,\theta ,\psi ,\pi )(\xi,x_3)$.

\begin{lem}\label{lem:0201}
Let $k\geq 0$, and $\mathbb{D}\subset \mathbb{A}^{\mm{g}}$ be a nonempty bounded set satisfying
the closure  $\bar{\mathbb{D}}\subset \mathbb{A}^{\mm{g}}$. Then for any $\xi\in \mathbb{D}$
there are positive constants $A_k$, $B_k$, $C_k$ and $D_k$, which may
depend on $g$, $M$, $\mu$, $\bar{\rho}$ and $\mathbb{D}$, such that
\begin{eqnarray}
&&\label{0252n} \|\psi(\xi)\|_{H^k(\mathbb{R})}\leq A_k,\\
&&\label{0253n}\|\varphi(\xi)\|_{H^k(\mathbb{R})}\leq B_k,\\
&&\label{0254n}\|\theta(\xi)\|_{H^k(\mathbb{R})}\leq C_k,\\
&&\label{0255n} \|\pi(\xi)\|_{H^k(\mathbb{R})}\leq D_k,
\end{eqnarray} where $(\varphi,\theta,\psi,\pi)(\xi)$ and $\lambda(\xi)$ be
constructed as in Theorem \ref{thm:0302}. Moreover,
\begin{equation}\label{0255}
\|\psi(\xi)\|_{L^2(\mathbb{R})}^2> 0.
\end{equation}
\end{lem}
\begin{pf} We still give the proof for the horizonal case only. Throughout this proof, we denote
 by $\tilde{c}$ a generic positive constant which may
vary from line to line, and may depend on $g$, $M$, $\mu$, $\bar{\rho}$ and $\mathbb{D}$.

We begin with the estimate \eqref{0252n}.
We first rewrite \eqref{0207} as
\begin{equation}\label{0256n}\begin{aligned}
\psi''''(\xi)=&\big[(\lambda^2\bar{\rho}+2\lambda\mu
|\xi|^2+(M\xi_1)^2)\psi''(\xi)+\lambda^2\bar{\rho}'\psi'(\xi)\\
&\ -|\xi|^2(\lambda^2\bar{\rho}
+\lambda\mu |\xi|^2+(M\xi_1)^2-g\bar{\rho}')\psi(\xi)\big]/\lambda\mu,
\end{aligned}\end{equation}
which yields
\begin{equation}\label{0257n}\begin{aligned}
\|\psi''(\xi)\|_{L^2(\mathbb{R})}=&(\lambda\mu)^{-1}\int_\mathbb{R}\big[
\lambda(\lambda\bar{\rho}+2\mu
|\xi|^2+(M\xi_1)^2)\psi''(\xi)+\lambda^2\bar{\rho}'\psi'(\xi)\\
&\ -|\xi|^2(\lambda^2\bar{\rho}
+\lambda\mu |\xi|^2+(M\xi_1)^2-g\bar{\rho}')\psi(\xi)\big]\psi(\xi)
\mm{d}x_3
\end{aligned}\end{equation}
and
\begin{equation}\label{0257nn}\begin{aligned}
\|\psi'''(\xi)\|_{L^2(\mathbb{R})}=&-(\lambda\mu)^{-1}\int_\mathbb{R}\big[
\lambda(\lambda\bar{\rho}+2\mu
|\xi|^2+(M\xi_1)^2)\psi''(\xi)+\lambda^2\bar{\rho}'\psi'(\xi)\\
&\ -|\xi|^2(\lambda^2\bar{\rho}
+\lambda\mu |\xi|^2+(M\xi_1)^2-g\bar{\rho}')\psi(\xi)\big]\psi''(\xi)
\mm{d}x_3.
\end{aligned}\end{equation}
Noting that since  $\psi\in \mathcal{A}$, the inequality \eqref{0255}
holds, and there is a constant $\tilde{c}$, such that
\begin{equation}\label{0258n}
\|\psi(\xi)\|_{H^1(\mathbb{R})}\leq \tilde{c}. \end{equation}
On the other hand, in view of Proposition \ref{pro:0208}, we have
\begin{equation}\label{0259n}
 \lambda(\xi) \geq \tilde{c}>0  \quad\mbox{for any }\xi\in \mathbb{D}.
\end{equation}
Thus, using Cauchy-Schwarz's inequality, \eqref{0258n} and
\eqref{0259n}, we deduce from \eqref{0257n} and \eqref{0257nn} that
\begin{equation*} \begin{aligned}
\|\psi''(\xi)\|_{H^1(\mathbb{R})}\leq \tilde{c},
\end{aligned}\end{equation*}
whence, by \eqref{0256n},
\begin{equation}\label{0262n}\begin{aligned}
\|\psi(\xi)\|_{H^4(\mathbb{R})}\leq \tilde{c}.
\end{aligned}\end{equation}
Consequently, differentiating \eqref{0256n} with respect to $x_3$ and
using \eqref{0262n}, we find, by induction on $k$, that \eqref{0252n}
holds for any $k\geq 0$.

Now we turn to the estimate of \eqref{0253n}. From \eqref{0250}, we have
$$\|\varphi'\|_{L^2(\mathbb{R})}^2+\sigma\|
\varphi\|_{L^2(\mathbb{R})}^2=\int_\mathbb{R}\omega\varphi\mm{d}x_3, $$
which, together with Cauchy-Schwarz's inequality and \eqref{0252n},
yields
\begin{equation*}\|\varphi\|_{H^1(\mathbb{R})}^2\leq \tilde{c}.
\end{equation*}
We further deduce from \eqref{0250} that
\begin{equation}\label{0263n}\|\varphi\|_{H^2(\mathbb{R})}^2\leq \tilde{c}.
\end{equation}
Thus, we can employ \eqref{0263n}
and  \eqref{0250} to deduce that \eqref{0253n}
holds for any $k\geq 0$. Using this fact, the estimates \eqref{0252n} and \eqref{0253n},
the expression of $\pi$ in \eqref{0249}
implies \eqref{0255n}. Finally, similarly to the deduction of \eqref{0253n}, we can show
that \eqref{0254n} holds. This completes the proof. \hfill $\Box$
\end{pf}

\subsection{Exponential growth rate}\label{sec:0204}

In this subsection we will construct a linear real-valued solution to the linearized problem
 \eqref{0108} along with \eqref{0107} which grows in-time in the Sobolev space of order $k$.
\begin{thm}\label{thm:0203}
Under the assumptions of Theorem \ref{thm:0101},
let \begin{equation}\label{0267}
\Lambda=\sup_{\xi\in \mathbb{A}^\mm{g}\cap (L^{-1}\mathbb{Z})^2}\lambda(\xi),
\end{equation}
then there exist a positive constant $\Lambda^*\in (2\Lambda/3, \Lambda]$, and a real-valued solution
$(\varrho,\mathbf{u},\mathbf{N},q)$  to the linearized problem
 \eqref{0107}--\eqref{0108}
defined on the horizontally periodic domain $\Omega$, such that
\begin{enumerate}[\quad (1)]
  \item For every $k\in \mathbb{N}$,
\begin{equation}\label{0262}
\|(\varrho,\mathbf{u},\mathbf{N},q)(0)\|_{H^k}<\infty¡£
\end{equation}
  \item For every $t>0$, $(\varrho,\mathbf{u},\mathbf{N},q)(t)\in H^k$, and
\begin{eqnarray}
&&\label{0263} e^{\Lambda^* t}\|(\varrho,\mf{N},q)(0)\|_{H^k}
=\|(\varrho,\mf{N},q)(t)\|_{H^k}= e^{t\Lambda^*}\|(\varrho,\mf{N},q)(0)\|_{H^k},\\[1mm]
&&\label{0264}e^{\Lambda^* t}\|u_i(0)\|_{H^k}= \|u_i(t)\|_{H^k}=
e^{t\Lambda^*}\|u_i(0)\|_{H^k},\quad i=1,2,3.
\end{eqnarray}
  \item Moreover, \begin{equation}\label{02div}
\mm{div}\,\mf{u}(0)=\mm{div}\,\mf{N}(0)=0, \end{equation}
and
\begin{equation}\label{0268}
\|({u}_1,u_2)(0)\|_{L^2}\|{u}_3(0)\|_{L^2}>0.   \end{equation}
\item In addition, if $\mf{\bar{M}}=M\mf{e}_3$, then
\begin{equation}\label{0268withN3} \|N_3(0)\|_{L^2}>0. \end{equation}
\end{enumerate}
\end{thm}
\begin{rem}
In view of the third conclusion in Proposition \ref{pro:0202} and \eqref{0241}, we see that
 $$0<\Lambda\leq \sup_{\xi\in \mathbb{A}^{\mm{g}}}\lambda(\xi)\leq
\sqrt{ g\left\|{\bar{\rho}'}/{\bar{\rho}}\right \|_{L^\infty(\mathbb{R})} }.$$
\end{rem}
\begin{pf}
Let
\begin{equation*}\mathbf{v}(\xi,x_3)=-\mathrm{i}\varphi(\xi,x_3)\mf{e}_1-\mathrm{i}\theta(\xi,x_3)
\mf{e}_2+\psi(\xi, x_3)\mf{e}_3,
\end{equation*}
where $(\varphi,\theta,\psi)$ with an associated growth rate $\lambda(\xi)$ is constructed
in Theorem \ref{thm:0302} for any given $\xi\in\mathbb{A}^{\mm{g}}$.
 Recalling the definition of supremum, and using Proposition \ref{pro:0202}, we find that
 there exist $\xi^1$ and $\xi^2:=-{\xi}^1$ such that
  $$\Lambda^*:=\lambda(\xi^1)=\lambda(\xi^2)\in (2\Lambda/3, \Lambda]\mbox{ and }
  \xi^i\in \mathbb{A}^{\mm{g}}\cap (L^{-1}\mathbb{Z})^2\mbox{ for }i=1\mbox{ and }2.$$
 In view of  Remark \ref{rem:0201},
\begin{eqnarray}\label{linearsolution}\left\{\begin{array}{l}
                 {\varrho} (t,\textbf{x})=-e^{\Lambda^*t}\bar{\rho}'\sum_{m=1}^2{v}_3(\xi^m,x_3)
e^{\mathrm{i}\mf{x}'\xi^m},\\[2mm]
                 {\mathbf{u}}(t,\textbf{x})=\Lambda^*e^{\Lambda^*t}\sum_{m=1}^2
\mathbf{v}(\xi^m,x_3)e^{\mathrm{i}\mf{x}'\xi^m},\\[2mm]
{q}(t,\textbf{x})=\Lambda^*e^{\Lambda^*t}\sum_{m=1}^2{\pi}(\xi^m,x_3)e^{\mathrm{i}\mf{x}'\xi^m},\\[2mm]
{\mf{N}}(t,\mf{x})=e^{\Lambda^*t}\sum_{m=1}^2\bar{\mf{M}}\cdot\nabla_{\mf{x}}(
\mathbf{v}(\xi^m,x_3)e^{\mathrm{i}\mf{x}'\xi^m})
                 \end{array}\right.
\end{eqnarray}
gives a horizontally periodic, real-value solution to \eqref{0106}--\eqref{0108} satisfying \eqref{0263}--\eqref{02div}.
Lemma \ref{lem:0201} immediately implies that $(\varrho,\mathbf{u},\mathbf{N},q)$
constructed in \eqref{linearsolution} also satisfies \eqref{0262}.
Finally, thanks to \eqref{0255}, we get
 \begin{equation*} \begin{aligned}
\|{u}_3(0)\|_{L^2}^2=&4(\Lambda^*)^2\int_{\mathbb{R}}\psi^2(\xi^1,x_3)
\mm{d}x_3\int_{(2\pi L\mathbb{T})^2}\cos^2(\mf{x}'\cdot\xi^1)\mathrm{d}\mf{x}'>0,
\end{aligned}\end{equation*}
which, together with \eqref{02div}, implies \eqref{0268}.

 When $\mf{\bar{M}}=M\mf{e}_3$, recalling the expression of $\mf{M}$,
one gets from a simple computation that
$$\|N_3(0)\|_{L^2}^2=4M^2
\int_{\mathbb{R}}(\partial_{x_3}\psi(\xi^1,x_3))^2\mm{d}x_3
\int_{(2\pi L \mathbb{T})^2}\cos^2(\mf{x}'\cdot\xi^1)\mm{d}\mf{x}'.$$
On the other hand, thanks to \eqref{0255} and $\psi\in H^\infty(\mathbb{R})$, we have
$\partial_{x_3}\psi(\xi^1,x_3)\not\equiv 0$. Consequently, \eqref{0268withN3} follows.
This completes the proof of Theorem \ref{thm:0203}.
 \hfill $\Box$
\end{pf}

\section{Nonlinear energy estimates for the perturbed problem}\label{sec:04}

 In this section, we derive some nonlinear energy estimates for the perturbed   problem
\eqref{0105}--\eqref{0107}, and the integrand form of Gronwall's inequality of the high-order energy
estimates.

First we shall give a local well-posedness result of the perturbed
problem. We mention that the local existence of classical solutions and global existence
of classical small solutions to the non-resistive MHD equations
  have been established by many authors, see \cite{CLFDSMCRJLR,KSS1P,LXLSNWDHL,LFZPGC,HXPLFHG} for example.
To our best knowledge, there is no existence result for the MHD equations \eqref{0101} in the horizontally periodic domain $\Omega$.
However, with the help of the usual approaches in the proof of local existence for fluid dynamical equations and some mathematical
techniques to deal with the horizontally periodic domain, we can establish the following local well-posedness result,
the proof of which will be offered in Section \ref{appendix} for the completeness.
\begin{pro} \label{pro:0401new}
For any given initial data $(\varrho_0,\mf{u}_0,\mf{N}_0)\in H^3\times H^4\times H^3$ satisfying
$\inf_{
\mathbf{x}\in\Omega}\{(\varrho_0+\bar{\rho})(\mathbf{x})\}
>0$ and
 $\mm{div}\mf{u}_0=\mm{div}\mf{N}_0=0$, there exist a $T>0$ and a classical unique solution
$(\varrho,\mf{u},\mf{N})\in C^0(\bar{I}_T,H^3\times H^4\times H^3)$ to the perturbed 
problem \eqref{0105}--\eqref{0107} with an associated pressure $q$.
\end{pro}

With Proposition \ref{pro:0401new} in hand, we further derive the integrand form of Gronwall's inequality of the high-order energy
estimates. It should be pointed out that the solution $(\varrho,\mf{u},\mf{N},q)$ constructed
in Proposition \ref{pro:0401new} possesses more regularity, see the proof in Section \ref{appendix}.
This additional regularity makes it sense to derive some identities   and high-order energy estimates later. In particular,
 we omit the standard regularization argument in the derivation of some identities, \eqref{uttinform} for example,
 we refer to \eqref{eq0510} or \eqref{fracrho} in Section \ref{appendix} for the proof.

In what follows, the notation $a\lesssim b$ means that $a\leq Cb$ for a universal constant $C>0$,
which may depend on some physical parameters in \eqref{0105}.  $C(\delta_0)$ means that the positive constant $C$ further depends
on $\delta_0$. We denote
$$\mathcal{E}:=\mathcal{E}(t):=\mathcal{E}(\varrho,\mathbf{u},\mathbf{N}):=
({\|(\varrho,\mathbf{N})\|^2_{H^3}+\|\mathbf{u}\|_{H^4}^2})^{{1}/{2}}$$
and $\mathcal{E}_0={\mathcal{E}}(\varrho_0,\mathbf{u}_0,\mathbf{N}_0)$.
$\alpha=(\alpha_1,\alpha_2,\alpha_3)$ denotes a multi-index of order
$|\alpha|=\alpha_1+\alpha_2+\alpha_3$.

\subsection{Estimate of the perturbation density}

We first note that
by the classical Sobolev embedding results (see \cite[Chapter 5]{ARAJJFF}),
we have
 \begin{eqnarray}
&&\label{0403} \|u\|_{L^4}\lesssim  \|u\|_{L^2}^{\frac{1}{4}}\| u\|_{H^1}^{\frac{3}{4}}\lesssim \|u\|_{H^{1}},\\
 && \label{n0404}\|u\|_{L^\infty}\lesssim \|  u\|_{L^2}^{\frac{1}{4}}\|  u\|_{H^2}^{\frac{3}{4}}\lesssim\|  u\|_{H^2},\\
 && \label{0405} \|u\|_{H^{j}}\lesssim \|u\|_{L^2}^{\frac{1}{j+1}}\|u\|_{H^{j+1}}
 ^{\frac{j}{j+1}}.
\end{eqnarray}
With the help of the above estimates, we can bound the perturbation density $\varrho$.
In fact, applying $\partial^\alpha$ to \eqref{0105}$_1$, multiplying the resulting identity
by $\partial^\alpha\varrho$ in $L^2(\Omega)$, we get
 \begin{equation}\label{0406}\begin{aligned}
\frac{1}{2}\frac{d}{dt}\sum_{0\leq |\alpha |\leq 3}\int_{\Omega}(\partial^\alpha\varrho)^2\mm{d}\mf{x}=&
-\sum_{0\leq |\alpha |\leq 3}\int_{\Omega}\partial^\alpha(\bar{\rho}' u_3)\partial^\alpha \varrho\mm{d}\mf{x}+
\sum_{0\leq |\alpha |\leq 3}\int_{\Omega} \partial^\alpha (\mf{u}\cdot\nabla \varrho)
\partial^\alpha\varrho\mm{d}\mf{x}\\
:=& I_1+I_2,
\end{aligned} \end{equation}
where we have defined
$$\partial^\alpha (\mf{u}\cdot\nabla \varrho)= \sum_{\alpha=\beta+\gamma\atop |\gamma|\leq 2}
\partial^\beta\mf{u}\cdot\partial^\gamma\nabla \varrho\mbox{ for }|\alpha|=3,$$
so that \eqref{0406} makes sense (cf. Lemma \ref{lem:0501}).

Using \eqref{0405}, H\"older's and Young's inequalities,
$I_1$ can be bounded as follows:
 \begin{equation}\label{0407}\begin{aligned}
I_1\lesssim  &\sum_{0\leq |\alpha |\leq 3}\|\partial^\alpha(\bar{\rho}'u_3)\|_{L^2}\|\partial^\alpha \varrho\|_{L^2}\\
\lesssim &\|\varrho\|_{L^2}\|u_3\|_{L^2} +\|\nabla \varrho\|_{L^2}(\|  u_3\|_{L^2}
+\|\nabla u_3\|_{L^2})\\
&+\|\nabla^2 \varrho\|_{L^2}(\| u_3\|_{L^2}+\|\nabla  u_3\|_{L^2}+ \|\nabla^2 u_3\|_{L^2}) \\
&+\|\nabla^3 \varrho\|_{L^2}( \| u_3\|_{L^2}+\|\nabla  u_3\|_{L^2}+
\|\nabla^2 u_3\|_{L^2}+\|\nabla^3 u_3\|_{L^2})\\
\lesssim &\|\varrho\|_{H^3}\|u_3\|_{L^2}+\|\nabla \varrho\|_{H^2}
\| u_3\|_{L^2}^{\frac{1}{2}}\| u_3\|_{H^2}^{\frac{1}{2}}
+\|\nabla^2 \varrho\|_{H^1}\|  u_3\|_{L^2}^{\frac{1}{3}}
\| u_3\|_{H^3}^{\frac{2}{3}}\\
&+\|\nabla^3 \varrho\|_{L^2}\|  u_3\|_{L^2}^{\frac{1}{4}}
\| u_3\|_{H^3}^{\frac{3}{4}}\lesssim   C(\delta_0)\|u_3\|_{L^2}^2
 +\delta_0\mathcal{E}^2. \end{aligned} \end{equation}
On the other hand, we use \eqref{0403}, \eqref{n0404} and
H\"older's inequality to control $I_2$ as follows.
\begin{equation}\label{0408}\begin{aligned}
I_2\lesssim  &\|\nabla\mf{u}\|_{L^4}\|\nabla \varrho
\|_{L^4}\|\nabla \varrho\|_{L^2}+\|\nabla^2\mf{u}\|_{L^4}\|\nabla\varrho\|_{L^4}
\|\nabla^2\varrho\|_{L^2}+\|\nabla \mf{u}\|_{L^\infty}\|\nabla^2\varrho\|_{L^2}^2  \\
&  +\|\nabla^3\mf{u}\|_{L^4}\|\nabla\varrho\|_{L^4}\|\nabla^3\varrho\|_{L^2}
+\|\nabla^2 \mf{u}\|_{L^\infty}\|\nabla^2\varrho\|_{L^2}\|\nabla^3\varrho\|_{L^2} +\|\nabla \mf{u}\|_{L^\infty}\|\nabla^3\varrho\|_{L^2}^2 \\
\lesssim  &\|\nabla\mf{u}\|_{H^1}\|\nabla \varrho\|_{H^1}\|\nabla \varrho\|_{L^2}+\|\nabla\mf{u}\|_{H^2}
\|\nabla\varrho\|_{H^1}\|\nabla^2\varrho\|_{L^2}+\|\nabla \mf{u}\|_{H^2}\|\nabla^2\varrho\|_{L^2}^2\\
&+\|\nabla \mf{u}\|_{H^3}\|\nabla\varrho\|_{H^2}^2\lesssim \mathcal{E}^3.
\end{aligned} \end{equation}
Hence, one gets from \eqref{0406}--\eqref{0408} that
  \begin{equation}\label{density}\begin{aligned}
\frac{1}{2}\frac{d}{dt}\|\varrho(t)\|^2_{H^3}\lesssim
C(\delta_0)\|u_3(t)\|_{L^2}^2+ \delta_0\mathcal{E}^2(t)+ \mathcal{E}^3(t).
\end{aligned} \end{equation}
In addition, we can deduce from \eqref{0105}$_1$ that
\begin{eqnarray}
&& \label{0410}\|\varrho_t\|_{L^2} \lesssim
 \| {u}_3\|_{L^2}  + \mathcal{E}^2,\\
&& \label{0411}
\|\nabla \varrho_t\|_{L^2}\lesssim \| {u}_3\|_{L^2} +\|\nabla  {u}_3\|_{L^2}  +\mathcal{E}^2
\lesssim  C(\delta_0)\| {u}_3\|_{L^2}  + \delta_0\mathcal{E}+\mathcal{E}^2,\\
 && \label{0411nn}\|\varrho_{tt}\|_{L^2}\lesssim \|\mf{u}_t\|_{L^2}(1+\mathcal{E})+\mathcal{E}^2+\mathcal{E}^3.
\end{eqnarray}

\subsection{Estimate of the perturbation magnetic field}
We continue to bound the perturbation magnetic field $\mf{N}$.
 Applying $\partial^\alpha$ to \eqref{0105}$_3$, multiplying the resulting identity by
 $\partial^\alpha\mathbf{N}$ in $L^2$, we obtain
 \begin{equation*} \begin{aligned}
&\frac{1}{2}\frac{d}{dt}\sum_{0\leq |\alpha |\leq 3}\int_{\Omega}|\partial^\alpha \mf{N}|^2\mm{d}\mf{x}\\
&=
\sum_{0\leq |\alpha |\leq 2}\int_{\Omega}\partial^\alpha(\bar{\mf{M}}
\cdot\nabla\mf{u})\cdot\partial^\alpha \mf{N} \mm{d}\mf{x}+
\sum_{0\leq |\alpha |\leq 3}\int_{\Omega} [\partial^\alpha
(\mf{N}\cdot\nabla \mf{u})-\partial^\alpha(\mf{u}\cdot\nabla \mf{N})]\cdot\partial^\alpha\mf{N}\mm{d}\mf{x}\\
&\quad  +\sum_{ |\alpha |= 3}\int_{\Omega}\partial^\alpha(\bar{\mf{M}}
\cdot\nabla\mf{u})\cdot\partial^\alpha \mf{N}
\mm{d}\mf{x}:= J_1+J_2+J_3,
\end{aligned} \end{equation*}
where we have used the facts $\nabla \times (\mf{u}\times \bar{\mf{M}})=\bar{\mf{M}}\cdot\nabla \mf{u}$ and
$\nabla \times (\mf{u}\times \mf{N})=\mf{N}\cdot\nabla \mf{u}-\mf{u}\cdot\nabla \mf{N}$,
and defined
$$\partial^\alpha (\mf{u}\cdot\nabla\mf{N})= \sum_{\alpha=\beta+\gamma\atop |\gamma|\leq 2}
\partial^\beta\mf{u}\cdot\partial^\gamma\nabla \mf{N}\mbox{ for }|\alpha|=3.$$

Similarly to the derivation of \eqref{0407} and \eqref{0408}, we can control $J_1$ and $J_2$ as follows.
 \begin{equation*} \begin{aligned}
J_1\lesssim  &\sum_{0\leq |\alpha |\leq 2} \|\partial^\alpha\nabla \mf{u}\|_{L^2}\|\partial^\alpha\mf{N}\|_{L^2}\\
\lesssim & \|\nabla\mf{u}\|_{L^2}\|\mf{N}\|_{L^2}+
\|\nabla^2\mf{u}\|_{L^2}\|\nabla \mf{N}\|_{L^2}+
\|\nabla^3\mf{u}\|_{L^2}\|\nabla^2 \mf{N}\|_{L^2}\\
\lesssim & \|\mf{u}\|_{L^2}^{\frac{1}{2}}\|\mf{u}\|_{H^2}^{\frac{1}{2}}
\|\mf{N}\|_{L^2}+\|\mf{u}\|_{L^2}^{\frac{1}{3}}\|\mf{u}\|_{H^3}^{\frac{2}{3}}\|\nabla \mf{N}\|_{L^2}+
\|\mf{u}\|_{L^2}^{\frac{1}{4}}\|\mf{u}\|_{H^4}^{\frac{3}{4}}\|\nabla^2 \mf{N}\|_{L^2}\\
\lesssim  & C(\delta_0)\|\mf{u}\|_{L^2}^2+\delta_0\mathcal{E}^2,
\end{aligned} \end{equation*}
and
 \begin{equation*} \begin{aligned}
J_2\lesssim  &\|\nabla \mf{u}\|_{L^4}\|\mf{N}\|_{L^4}\|\mf{N}\|_{L^2}
+(\|\nabla\mf{u}\|_{L^4}\|\nabla \mf{N} \|_{L^4}+\|\nabla^2\mf{u}\|_{L^4}\|\mf{N}
\|_{L^4})\|\nabla\mf{N}\|_{L^2}\\
& +(\|\nabla\mf{u}\|_{L^\infty}\|\nabla^2 \mf{N}\|_{L^2}+\|\nabla^2\mf{u}\|_{L^4}\|\nabla \mf{N}
\|_{L^4}+\|\nabla^3\mf{u}\|_{L^2}\|\mf{N} \|_{L^\infty})\|\nabla^2\mf{N}\|_{L^2}\\
& +(\|\nabla\mf{u}\|_{L^\infty}\|\nabla^3 \mf{N}
\|_{L^2}+\|\nabla^2\mf{u}\|_{L^4}\|\nabla^2 \mf{N} \|_{L^4}\\
&\qquad +\|\nabla^3\mf{u}\|_{L^2}\|\nabla \mf{N}\|_{L^\infty}+\|\nabla^4\mf{u}\|_{L^2}\| \mf{N}
\|_{L^\infty})\|\nabla^3\mf{N}\|_{L^2}\\
\lesssim  &\|\nabla \mf{u}\|_{H^1}\|\mf{N}\|_{H^1}\|\mf{N}\|_{L^2}+\|\nabla\mf{u}\|_{H^2}\| \mf{N}
\|_{H^2}\|\nabla\mf{N}\|_{L^2} +\|\nabla\mf{u}\|_{H^2}\|  \mf{N}
\|_{H^2}\|\nabla^2\mf{N}\|_{L^2}\\
&+\|\nabla\mf{u}\|_{H^3}\| \mf{N} \|_{H^3}\|\nabla^3\mf{N}\|_{L^2}
\lesssim  \mathcal{E}^3.
\end{aligned} \end{equation*}
Hence, one gets
  \begin{equation}\label{0415}\begin{aligned}
\frac{1}{2}\frac{d}{dt}\|\mf{N}(t)\|^2_{H^3}\lesssim C(\delta_0)\|\mf{u}(t)\|_{L^2}^2+ \delta_0\mathcal{E}^2(t)
+\mathcal{E}^3(t) +J_3.
\end{aligned} \end{equation}
In addition, we infer from \eqref{0105}$_4$ that
\begin{eqnarray}&& \label{0416} \begin{aligned}
\| \mf{N}_t\|_{L^2}\lesssim & \|\nabla\mf{u} \|_{L^2}+ \mathcal{E}^2
\lesssim C(\delta_0)\| \mf{u} \|_{L^2}+\delta_0\mathcal{E}
+\mathcal{E}^2,\end{aligned}\\
&&\label{0417}
\|\nabla\mf{N}_t\|_{L^2}\lesssim  \|\nabla^2\mf{u} \|_{L^2}+ \mathcal{E}^2
\lesssim C(\delta_0)\|\mf{u} \|_{L^2}+\delta_0\mathcal{E}+\mathcal{E}^2,\\
&&\label{ttN0417}
\|\mf{N}_{tt}\|_{L^2}\lesssim  \|\mf{u}_t \|_{H^1}(1+ \mathcal{E})+\mathcal{E}^2+\mathcal{E}^3.
\end{eqnarray}

\subsection{Estimate of the perturbation velocity}
First we restrict the density $\rho=\varrho+\bar{\rho}$ so that it has a positive lower bound.
By virtue of \eqref{n0404}, there is a constant $\delta_0'\in (0,1)$, such that
\begin{equation}\label{conditiondensity}
\|\varrho_0\|_{L^\infty}\leq \frac{\inf_{{x}_3\in\mathbb{R}}\{\bar{\rho}{(x_3)}\}}{2}\;\;
\mbox{ for any }\|\varrho_0\|_{H^2}\leq \delta_0'.
\end{equation}
Consequently, in view of \eqref{0105}$_1$ and $\mm{div}\mf{u}=0$, we find that for any
$(t,\mf{x})\in [0,T)\times \Omega$ by the method of characteristics,
\begin{equation*}\label{0418}\frac{\inf_{{x}_3\in\mathbb{R}}\{\bar{\rho}({x_3})\}}{2}
\leq \inf_{\mathbf{x}\in\Omega}\{\rho_0(t,\mathbf{x})\}
\leq \rho(t,\mf{x})\leq \sup_{\mathbf{x}\in\Omega}\{{\rho}_0(t,{\mathbf{x}})\}
\leq \frac{3\sup_{{x}_3\in\mathbb{R}}\{\bar{\rho}({x_3})\}}{2},
\end{equation*}
where $\rho_0=\varrho_0+\bar{\rho}$.
From now on, we always assume that $\mathcal{E}\leq \delta_0\leq \delta_0'<1$. We remark that
the upper and lower boundedness of $\rho$ will be repeatedly used below.

We now estimate the time-derivative of the perturbation velocity.
Multiplying \eqref{0105}$_2$ by $\mf{u}$ in $L^2$, utilizing \eqref{0105}$_1$
and integrating by parts, we obtain
\begin{equation*} \begin{aligned} \frac{1}{2}\frac{{d}}{{d}t}
\int_{\Omega}(\varrho+\bar{\rho} )|\mf{u}|^2\mm{d}\mf{x}+
\mu\int_{\Omega}|\nabla \mf{u}|^2\mm{d}\mf{x}
=&\int_{\Omega}[(\nabla\times \mf{N})\times (\mf{N}+\bar{\mf{M}})
-g \varrho \mf{e}_3]\cdot\mf{u}\mm{d}\mf{x} \\
:=& \int_{\Omega}\mf{Y}_0\cdot\mf{u}\mm{d}\mf{x},
\end{aligned}\end{equation*}
where $\mf{Y}_0:=(\nabla\times \mf{N})\times (\mf{N}+\bar{\mf{M}})-g \varrho \mf{e}_3$.
Applying $\partial_t$ to \eqref{0105}$_2$, multiplying the resulting identity by $\mf{u}_t$
in $L^2$, we make use of \eqref{0105}$_1$,
and integrate by parts to deduce
\begin{equation*} \begin{aligned}&\frac{1}{2}\frac{{d}}{{d}t}
\int_{\Omega}(\varrho+\bar{\rho} )|\mf{u}_t|^2\mm{d}\mf{x}+
\mu\int_{\Omega}|\nabla \mf{u}_t|^2\mm{d}\mf{x}\\
& = \int_{\Omega}\Big[ (\nabla\times \mf{N}_t)\times (\mf{N}+\bar{\mf{M}})
+(\nabla\times \mf{N})\times \mf{N}_t-\varrho_t\mf{u}\cdot\nabla \mf{u}\\
&\quad  -(\varrho+\bar{\rho})\mf{u}_t\cdot\nabla \mf{u}-g \varrho_t \mf{e}_3
-\varrho_t\mf{u}_t\Big]\cdot\mf{u}_t
\mm{d}\mf{x}:=\int_{\Omega}\mf{Y}_1\cdot\mf{u}_t\mm{d}\mf{x}.
\end{aligned}\end{equation*}
Applying $\partial_{t}^2$ to \eqref{0105}$_2$ again, multiplying the resulting identity by $\mf{u}_t$,
employing \eqref{0105}$_1$ and
\begin{equation*}
(\nabla\times \mf{N})\times (\mf{N}+\bar{\mf{M}})=(\mf{N}+\bar{\mf{M}})
\cdot \nabla \mf{N}-\frac{1}{2}\nabla|(\mf{N}+\bar{\mf{M}})|^2,
\end{equation*}
we have
\begin{align}
 &\frac{1}{2}\frac{{d}}{{d}t}
\int_{\Omega}(\varrho+\bar{\rho} )|\mf{u}_{tt}|^2\mm{d}\mf{x}+
\mu\int_{\Omega}|\nabla \mf{u}_{tt}|^2\mm{d}\mf{x}\nonumber\\
&=\int_{\Omega}\Big\{-\mf{N}_{tt}\cdot\nabla\mf{u}_{tt}\cdot\mf{N}
-\mf{N}_{t}\cdot\nabla\mf{u}_{tt}\cdot\mf{N}_t-(\mf{N}
+\bar{\mf{M}})\cdot\nabla \mf{u}_{tt}\cdot\mf{N}_{tt}-\Big[g\varrho_{tt} \mf{e}_3-\varrho_{tt}\mf{u}_t \nonumber\\
&\quad  -\varrho_{t}\mf{u}_{tt}-\varrho_{tt}\mf{u}\cdot\nabla \mf{u}
-\varrho_{t}\mf{u}_t\cdot\nabla \mf{u} -\varrho_{t}\mf{u}\cdot\nabla \mf{u}_t
-(\varrho+\bar{\rho})\mf{u}_t\cdot\nabla \mf{u}_t
-(\varrho+\bar{\rho})\mf{u}_{tt}\cdot\nabla \mf{u}\Big]\cdot\mf{u}_{tt}\Big\} \mm{d}\mf{x}\nonumber\\
&\label{uttinform}:=-\sum_{i=0}^2
\int_{\Omega}\partial_t^i(\mf{N}+\bar{\mf{M}})\cdot\nabla\mf{u}_{tt}\cdot\partial_t^{2-i}\mf{N}
\mm{d}\mf{x}+\int_{\Omega}\mf{Y}_2\cdot\mf{u}_{tt}\mm{d}\mf{x}.
\end{align}
Adding the above three equalities, we get
\begin{equation}\label{0421}\begin{aligned}&\frac{1}{2}\frac{{d}}{{d}t}\sum_{i=0}^2
\int_{\Omega}(\varrho+\bar{\rho} )|\partial_t^i\mf{u}|^2\mm{d}\mf{x}+\mu\sum_{i=0}^2
\int_{\Omega}|\nabla \partial_t^i\mf{u}|^2\mm{d}\mf{x}\\
=&\sum_{i=0}^2
\int_{\Omega}\mf{Y}_i\cdot\partial_t^i\mf{u}\mm{d}\mf{x}-\sum_{i=0}^2
\int_{\Omega}\partial_t^i(\mf{N}+\bar{\mf{M}})\cdot\nabla\mf{u}_{tt}\cdot\partial_t^{2-i}\mf{N}
\mm{d}\mf{x}.
\end{aligned}\end{equation}

Obviously,
\begin{equation}\label{0422}\begin{aligned}\|\mf{Y}_0\|_{L^2}\lesssim & \|\varrho\|_{L^2}
+\|\mf{N}\|_{L^2}^{\frac{1}{2}}\| \mf{N}\|_{H^2}^{\frac{1}{2}}\|\mf{N}\|_{L^\infty}
\lesssim C(\delta_0)\|(\varrho,\mf{N})\|_{L^2}+\delta_0\mathcal{E}^2.
\end{aligned}\end{equation}
Recalling that $\mathcal{E}\leq \delta<1$, we use \eqref{0403}, \eqref{n0404}, \eqref{0410},
\eqref{0411}, \eqref{0416} and \eqref{0417} to arrive at
\begin{equation}\label{0423}\begin{aligned}\|\mf{Y}_1\|_{L^2}\lesssim &
\|\varrho_t\|_{L^2} +\|\varrho_t\|_{L^2}\|\mf{u}\|_{L^\infty}\|\nabla \mf{u}\|_{L^\infty}
+ \|\varrho+\bar{\rho}\|_{L^\infty}\| \mf{u}_t\|_{L^2}\|\nabla \mf{u}\|_{L^\infty} \\
& +(1+\| \mf{N}\|_{L^\infty})\| \nabla \mf{N}_t\|_{L^2}+\|\nabla \mf{N}\|_{L^4}\|\mf{N}_t\|_{L^4}
+\|\varrho_t\|_{L^4}\|\mf{u}_t\|_{L^4}\\
\lesssim &\|\varrho_t\|_{L^2} +\| \mf{u}\|_{H^3}^2\|\varrho_t\|_{L^2}
+\|\varrho_t\|_{H^1}\| \mf{u}_t\|_{H^1}+\|\nabla \mf{u}\|_{H^2}\|\mf{u}_t\|_{L^2}  \\
 & +\| \nabla \mf{N}_t\|_{L^2}+ \|\mf{N}\|_{H^2}\| \mf{N}_t\|_{H^1}  \\
\lesssim & C(\delta_0)\|\mf{u} \|_{L^2}+\delta_0 \mathcal{E}+\mathcal{E}\| \mf{u}_t\|_{H^1}.
\end{aligned}\end{equation}
To bound $\mf{Y}_2$, we argue in a way similar to that used for
\eqref{0423}, with additional estimates \eqref{0411nn} and \eqref{0405}, to infer that
\begin{equation}\label{0423nb}\begin{aligned}\|\mf{Y}_2\|_{L^2}\lesssim &
\|\varrho_{tt}\|_{L^2}(1+\mathcal{E}^2+\|\mf{u}_t\|_{H^2})
+ \|\varrho_t\|_{H^1}(\|\mf{u}_{tt}\|_{L^2}^{\frac{1}{4}}\| \mf{u}_{tt}\|_{H^1}^{\frac{3}{4}}
+\mathcal{E}\|\mf{u}_{t}\|_{H^2})  \\[1mm]
& +\| \mf{u}_t\|_{H^1}\| \nabla \mf{u}_t\|_{H^1}+\mathcal{E}\| \mf{u}_{tt}\|_{L^2}\\
\lesssim & [\|\mf{u}_t\|_{H^1}(1+\mathcal{E})+\mathcal{E}^2](1+\mathcal{E}^2
+\|\mf{u}_t\|_{H^2})+ \mathcal{E} (\|\mf{u}_{tt}\|_{L^2}^{\frac{1}{4}}
\| \mf{u}_{tt}\|_{H^1}^{\frac{3}{4}} +\mathcal{E}\|\mf{u}_{t}\|_{H^2}).
\end{aligned}\end{equation}
Finally, using \eqref{0416}--\eqref{ttN0417}, we obtain
\begin{equation}\label{0423nb1}\begin{aligned}
&\left|
\sum_{i=0}^2
\int_{\Omega}\partial_t^i(\mf{N}+\bar{\mf{M}})\cdot\nabla\mf{u}_{tt}\cdot\partial_t^{2-i}\mf{N}
\mm{d}\mf{x}\right| \\[1mm]
&\lesssim  [(1+\mathcal{E})\|\mf{N}_{tt}\|_{L^2}+\|\mf{N}_t\|_{H^1}^2]
\|\nabla\mf{u}_{tt}\|_{L^2}\lesssim  [(1+\mathcal{E})\|\mf{u}_t\|_{H^1}
+ \mathcal{E}^2]\|  \mf{u}_{tt}\|_{H^1}.
\end{aligned}\end{equation}
Putting \eqref{0421}--\eqref{0423nb1} together and using Cauchy-Schwarz's inequality, we conclude
\begin{equation}\label{04041n}\begin{aligned}
\frac{1}{2} & \frac{d}{dt} \sum_{i=0}^2 (\varrho+\bar{\rho})
\|\partial_t^i\mathbf{u}(t)\|^2_{L^2}+\sum_{i=0}^2\|\nabla \partial_t^i \mathbf{u}(t)\|^2_{L^2}
\lesssim C(\delta_0)\|(\varrho,\mf{u},\mf{N})(t) \|_{L^2}^2+\delta_0\mathcal{E}^2(t) \\
& +(C(\delta_0)\|\mf{u}(t) \|_{L^2}+\delta_0 \mathcal{E}(t) +\mathcal{E}^2(t)
+\mathcal{E}(t)\|\mf{u}_t(t)\|_{H^1})\|\mf{u}_t(t) \|_{L^2}  \\
& + [\|\mf{u}_t(t)\|_{H^1}(1+\mathcal{E}(t))+\mathcal{E}^2(t)](1+\mathcal{E}^2(t)+\|\mf{u}_t(t)\|_{H^2})\|
 \mf{u}_{tt}(t) \|_{L^2} \\
& + \mathcal{E}(t) (\|\mf{u}_{tt}(t)\|_{L^2}^{\frac{1}{4}}\| \mf{u}_{tt}(t)\|_{H^1}^{\frac{3}{4}}
+\mathcal{E}(t)\|\mf{u}_{t}(t)\|_{H^2})\|\mf{u}_{tt}(t) \|_{L^2}
\\
& +[(1+\mathcal{E}(t))\|\mf{u}_t(t)\|_{H^1}+ \mathcal{E}^2(t)]^2.
\end{aligned}\end{equation}

We proceed to estimate higher derivatives of the perturbation velocity. Rewriting
 \eqref{0105}$_2$ as
 \begin{equation}\label{0424}
-\mu \Delta \mathbf{u}+ (\varrho+\bar{\rho}){\bf u}_t+\nabla q =(\nabla\times \mf{N})\times
(\mf{N}+\bar{\mf{M}})-g \varrho \mf{e}_3-(\varrho+\bar{\rho}){\bf u}\cdot\nabla {\bf u},
\end{equation}
and multiplying \eqref{0424} by $\mf{u}_t$ in $L^2$, we have
$$ \frac{\mu}{2}\frac{d}{dt}
\int_{\Omega}|\nabla \mathbf{u}(t)|^2\mm{d}\mf{x}+\int_{\Omega}(\varrho+\bar{\rho})| \mathbf{u}_t|^2
\mm{d}\mf{x}=\int_{\Omega}[\mf{Y}_0-(\varrho+\bar{\rho}){\bf u}\cdot\nabla
{\bf u}]\cdot\mathbf{u}_t\mm{d}\mf{x}. $$
Differentiating \eqref{0424} with respect to $t$ and multiplying the resulting equations by $\mf{u}_{tt}$
in $L^2$, we obtain
 \begin{equation*} \begin{aligned}
 & \frac{\mu}{2}\frac{d}{dt} \int_{\Omega}|\nabla\mathbf{u}_t(t)|^2\mm{d}\mf{x}
+\int_{\Omega}(\varrho+\bar{\rho})| \mathbf{u}_{tt}|^2
\mm{d}\mf{x} =\int_{\Omega}[\mf{Y}_1-(\varrho+\bar{\rho}){\bf u}\cdot\nabla
{\bf u}_t]\cdot\mathbf{u}_{tt} \mm{d}\mf{x}.
\end{aligned}\end{equation*}
Adding the above two equalities, we get
 \begin{equation}\label{nn425}\begin{aligned}
& \frac{\mu}{2}\frac{d}{dt}\sum_{i=0}^1
\int_{\Omega}|\nabla \partial_t^i\mathbf{u}(t)|^2\mm{d}\mf{x}
+\sum_{i=0}^1\int_{\Omega}(\varrho+\bar{\rho})|\partial_t^{i+1} \mathbf{u}|^2\mm{d}\mf{x}\\
& =\int_{\Omega}[\mf{Y}_0-(\varrho+\bar{\rho}){\bf u}\cdot\nabla
{\bf u}]\cdot \mathbf{u}_t\mm{d}\mf{x}+\int_{\Omega}\left[\mf{Y}_1
- (\varrho+\bar{\rho}){\bf u}\cdot\nabla {\bf u}_t\right]\cdot\mathbf{u}_{tt} \mm{d}\mf{x}.
\end{aligned}\end{equation}
Clearly,
 \begin{equation}\label{0422n}\begin{aligned}
 \|(\varrho+\bar{\rho}){\bf u}\cdot\nabla
{\bf u}\|_{L^2}+\| (\varrho+\bar{\rho}){\bf u}\cdot\nabla
{\bf u}_t\|_{L^2}\lesssim  \mathcal{E}^2+\mathcal{E}\|\nabla\mf{u}_t\|_{L^2}. \end{aligned}\end{equation}
Consequently, in view of \eqref{0423}, \eqref{0423nb} and \eqref{0422n}, using Cauchy-Schwarz's inequality,
we get from \eqref{nn425} that
 \begin{equation}\label{highvelcocity}\begin{aligned}
 & \frac{d}{dt}\sum_{i=0}^1 \|\nabla \partial_t^i\mathbf{u}(t)\|^2_{L^2}
+\sum_{i=0}^1\|\partial_t^{i+1} \mathbf{u}(t)\|^2_{L^2} \\
&\lesssim C(\delta_0)\|(\varrho,\mf{u},\mf{N})(t)\|_{L^2}^2
+\delta_0\mathcal{E}^2(t)+\mathcal{E}^2(t)\|\mf{u}_t(t)\|_{H^1}^2
.
\end{aligned} \end{equation}

 Next, we continue to derive more estimates of higher-order derivatives of the perturbation velocity.
Multiplying \eqref{0105}$_2$ by $\mf{u}_t$, integrating and using \eqref{0405}, one gets
 \begin{equation}\label{l2uestimate}\begin{aligned}
\|\mf{u}_t\|_{L^2}^2\lesssim&\|\sqrt{\varrho+\bar{\rho}}\mf{u}_t\|_{L^2}^2\lesssim
\|\varrho\|_{L^2}^2+\|\Delta\mf{u}\|_{L^2}^2+\|\nabla\mf{N}\|_{L^2}^2+\mathcal{E}^4\\
\lesssim
&C(\delta_0)\|(\varrho,\mf{u},\mf{N})\|_{L^2}^2+\delta_0\mathcal{E}^2,
\end{aligned}
\end{equation}
while applying $\partial_{i}$ to \eqref{0105}$_2$, multiplying the resulting equations
by $\partial_{i}\mf{u}_t$ in $L^2$, and using \eqref{0405}, we have
 \begin{equation}\label{utone}\begin{aligned}
\|\partial_{i}\mf{u}_t\|_{L^2}^2\lesssim &\|\sqrt{\varrho+\bar{\rho}}
\partial_{i}\mf{u}_t\|_{L^2}^2\lesssim C(\delta_0)\|(\varrho,\mf{u},\mf{N})\|_{L^2}^2+\delta_0\mathcal{E}^2+(1+\mathcal{E}^2)
\|\mf{u}_t\|_{L^2}^2\\
\lesssim & C(\delta_0)\|(\varrho,\mf{u},\mf{N})\|_{L^2}^2+\delta_0\mathcal{E}^2.
\end{aligned}
\end{equation}
If we take $\partial_{i}\partial_{j}$ to \eqref{0105}$_2$ and multiply the resulting equations
with $\partial_{i}\mf{u}_t$ in $L^2$, we see that
 \begin{equation*} \begin{aligned}
\|\partial_{i}\partial_{j}\mf{u}_t\|_{L^2}^2\lesssim \|\sqrt{\varrho +\bar{\rho}}
\partial_{i}\partial_{j}\mf{u}_t\|_{L^2}^2 \lesssim\mathcal{E}^2+(1
+\mathcal{E}^2)\|\mf{u}_t\|_{H^1}^2\lesssim\mathcal{E}^2.
\end{aligned}
\end{equation*}
In particular, summing up the above three estimates, we conclude
 \begin{equation}\label{partialu}
\|{\bf u}_t\|^2_{H^2} \lesssim \mathcal{E}^2. \end{equation}

Finally, adding \eqref{04041n} to \eqref{highvelcocity}, and utilizing \eqref{l2uestimate}--\eqref{partialu},
we arrive at
  \begin{equation}\label{0404}\begin{aligned}
&\frac{d}{dt}\left(\sum_{i=0}^1\|\nabla \partial_t^i\mathbf{u}(t)\|^2_{L^2}+\sum_{i=0}^2
\|\sqrt{\varrho +\bar{\rho}}\partial_t^i\mf{u}(t)\|^2_{L^2}\right)  \\
&\qquad +\sum_{i=0}^1 \|\partial_t^i(\nabla\mf{u},\mf{u}_t,\nabla\mf{u}_t)(t)\|^2_{L^2}\lesssim
C(\delta_0)\| (\varrho,\mf{u},\mf{N})(t) \|_{L^2}^2+\delta_0\mathcal{E}^2(t),
\end{aligned} \end{equation}
provided $\delta_0$ is sufficiently small.

\subsection{Energy estimates}\label{sec:0404}
Now, we sum up the previous estimates \eqref{density}, \eqref{0415} and \eqref{0404},
and use Cauchy-Schwarz's inequality to find that
  \begin{equation}\label{densityenegr}\begin{aligned}
&\frac{d}{dt}\left[\|(\varrho,\mf{N})(t)\|^2_{H^3}
+\sum_{i=0}^1\|\nabla \partial_t^i\mathbf{u}(t)\|^2_{L^2} +\sum_{i=0}^2
\|\sqrt{\varrho+\bar{\rho}}\,\partial_t^i\mf{u}(t)\|^2_{L^2}\right] \\
&\qquad +\sum_{i=0}^1 \| \partial_t^i(\nabla\mf{u},\mf{u}_t,\nabla\mf{u}_t)(t)\|^2_{L^2}\lesssim
C(\delta_0)\| (\varrho,\mf{u},\mf{N})(t) \|_{L^2}^2 +J_3+\delta_0\mathcal{E}^2(t),
\end{aligned} \end{equation}
provided $\delta_0$ is sufficiently small.

To deal with the term $J_3$, we shall make use of the momentum equations \eqref{0105}$_2$.
Let $\alpha:=(\alpha_1,\alpha_2,\alpha_3)$ be given and satisfy $|\alpha|=3$.
Without loss of generality, assume $\alpha_1\neq 0$. By employing a partial integration,
one sees
  \begin{equation}\label{mageequality}
\int_{\Omega}\partial^\alpha(\bar{\mf{M}}\cdot\nabla\mf{u})\cdot\partial^\alpha\mf{N}
\mm{d}\mf{x}=\int_{\Omega}\bar{\mf{M}}\cdot \nabla \partial^\gamma\mf{N}\cdot
\partial^\beta \mf{u}\mm{d}\mf{x},
\end{equation}
where $\gamma:=(\alpha_1-1, \alpha_2,\alpha_3)$ and $\beta:=(\alpha_1+1, \alpha_2,\alpha_3)$
satisfying $|\gamma|=2$ and $|\beta|=4$. Letting $\partial^\gamma:=\partial_{i}\partial_{j}$,
and applying $\partial^\gamma$ to \eqref{0105}$_2$, keeping in mind that
\begin{equation*} \begin{aligned}
(\nabla\times \mf{N})\times\bar{\mf{M}}=\bar{\mf{M}}\cdot \nabla \mf{N}
-\nabla (\mf{N}\cdot \bar{\mf{M}}),
\end{aligned} \end{equation*}
we find that
\begin{equation}\begin{aligned}\label{0426}
&({\varrho}+\bar{\rho})\partial^\gamma\mathbf{u}_t  +\nabla  \partial^\gamma q
+g \partial^\gamma\varrho \mf{e}_3-\mu \Delta  \partial^\gamma\mathbf{u}\\
& = -\mathbf{u}_t \partial^\gamma({\varrho}+\bar{\rho})- \partial_{i}({\varrho}+\bar{\rho})
\partial_{j}\mathbf{u}_t- \partial_{j}({\varrho}+\bar{\rho})\partial_{i}\mathbf{u}_t
-\partial^\gamma[(\varrho+\bar{\rho}){\bf u}\cdot\nabla {\bf u}] \\
 &\quad +\partial^\gamma[(\nabla\times \mf{N})\times \mf{N}]-\nabla (\partial^\gamma\mf{N}\cdot
 \bar{\mf{M}})+\bar{\mf{M}}\cdot \nabla \partial^\gamma\mf{N}.
\end{aligned}\end{equation}
Multiplying \eqref{0426} by $\partial^\beta\mf{u}$ with $|\beta|=4$ in $L^2$, we get
\begin{equation*}\begin{aligned}
&\frac{1}{2}\frac{d}{dt}\int_{\Omega}({\varrho}+\bar{\rho})|\partial^\alpha\mathbf{u}|^2\mm{d}\mf{x}
+\mu \int_{\Omega}|\nabla  \partial^\alpha\mathbf{u}|^2\mm{d}\mf{x}
\\
=&-\int_{\Omega}\partial_{1}({\varrho}+\bar{\rho})\partial^\gamma\mathbf{u}_t
\cdot \partial^\alpha\mathbf{u}\mm{d}\mf{x} +\int_{\Omega}\mathbf{u}_t
\partial^\gamma({\varrho}+\bar{\rho})\cdot \partial^\beta\mathbf{u}\mm{d}\mf{x} \\
&+\int_{\Omega}
 \partial_{i}({\varrho}+\bar{\rho})\partial_{j}\mathbf{u}_t
 \cdot\partial^\beta\mathbf{u}\mm{d}\mf{x}+ \int_{\Omega}
\partial_{j}({\varrho}+\bar{\rho})\partial_{i}\mathbf{u}_t
 \cdot\partial^\beta\mathbf{u}\mm{d}\mf{x} \\
& +\int_{\Omega} \partial^\gamma[(\varrho+\bar{\rho}){\bf u}\cdot\nabla
{\bf u}]\cdot\partial^\beta\mathbf{u}\mm{d}\mf{x} -\int_{\Omega}
\partial^\gamma[(\nabla\times \mf{N})\times \mf{N}]\cdot\partial^{\beta}\mathbf{u}\mm{d}\mf{x}\\
&-g\int_{\Omega}\partial^\alpha\varrho \partial^\alpha{u}_3\mm{d}\mf{x}
+\int_{\Omega}(\varrho+\bar{\rho})\mf{u}\cdot\nabla\partial^\alpha\mf{u}\cdot\partial^\alpha\mf{u}\mm{d}\mf{x}-\int_{\Omega}\bar{\mf{M}}\cdot \nabla \partial^\gamma\mf{N}\cdot\partial^{\beta}
\mathbf{u}\mm{d}\mf{x}=\sum_{n=1}^{9} I_n
\end{aligned}\end{equation*}
Recalling that $\mathcal{E}(t)<\delta_0\ll 1$, it is easy to verify that
\begin{equation*}\begin{aligned}
\left|\sum_{n=1}^8 I_n\right|\lesssim &
\|\partial^\alpha\varrho\|_{L^2}\|\partial^\alpha u_3\|_{L^2}+(\|\partial^\alpha
\mf{u}\|_{L^2}+\mathcal{E}^2)\|\mf{u}_t\|_{H^2}+ \mathcal{E}^3 \\
 \lesssim & C(\delta_0) \|\mf{u}\|_{L^2}^2+\delta_0\mathcal{E}^{\frac{6}{7}}\|\mf{u}_t\|_{H^2}^{\frac{8}{7}}
+\mathcal{E}^2\|\mf{u}_t\|_{H^2}+ \delta_0\mathcal{E}^2.
\end{aligned}\end{equation*}
Therefore,
\begin{equation}\begin{aligned}\label{0427}
\frac{1}{2} & \frac{d}{dt}\int_{\Omega}({\varrho}+\bar{\rho})|\partial^\alpha\mathbf{u}(t)|^2
\mm{d}\mf{x}+\mu \int_{\Omega}|\nabla  \partial^\alpha\mathbf{u}(t)|^2\mm{d}\mf{x}  \\
& \lesssim C(\delta_0) \|\mf{u}(t)\|_{L^2}^2+\delta_0\mathcal{E}^{\frac{6}{7}}(t)
\|\mf{u}_t(t)\|_{H^2}^{\frac{8}{7}} +\mathcal{E}^2(t)\|\mf{u}_t(t)\|_{H^2}+
\delta_0\mathcal{E}^2(t)+I_{9} \\
& \lesssim C(\delta_0) \|\mf{u}(t)\|_{L^2}^2+\delta_0\mathcal{E}^2(t)+I_{9},
\end{aligned}\end{equation}
where we have used \eqref{partialu} in the last inequality.

Adding \eqref{0427} to \eqref{densityenegr} and using \eqref{mageequality}, we have
 \begin{equation*}\begin{aligned}
 &\frac{d}{dt}\Big[\|(\varrho,\mf{N})(t)\|^2_{H^3}
 +\sum_{|\alpha|=3}  \|\sqrt{\varrho+\bar{\rho}}\partial^\alpha\mf{u}(t)
 \|_{L^2}^2+\sum_{i=0}^1\|\nabla \partial_t^i\mathbf{u}(t)\|^2_{L^2}+\sum_{i=0}^2
\|\sqrt{\varrho+\bar{\rho}}\partial_t^i\mf{u}(t)\|^2_{L^2}\Big]\\
&+\sum_{i=0}^1\|\partial_t^i(\nabla  \mathbf{u}, \mathbf{u}_t,\nabla\mathbf{u}_t)(t)\|^2_{L^2}
+\|\nabla^4\mf{u}(t)\|_{L^2}^2  \lesssim C(\delta_0)\|
 (\varrho,\mf{u},\mf{N})(t) \|_{L^2}^2+ \delta_0\mathcal{E}^2(t).
 \end{aligned}\end{equation*}
In particular, integrating the above inequality over $(0,t)$, we get immediately
   \begin{equation}\label{enggyn}\begin{aligned}
&\|(\varrho,\mf{N})(t)\|^2_{H^3}+\|\mathbf{u}(t)\|_{H^1}^2+\|\nabla^3\mathbf{u}(t)\|_{L^2}^2
+\|\mathbf{u}_t(t)\|_{H^1}^2+\|\mathbf{u}_{tt}(t)\|^2_{L^2}  \\
&\quad +\int_0^t\Big[\sum_{i=0}^1\|\partial_s^i(\nabla  \mathbf{u}, \mathbf{u}_s,\nabla\mf{u}_s)(s)\|^2_{L^2}
+\|\nabla^4\mf{u}(s)\|_{L^2}^2\Big]\mm{d}s  \\
& \lesssim \|(\varrho_0,\mf{N}_0)\|^2_{H^3}+\|\mathbf{u}_0\|_{H^3}^2+\|
\partial_t\mathbf{u}_0\|^2_{H^2}+\|
\partial_{t}^2\mathbf{u}_0\|^2_{L^2} \\
&\qquad +\int_0^t\Big[C(\delta_0)\|
 (\varrho,\mf{u},\mf{N})(s) \|_{L^2}^2 + \delta_0\mathcal{E}^2(s)\Big]\mm{d}s.
\end{aligned} \end{equation}

Applying the classical regularity theory on the Stokes equations to \eqref{0105}$_2$ (referring to the estimate \eqref{stokestimates}),
one finds that
 \begin{equation}\label{regularity}\begin{aligned}\|\nabla^2\mf{u}_t\|_{L^2}+\|\nabla  q_t\|_{L^2}
\lesssim  & \|\partial_t\big[(\nabla\times \mf{N})\times (\mf{N}+\bar{\mf{M}})
-g \varrho \mf{e}_3-(\varrho+\bar{\rho}){\bf u}\cdot\nabla
{\bf u}-(\varrho+\bar{\rho}){\bf u}_t\big]\|_{L^2} \\
\lesssim & (1+\mathcal{E}^2)\|\varrho_t\|_{L^2}+(\mathcal{E}+\|\varrho_t\|_{L^2})\|\mf{u}_t\|_{H^1}+(1+\mathcal{E})\|\mf{N}_t\|_{H^1}+ \|\mf{u}_{tt}\|_{L^2}\\
\lesssim & \mathcal{E}+\|\mf{u}_t\|_{H^1}+ \|\mf{u}_{tt}\|_{L^2},
\end{aligned}\end{equation}
and
 \begin{equation*}\begin{aligned}
 \|\nabla^2\mf{u}\|_{H^2}+\|\nabla  q\|_{H^2} \lesssim  & \|(\nabla\times \mf{N})\times
(\mf{N}+\bar{\mf{M}})-g \varrho \mf{e}_3-(\varrho+\bar{\rho}){\bf u}\cdot\nabla
{\bf u}-(\varrho+\bar{\rho}){\bf u}_t\|_{H^2}\\
\lesssim  &\mathcal{E}
+\|{\bf u}_t\|_{H^2}.
\end{aligned}\end{equation*}
On the other hand, letting $t\rightarrow 0$ in \eqref{partialu}, we infer that
 \begin{equation}\label{initial}\| \partial_t{\bf u}_0\|^2_{H^2} \lesssim \mathcal{E}^2_0.
 \end{equation}
 Similarly, it is easy to infer that
 \begin{equation}\label{initial123}\| \partial_t^2{\bf u}_0\|^2_{L^2} \lesssim \mathcal{E}^2_0.
 \end{equation}
Consequently, combining \eqref{enggyn} with \eqref{regularity}--\eqref{initial123}, we conclude
   \begin{equation*}\begin{aligned}
& \mathcal{E}^2(t)+\|\mathbf{u}_t(t)\|^2_{H^2}+\|\nabla q_t\|_{L^2}^2
+\|\nabla q\|_{H^2}^2+\|\mf{u}_{tt}\|_{L^2}^2 \\
& \quad +\int_0^t\left[\sum_{i=0}^1\|\partial_s^i(\nabla  \mathbf{u},\mathbf{u}_s,\nabla
\mathbf{u}_s)(s)\|^2_{L^2}+\|\nabla^4\mf{u}(s)\|_{L^2}^2\right]\mm{d}s  \\
& \lesssim \mathcal{E}^2_0 +\int_0^t\left[C(\delta_0)\|(\varrho,\mf{u},\mf{N})(s) \|_{L^2}^2
 + \delta_0\mathcal{E}^2(s)\right]\mm{d}s,\ \mbox{ provided }\delta_0\mbox{ is sufficiently small}.
\end{aligned} \end{equation*}
In addition, from \eqref{0410}, \eqref{0411}, \eqref{0416}, \eqref{0417} and \eqref{partialu} we get the following estimate, which will be used in Section \ref{sec:05}.
 \begin{equation*} \|F\|_{L^2}^2+  \|\mf{G}_t\|_{L^2}^2
+ \|\mf{H}\|_{H^1}^2\lesssim  \mathcal{E}^2+\|\mf{u}_{tt}\|^2_{L^2},
 \end{equation*}
where
 $F=-{{\mathbf{u}}}\cdot \nabla\varrho$, $\mathbf{G}:= (\nabla\times \mf{N})\times \mf{N}
 - (\varrho+\bar{\rho})\mathbf{u}\cdot\nabla \mathbf{u}-\varrho{\mathbf{u}}_t$,
and $\mf{H}:=\nabla\times (\mf{u}\times\mf{N})$.

Finally, summing up the above estimates and the existence statement, we arrive at the following conclusion:
\begin{pro} \label{pro:0401}
Let
$(\varrho,\mf{u},\mf{N})\in C^0(\bar{I}_T,H^3\times H^4\times H^3)$ be constructed in Proposition \ref{pro:0401new}, and $\delta_0'$ be chose as in \eqref{conditiondensity}.
 If $\varrho_0$ further  satisfies $\|\varrho_0\|_{H^2}\leq \delta_0'$,
 then there is a constant $\delta_0\in (0,\delta_0']$, such that
if $\mathcal{E}(t)\leq \delta_0$ on $\bar{I}_T$, then the classical solution satisfies
   \begin{equation}\label{energyinequality}\begin{aligned}
&\mathcal{E}^2(t)+\|\mathbf{u}_t(t)\|^2_{H^2}+\|(\nabla q,\mf{u}_t)_t\|_{L^2}^2+\|\nabla q\|_{H^2}^2 \\
& \qquad +\int_0^t\Big[\sum_{i=0}^1\|\partial_s^i(\nabla  \mathbf{u}, \mathbf{u}_s, \nabla\mathbf{u}_s)(s)\|^2_{L^2}
+\|\nabla^4\mf{u}(s)\|_{L^2}^2\Big]\mm{d}s  \\
&\quad \leq C(\delta_0)\mathcal{E}^2_0 +\int_0^t\left[C(\delta_0)\|
 (\varrho,\mf{u},\mf{N})(s) \|_{L^2}^2+\Lambda^*\mathcal{E}^2(s)\right]\mm{d}s,
\end{aligned} \end{equation}
and
\begin{equation}\label{energyinequality1}\begin{aligned}
\|F\|_{L^2}^2+  \|\mf{G}_t\|_{L^2}^2
+ \|\mf{H}\|_{H^1}^2\leq C(\delta_0) (\mathcal{E}^2+\|\mf{u}_{tt}\|^2_{L^2}),
\end{aligned} \end{equation}
where $\Lambda^*$ is provided by Theorem \ref{thm:0203}.
\end{pro}

\section{Proof of Theorem \ref{thm:0101}}\label{sec:05}
Now we are in a position to prove Theorem \ref{thm:0101}. First,
 in view of Theorem \ref{thm:0203}, we can construct a linear solution \begin{equation}\label{0501}
\left(\varrho^\mm{l},
{\mathbf{u}}^\mm{l},{\mathbf{N}}^\mm{l}
\right)= e^{\Lambda^* t}
\left(\bar{\varrho}_0,
\bar{\mathbf{u}}_0, \bar{\mathbf{N}}_0
\right)\mbox{ for some } \Lambda^*\in (2\Lambda/3,\Lambda],
\end{equation}
 to
\eqref{0108} with the initial data $(\bar{\varrho}_0,\bar{\mathbf{u}}_0,\bar{\mathbf{N}}_0)
\in H^\infty$ and $\mm{div}\bar{\mf{u}}_0=\mm{div}\bar{\mf{N}}_0=0$; moreover, the linear solution satisfies

\begin{equation}\label{n0502}\begin{aligned}
&0< m_0:=\left\{
                               \begin{array}{l}
                               \min\{\|{\varrho}_0\|_{L^2},\|(\bar{u}_{01},\bar{u}_{02})
\|_{L^2},\|\bar{u}_{03}\|_{L^2}\}, \\
                               \min\{\|{\varrho}_0\|_{L^2},\|({\bar{u}}_{01},{\bar{u}}_{02})
\|_{L^2},\|{\bar{u}}_{03}\|_{L^2},\|{\bar{N}}_{03}\|_{L^2}\} ,  \hbox{ if }  \bar{\bf{M}}=M\mf{e}_3,
                               \end{array}
                             \right.\end{aligned}
\end{equation}
where the sharp linear growth rate $\Lambda>0$ is defined by  \eqref{0267}, and $\bar{u}_{0i}$ and $\bar{N}_{03}$ stand for the $i$-th component of $\bar{\mf{u}}_0$ and the third component of $\bar{\mf{N}}_0$,
respectively.

Denote $C_0:=\|(\bar{\varrho}_0,\bar{\mf{u}}_0,\bar{\mf{N}}_0)\|_{\mathcal{E}}
=
\mathcal{E}(\bar{\varrho_0},\bar{\mf{u}}_0,\bar{\mf{N}}_0)$,
$(\varrho_0^\delta,\mf{u}_0^\delta,\mf{N}^\delta_0)
:=\delta (\bar{\varrho}_0,\bar{\mf{u}}_0$, $\bar{\mf{N}}_0)$, and
$C_1:=\|(\bar{\varrho_0},\bar{\mf{u}}_0,\bar{\mf{N}}_0)\|_{L^2}$.
By virtue of Proposition \ref{pro:0401}, there exists a unique local classical solution
$(\varrho^\delta,\mathbf{u}^\delta,\mathbf{N}^\delta)\in C^0(\bar{I}_T,H^3\times H^4\times H^3)$ to \eqref{0105},
emanating from the initial data $(\varrho_0^\delta,\mf{u}_0^\delta,\mf{N}^\delta_0)$
with $\|(\varrho_0^\delta,\mf{u}_0^\delta,\mf{N}^\delta_0)\|_{\mathcal{E}}=C_0\delta$.
Let us now choose $\delta_0\in (0,1)$ as small as in Proposition \ref{pro:0401} and let $C(\delta_0)>0$ be the
constant appearing in Proposition \ref{pro:0401} for the fixed choice of $\delta_0$.
Let $\iota=\min\{\delta_0,\varepsilon_0\}$, and $\delta\in (0,\iota)$ satisfy
 \begin{equation}\label{times}
 T^{\delta}=\frac{1}{\Lambda^*}\mm{ln}\frac{2\varepsilon_0}{\delta}>0,\quad\mbox{i.e. }\;
 \delta e^{\Lambda^* T^\delta}=2\varepsilon_0, \end{equation}
where the value of $\varepsilon_0$ will be fixed later to be sufficiently
small and independent of $\delta$.

We then define
 \begin{equation*}
T^*=\sup\left\{t\in (0,T^{\max})\left|~\left\|\left(\varrho^\delta,
{\mathbf{u}}^\delta,{\mathbf{N}}^\delta
\right)(t)\right\|_{\mathcal{E}}\leq {\delta_0}\right\}\right.\end{equation*}
  and
   \begin{equation*}T^{**}=\sup\left\{t\in (0,T^{\max})\left|~\left\|\left(\varrho^\delta,
{\mathbf{u}}^\delta,{\mathbf{N}}^\delta
\right)(t)\right\|_{{L}^2}\leq 2\delta C_1e^{\Lambda^* t}\right\}\right., \end{equation*}
where $T^{\mm{max}}$ denotes the maximal time of existence. Obviously, $T^*T^{**}>0$, and furthermore,
 \begin{equation}\label{0502n1}
\begin{aligned}&\left\|\left(\varrho^\delta,
{\mathbf{u}}^\delta,{\mathbf{N}}^\delta
\right)(T^*)\right\|_{\mathcal{E}}={\delta_0} \mbox{ if }T^*<\infty,\\
&\left\|\left(\varrho^\delta,
{\mathbf{u}}^\delta,{\mathbf{N}}^\delta
\right)(T^{**})\right\|_{{L}^2}=2\delta C_1e^{\Lambda^* T^{**}} \mbox{ if }T^{**}<T^{\max}.
\end{aligned}
\end{equation}
Then for all $t\leq \min\{T^\delta,T^*,T^{**}\}$,
we deduce from the estimate \eqref{energyinequality}, and the
definitions of $T^*$ and $T^{**}$ that
 \begin{equation}\begin{aligned}\label{0503}
&\left\|\left(\varrho^\delta,
{\mathbf{u}}^\delta,{\mathbf{N}}^\delta
 \right)(t)\right\|_{\mathcal{E}}^2+\|\mf{u}_{tt}^\delta\|^2_{L^2}\\
\leq &C(\delta_0) \delta^2\left\|\left( \bar{\varrho}_0,
\bar{\mathbf{u}}_0,\bar{\mathbf{N}}_0
 \right)\right\|_{\mathcal{E}}
 ^2+\Lambda^*\int_0^t \left\|\left(\varrho^\delta,
{\mathbf{u}}^\delta,{\mathbf{N}}^\delta
 \right)(s)\right\|_{\mathcal{E}}^2\mm{d}s
\\
&\ +C(\delta_0)\int_0^t\left\|\left(\varrho^\delta,
{\mathbf{u}}^\delta,{\mathbf{N}}^\delta
 \right)(s)\right\|_{L^2}^2\mm{d}s\\
\leq & {\Lambda^*}\int_0^t \left\|\left( \varrho^\delta,
{\mathbf{u}}^\delta,{\mathbf{N}}^\delta
 \right)(s)\right\|_{\mathcal{E}}^2\mm{d}s
+C_0^2C(\delta_0 )\delta^2+2 C(\delta_0 )C_1^2\delta^2e^{2\Lambda^* t}/\Lambda^*\\
\leq & \Lambda^*\int_0^t \left\|\left( \varrho^\delta,
{\mathbf{u}}^\delta,{\mathbf{N}}^\delta
\right)(s)\right\|_{\mathcal{E}}^2\mm{d}s  +
C_2\delta^2e^{2\Lambda^* t}
   \end{aligned}
 \end{equation}
for some constant $C_2:=C_2(\delta_0)>0$ independent of $\delta$.
 Then applying Gronwall's inequality to the above estimate, one obtains
  \begin{equation}\begin{aligned}\label{0504}
 \left\|\left( \varrho^\delta,
{\mathbf{u}}^\delta, {\mathbf{N}}^\delta
 \right)(t)\right\|_{\mathcal{E}}^2+\|\mf{u}_{tt}^\delta\|^2_{L^2}
\leq &C_2\delta^2e^{2\Lambda^* t}+C_2\delta^2 e^{\Lambda^* t}\int_0^t
\Lambda^* e^{\Lambda^* s}\mm{d}s\\
\leq & C_2\delta^2e^{2\Lambda^* t}+C_2\delta^2 e^{2\Lambda^* t}=
2C_2\delta^2e^{2\Lambda^* t}
   \end{aligned}
 \end{equation}
for any $t\leq \min\{T^\delta,T^*,T^{**}\}$.

Let $(\varrho^{\mathrm{d}},
{\mathbf{u}}^{\mathrm{d}},{\mathbf{N}}^{\mathrm{d}})=(\varrho^{\delta},
{\mathbf{u}}^{\delta},{\mathbf{N}}^{\delta})-\delta(\varrho^{\mathrm{l}},
{\mathbf{u}}^{\mathrm{l}},{\mathbf{N}}^{\mathrm{l}})$.
Then $(\varrho^\mm{a}_\delta,\mf{u}^{\mm{a}}_\delta,\mf{N}^{\mm{a}}_\delta):=
\delta(\varrho^{\mm{l}},\mf{u}^{\mm{l}},\mf{N}^{\mm{l}})$ is also a linear solution to
\eqref{0108} with the initial data
$(\varrho_0^\delta,\mathbf{u}_0^\delta,\mathbf{N}_0^\delta)\in H^\infty$.
Moreover, $(\varrho^{\mathrm{d}}, {\mathbf{u}}^{\mathrm{d}},{\mathbf{N}}^{\mathrm{d}})$
satisfies the following non-homogenous equations:
\begin{equation}\label{newfor1232}\left\{\begin{array}{ll}
  \varrho_t^{\mathrm{d}}+\bar{\rho}'{u}_3^{\mathrm{d}}=F^\delta, \\[1mm]
  \bar{\rho}\mathbf{u}_t^{\mathrm{d}}+\nabla q^{\mathrm{d}}+g\varrho^{\mathrm{d}} \mf{e}_3
   -\mu \Delta \mathbf{u}^{\mathrm{d}}-(\nabla\times \mf{N}^{\mathrm{d}})\times\bar{\mf{M}}=\mf{G}^\delta,
   \\[1mm]
\mf{N}_t^{\mathrm{d}}- \nabla\times (\mf{u}^{\mathrm{d}}\times \bar{\mf{M}})=\mf{H}^\delta,\\[1mm]
 \mathrm{div}\mathbf{u}^{\mathrm{d}}={0},\ \mathrm{div}\mathbf{N}^{\mathrm{d}}={0},
 \end{array}\right.\end{equation}
 where
 \begin{equation*} \displaystyle\left(\begin{array}{c}
 F^\delta\\ \mathbf{G}^\delta\\
 \mf{H}^\delta
\end{array}\right):=\left(\begin{array}{c}
                             \displaystyle     -{{\mathbf{u}}}^{\delta}\cdot \nabla\varrho^{\delta} \\
                             \displaystyle (\nabla\times \mf{N}^{\delta})\times \mf{N}^{\delta}
                             -( \varrho^{\delta}+\bar{\rho})\mathbf{u}^{\delta}\cdot\nabla
                             \mathbf{u}^{\delta}-\varrho^\delta\mathbf{u}^{\delta}_t\\[1mm]
                        \nabla\times (\mf{u}^{\delta}\times\mf{N}^{\delta})
                                \end{array}\right)
                           \end{equation*}
 with  initial conditions
 \begin{equation*} (\varrho^{\mathrm{d}}(0),
{\mathbf{u}}^{\mathrm{d}}(0),
{\mathbf{N}}^{\mathrm{d}}(0))=
{\mathbf{0}},\quad \mathrm{div}\mf{u}_0^{\mathrm{d}}=0\mbox{ and }\mathrm{div}\mf{N}_0^{\mathrm{d}}=0.
   \end{equation*}
From the estimates \eqref{energyinequality1} and \eqref{0504}, it follows that
 \begin{equation}\label{FGHestimates}
 \|F^\delta\|_{L^2}^2+  \|\mf{G}_t^\delta\|_{L^2}^2
+ \|\mf{H}^\delta\|_{H^1}^2
\leq 2C_2C(\delta_0)\delta^2e^{2\Lambda^* t}\mbox{ for any }t\leq \min\{T^\delta,T^*,T^{**}\}.
 \end{equation}
Moreover, we have the following error estimate for $\left(\varrho^{\mathrm{d}},
{\mathbf{u}}^{\mathrm{d}},{\mathbf{N}}^{\mathrm{d}}\right)$:
\begin{lem}
There exists a constant $C_3$ such that
\begin{equation}\begin{aligned}\label{0508}
 \sqrt{\left\|\left(\varrho^{\mathrm{d}},
{\mathbf{u}}^{\mathrm{d}}\right)(t)\|_{H^1}^2+\|{\mathbf{N}}^{\mathrm{d}}
  (t)\right\|_{L^2}^2}
\leq C_3(\delta e^{\Lambda^* t})^{\frac{3}{2}}.    \end{aligned}
 \end{equation}
\end{lem}
\begin{pf}
The proof is divided into two steps, in which we omit the superscript ``$\mm{d}$" in $\left(\varrho^{\mathrm{d}},
{\mathbf{u}}^{\mathrm{d}},{\mathbf{N}}^{\mathrm{d}}
 \right)$   for simplicity.

\emph{Step 1}: \emph{We use the definition of $\Lambda$ to deduce the following estimate } \begin{equation}\label{0302}\int_{\Omega} \left(g\bar{\rho}'|u_3 |^2
-M^2|\partial_i\mf{u} |^2\right)\mm{d}\mf{x}\leq {\Lambda^2}\int_{\Omega}\bar{\rho}|\mf{u} |^2
\mm{d}\mf{x}+
{\Lambda}\mu\int_{\Omega}
|\nabla\mathbf{u}|^2 \mf{x}\quad\mbox{ for }i=1,3.
\end{equation}

 Let $\hat{f}$ be the horizontal Fourier transform of $f$, i.e.,
 $$\hat{f}(\xi,x_3)=\int_{(2\pi L\mathbb{T})^2}f(\mf{x}',x_3)e^{-\mm{i}\mf{x}'\cdot\xi}
 \mm{d}\mf{x}'
,$$
  and
$$ \hat{u}_1(\xi,x_3)=-i\varphi(\xi,x_3),\;\;
\hat{u}_2(\xi,x_3)=-i\theta(\xi,x_3),\;\; \hat{u}_3(\xi,x_3)=\psi(\xi,x_3). $$
Then  $\xi_1\varphi+\xi_2\theta+\psi'=0$ because of $\mathrm{div}\mathbf{u}=0$. Moreover,
 \begin{equation*}
\widehat{\nabla \mathbf{u}}=(\widehat{\partial_i u_j})=\left(
                                                     \begin{array}{ccc}
 \xi_1\varphi &  \xi_1\theta & i\xi_1\psi \\
 \xi_2\varphi &  \xi_2\theta  & i\xi_2\psi  \\
 -i\varphi' &  -i\theta'   & \psi'  \\
 \end{array}
                                                   \right).
\end{equation*}

 By the Fubini and Parseval theorems (see \cite[Proposition 3.1.16]{grafakos2008classical}), we have (see \cite[Section 3.3]{GYTI2})
\begin{equation}\label{0305}\begin{aligned}
\int_{\Omega} g\bar{\rho}'|u_3|^2\mm{d}\mf{x}
=&\frac{1}{4\pi^2 L^2}\sum_{\xi\in (L^{-1}\mathbb{Z})^2}\int_{\mathbb{R}}g(\bar{\rho}'1_{\{\bar{\rho}'\geq 0\}}
+\bar{\rho}'1_{\{\bar{\rho}'< 0\}})|\hat{u}_3|^2\mm{d}x_3\\
= &\frac{1}{4\pi^2 L^2}\sum_{\xi\in (L^{-1}\mathbb{Z})^2}\int_{\mathbb{R}}g(\bar{\rho}'1_{\{\bar{\rho}'\geq 0\}}
+\bar{\rho}'1_{\{\bar{\rho}'< 0\}})|\psi|^2\mm{d}x_3\\
=&\frac{1}{4\pi^2 L^2}\sum_{\xi\in (L^{-1}\mathbb{Z})^2}\int_{\mathbb{R}}g\bar{\rho}'|\psi|^2\mm{d}x_3,
\end{aligned}\end{equation}
where $1_{\{\bar{\rho}'\geq 0\}}$ and $1_{\{\bar{\rho}'< 0\}}$ denote the characteristic functions.
Obviously,
\begin{equation}\label{additional0305}\begin{aligned}
\lim_{|\xi|\rightarrow 0}\int_{\mathbb{R}}g\bar{\rho}'|\psi|^2\mm{d}x_3=\int_{\mathbb{R}}g\bar{\rho}'|\psi(0,x_3)|^2\mm{d}x_3.
\end{aligned}\end{equation}
 In addition, we can control the following terms involved with $\psi$ as follows.
\begin{equation}\begin{aligned}
&\frac{1}{4\pi^2 L^2}\left[\liminf_{|\xi|\rightarrow 0}\int_{\mathbb{R}}\frac{\xi_1^2}{|\xi|^2}(|{\xi}|^2|\psi|^2
+|\psi'|^2)\mm{d}x_3+\sum_{\xi\in (L^{-1}\mathbb{Z})^2\setminus\{\mf{0}\}}\int_{\mathbb{R}}\frac{\xi_1^2}{|\xi|^2}(|{\xi}|^2|\psi|^2
+|\psi'|^2)\mm{d}x_3\right]\\
&=\frac{1}{4\pi^2 L^2}\sum_{\xi\in (L^{-1}\mathbb{Z})^2\setminus\{\mf{0}\}}\int_{\mathbb{R}}\frac{\xi_1^2}{|\xi|^2}(|\xi|^2|\psi|^2+
|\xi_1\varphi+\xi_2\theta|^2)\mm{d}x_3\\
&\leq\frac{2}{4\pi^2 L^2}\sum_{\xi\in (L^{-1}\mathbb{Z})^2}\int_{\mathbb{R}}|\xi_1|^2
(|\varphi|^2+|\theta|^2+|\psi|^2)\mm{d}x_3\\
&=\frac{2}{4\pi^2 L^2}\sum_{\xi\in (L^{-1}\mathbb{Z})^2}\sum_{j=1}^3\int_{\mathbb{R}}
|\widehat{\partial_1{u}_j}|^2\mm{d}x_3= 2\int_{\Omega} |{\partial_1 \mathbf{u}}|^2\mm{d}\mf{x},
\end{aligned}\end{equation}
\begin{equation}\begin{aligned}
&\frac{1}{4\pi^2 L^2}\left[\liminf_{|\xi|\rightarrow 0}\int_{\mathbb{R}}\frac{(|{\xi}|^2|\psi'|^2
+|\psi''|^2)}{|\xi|^2} \mm{d}x_3+
\sum_{\xi\in (L^{-1}\mathbb{Z})^2\setminus\{\mf{0}\}}\int_{\mathbb{R}}\frac{(|{\xi}|^2|\psi'|^2
+|\psi''|^2)}{|\xi|^2} \mm{d}x_3\right]\\
&\leq\frac{2}{4\pi^2 L^2}\sum_{\xi\in (L^{-1}\mathbb{Z})^2}\int_{\mathbb{R}}
(|\varphi'|^2+|\theta'|^2+|\psi'|^2) \mm{d}x_3= 2\int_{\Omega} |{\partial_3\mathbf{u}}|^2\mm{d}\mf{x},
\end{aligned}\end{equation}
\begin{equation}\begin{aligned}
& \frac{1}{4\pi^2 L^2}\left[\limsup_{|\xi|\rightarrow 0}\int_{\mathbb{R}}\frac{\bar{\rho}
(|{\xi}|^2|\psi|^2+|\psi'|^2)}{|\xi|^2}\mm{d}x_3 +\sum_{\xi\in (L^{-1}\mathbb{Z})^2\setminus\{\mf{0}\}}\int_{\mathbb{R}}\frac{\bar{\rho}
(|{\xi}|^2|\psi|^2+|\psi'|^2)}{|\xi|^2}\mm{d}x_3 \right] \\
& \leq \frac{2}{4\pi^2 L^2}\sum_{\xi\in (L^{-1}\mathbb{Z})^2}\int_{\mathbb{R}}\bar{\rho}
(|\psi|^2+|\varphi|^2+|\theta|^2)\mm{d}x_3 =2\int_{\Omega}\bar{\rho} |\mf{u}|^2
\mathrm{d}\mf{x},\end{aligned}
\end{equation}
and
\begin{equation}\label{0308}\begin{aligned}
&\frac{1}{4\pi^2 L^2}\left[\limsup_{|\xi|\rightarrow 0}\int_{\mathbb{R}}\frac{(4|{\xi}|^2|\psi'|^2
+||{\xi}|^2\psi+\psi''|^2)}{|\xi|^2}\mm{d}x_3
\right.\\
&\qquad\qquad  \left.+\sum_{\xi\in (L^{-1}\mathbb{Z})^2\setminus\{\mf{0}\}}\int_{\mathbb{R}}\frac{(4|{\xi}|^2|\psi'|^2
+||{\xi}|^2\psi+\psi''|^2)}{|\xi|^2}\mm{d}x_3\right]\\
&=\frac{1}{4\pi^2 L^2}
\left[\limsup_{|\xi|\rightarrow 0}\int_{\mathbb{R}}
\frac{2|{\xi}|^2|\psi'|^2+|{\xi}|^4|\psi|^2+|\psi''|^2}{|\xi|^2}\mm{d}x_3
\right.\\
&\qquad\qquad  \left.+\sum_{\xi\in (L^{-1}\mathbb{Z})^2\setminus\{\mf{0}\}}\int_{\mathbb{R}}
\frac{2|{\xi}|^2|\psi'|^2+|{\xi}|^4|\psi|^2+|\psi''|^2}{|\xi|^2}\mm{d}x_3\right]
\\
&\leq  \frac{2}{4\pi^2 L^2}\sum_{\xi\in (L^{-1}\mathbb{Z})^2}\int_{\mathbb{R}}
\left[|\xi|^2(|\varphi|^2+|\theta|^2+|\psi|^2)+|\varphi'|^2+|\theta'|^2+|\psi'|^2\right]\mm{d}x_3\\
& = \frac{2}{4\pi^2 L^2}\sum_{\xi\in (L^{-1}\mathbb{Z})^2} \sum_{1\leq i,j\leq 3}
\int_{\mathbb{R}} |\widehat{\partial_i {u}_j}|^2\mm{d}x_3= 2\int_{\Omega} |{\nabla\mathbf{u}}|^2\mm{d}\mf{x}.
\end{aligned}\end{equation}

Now, let
\begin{equation}Z(\psi,\xi)=
\int_{\mathbb{R}}g\bar{\rho}'|\psi|^2\mathrm{d}x_3-M^2\left\{  \begin{array}{ll}
\int_{\mathbb{R}}\xi_1^2\left(|\psi|^2+\frac{|\psi'|^2}{|\xi|^2}\right)\mm{d}x_3,
& \hbox{ for }\bar{\mf{M}}={{M}}e_1; \\[3mm]
\int_{\mathbb{R}}\left(|\psi'|^2+\frac{|\psi''|^2}{|\xi|^2}\right)\mm{d}x_3,
& \hbox{ for }\bar{\mf{M}}={M}e_3.
               \end{array}\right.
\end{equation}
By splitting $Z(\psi,\xi)=Z(\mathfrak{R}(\psi),\xi)+Z(\mathfrak{I}(\psi),\xi)$,
we may reduce the boundedness of $Z(\psi,\xi)$ to the boundedness of $Z(\mathfrak{R}(\psi),\xi)$, since
the imaginary part can be dealt with in the same manner. Hence,
without loss of generality, assume that $\psi$ is real-valued.

For ${\xi}\in \mathbb{A}^\mm{g}\cap (L^{-1}\mathbb{Z})^2\neq \emptyset$, we make use of \eqref{0210}
and \eqref{0214} with $s=\lambda({\xi})$ to deduce that
\begin{equation}\label{0311}\begin{aligned}
Z(\psi,\xi)\leq &
\frac{\lambda^2({\xi})}{|{\xi}|^2}
\int_{\mathbb{R}}\bar{\rho}(|{\xi}|^2|\psi|^2+
|\psi'|^2)\mathrm{d}x_3+\frac{\mu\lambda({\xi})}{|{\xi}|^2}
\int_{\mathbb{R}}\left(4|{\xi}|^2|\psi'|^2
+||{\xi}|^2\psi+\psi''|^2\right)\mathrm{d}x_3\\
 \leq & \frac{\Lambda}{|{\xi}|^2}\left[\Lambda\int_{\mathbb{R}}\bar{\rho}(|{\xi}|^2|\psi|^2+
|\psi'|^2)\mathrm{d}x_3 +\mu\int_{\mathbb{R}}(4|{\xi}|^2|\psi'|^2
+||{\xi}|^2\psi+\psi''|^2)\mathrm{d}x_3\right].
\end{aligned}\end{equation}
On the other hand, by Proposition \ref{pro:0203}, we infer that for any $\xi\not\in\mathbb{A}^{\mm{g}}\setminus\{\mf{0}\}$,
\begin{equation}\label{xion0312}
Z(\psi,\xi)\leq 0\leq \frac{\Lambda}{|{\xi}|^2}\left[\Lambda\int_{\mathbb{R}}\bar{\rho}(|{\xi}|^2|\psi|^2+
|\psi'|^2)\mathrm{d}x_3 +\mu\int_{\mathbb{R}}(4|{\xi}|^2|\psi'|^2
 +||{\xi}|^2\psi+\psi''|^2)\mathrm{d}x_3\right].
\end{equation}
Consequently, putting \eqref{0305}--\eqref{xion0312} together, we can conclude that \eqref{0302} holds.
In fact, taking $i=1$ for example, we have
\begin{equation*}\begin{aligned}
& 2\int_{\Omega}\left( g\bar{\rho}'|u_3|^2-M^2|\partial_1\mf{u}|^2\right)\mm{d}\mf{x}\leq  \frac{1}{4\pi^2 L^2}\sum_{\xi\in (L^{-1}\mathbb{Z})^2}\int_{\mathbb{R}}g\bar{\rho}'|\psi|^2\mm{d}x_3\\
&- \frac{M^2 }{4\pi^2 L^2}\left[\liminf_{|\xi|\rightarrow 0}\int_{\mathbb{R}}\frac{\xi_1^2}{|\xi|^2}(|{\xi}|^2|\psi|^2
+|\psi'|^2)\mm{d}x_3+\sum_{\xi\in (L^{-1}\mathbb{Z})^2\setminus\{\mf{0}\}}\int_{\mathbb{R}}\frac{\xi_1^2}{|\xi|^2}(|{\xi}|^2|\psi|^2
+|\psi'|^2)\mm{d}x_3\right]
\\&  =\frac{1}{4\pi^2 L^2}\left\{\limsup_{|\xi|\rightarrow 0}\int_{\mathbb{R}}\left[g\bar{\rho}'|\psi|^2-\frac{\xi_1^2}{|\xi|^2}(|{\xi}|^2|\psi|^2
+|\psi'|^2)\right]\mm{d}x_3+\sum_{\xi\in (L^{-1}\mathbb{Z})^2\setminus\{\mf{0}\}} Z(\psi,\xi)\right\}     \\
 &\leq \frac{1}{4\pi^2 L^2}
\left\{\limsup_{|\xi|\rightarrow 0}\frac{\Lambda}{|{\xi}|^2}\left[\Lambda\int_{\mathbb{R}}\bar{\rho}(|{\xi}|^2|\psi|^2+
|\psi'|^2)\mathrm{d}x_3 +\mu\int_{\mathbb{R}}(4|{\xi}|^2|\psi'|^2 +||{\xi}|^2\psi+\psi''|^2)\mathrm{d}x_3\right]\right.\\
 &\left.+\sum_{\xi\in (L^{-1}\mathbb{Z})^2\setminus\{\mf{0}\}}\frac{\Lambda}{|{\xi}|^2}
 \left[\Lambda\int_{\mathbb{R}}\bar{\rho}(|{\xi}|^2|\psi|^2 + |\psi'|^2)\mathrm{d}x_3 +\mu\int_{\mathbb{R}}(4|{\xi}|^2|\psi'|^2
 +||{\xi}|^2\psi+\psi''|^2)\mathrm{d}x_3\right]  \right\}   \\
 &  \leq  2{\Lambda^2} \int_{\mathbb{R}^2}\rho|\mf{u}|^2\mm{d}\mf{x}+
 2{\Lambda} \mu\int_{\Omega}|\nabla\mathbf{u}|^2\mm{d}\mf{x},\end{aligned}
\end{equation*}
which gives \eqref{0302} for $i=1$.
\vspace{1mm}

\emph{Step 2: Proof of the error estimate \eqref{0508}}

In what follows, we still denote by $C$ a generic positive constant which may depend on physical parameters
and the known functions $\bar{\varrho}_0$, $\bar{\mathbf{u}}_0$, and $\bar{\mathbf{N}}_0$.

We differentiate \eqref{newfor1232}$_2$ with $\bar{\mf{M}}={M}e_i$ ($i=1,3$) in time, multiply the resulting
equation by $\mf{u}_t$ and integrate by part over $\Omega$. Then, using the first and third equations in \eqref{newfor1232}, we obtain
\begin{equation}\label{stpe20420}
\begin{aligned}
&\frac{d}{dt}
\int_{\Omega}\left(\bar{\rho}|\mathbf{u}_t|^2+M^2|\partial_i\mf{u}|^2
-g\bar{\rho}'{u}_3^2\right)\mathrm{d}\mathbf{x}+
2\mu\int_{\Omega}|\nabla\partial_t\mathbf{u}|^2\mm{d}\mf{x}\\
&=\int_\Omega \{\mf{G}_t^\delta+[(\nabla \times \mf{H}^\delta)\times \bar{\mf{M}}]-gF^\delta \mf{e}_3\}\cdot \mf{u}_t^{\mm{d}}\mm{d}\mf{x}.
\end{aligned}\end{equation}
Here, utilizing \eqref{0501}, \eqref{0504} and \eqref{FGHestimates}, we see that
\begin{equation*}
\begin{aligned}
\left|\int_\Omega \{\mf{G}_t^\delta+[(\nabla \times \mf{H}^\delta)\times \bar{\mf{M}}]
-gF^\delta \mf{e}_3\}\cdot \mf{u}_t^{\mm{d}}\mm{d}\mf{x}\right| \leq C(\delta e^{\Lambda^* t})^{3}.
\end{aligned}\end{equation*}
Noting that $\|\mf{u}_t(0)\|\leq C\delta^3$, we integrate \eqref{stpe20420} from $0$ to $t$ and use \eqref{0302} to infer that
\begin{equation}\label{0314}
\begin{aligned}
\|\sqrt{\bar{\rho}} \mathbf{u}_t(t)\|^2_{L^2}+2\mu\int_0^t\|\nabla\partial_s\mathbf{u}(s)\|^2_{L^2}\mm{d}s
\leq {\Lambda^2} \|\sqrt{\bar{\rho}}\mf{u}(t)\|_{L^2}^2+ {\Lambda}\mu\|\nabla\mathbf{u}(t)\|^2_{L^2}+C(\delta e^{\Lambda^* t})^{3}.
\end{aligned}\end{equation}

The following inequality follows easily from integration in time and Cauchy-Schwarz's inequality.
 \begin{equation}\begin{aligned}\label{0316}
\Lambda\|\nabla\mathbf{u}(t)\|^2_{L^2}=&
2\Lambda\int_0^t
\int_{\Omega}\nabla\partial_s\mathbf{u}(s):\nabla\mathbf{u}(s)\mm{d}\mf{x}\mathrm{d}s\\
\leq &\int_0^t\|\nabla\partial_s\mathbf{u}(s)\|^2_{L^2}\mathrm{d}s
+\Lambda^2\int_0^t\|\nabla\mathbf{u}(s)\|^2_{L^2}\mathrm{d}s.
\end{aligned}\end{equation}
On the other hand,
\begin{equation}\begin{aligned}\label{0317}
\Lambda\partial_t\|\sqrt{\bar{\rho}}\mathbf{u}(t)\|^2_{L^2}=2\Lambda\int_{\Omega}
\bar{\rho}\mathbf{u}\cdot\partial_t\mathbf{u}(t)\mm{d}\mf{x}\leq\|\sqrt{\bar{\rho}}\partial_t\mathbf{u}(t)\|^2_{L^2}
+\Lambda^2\|\sqrt{\bar{\rho}}\mathbf{u}(t)\|^2_{L^2}.
\end{aligned}\end{equation}
Hence, putting \eqref{0314}--\eqref{0317} together, we obtain the differential inequality:
$$ \partial_t\|\sqrt{\bar{\rho}}\mathbf{u}(t)\|^2_{L^2}+\mu\|\nabla\mathbf{u}(t)\|^2_{L^2}
\leq 2\Lambda\left(\|\sqrt{\bar{\rho}}\mathbf{u}(t)\|^2_{L^2}+\mu\int_0^t\|\nabla\mathbf{u}(s)\|^2_{L^2}
\mathrm{d}s\right)+C(\delta e^{\Lambda^* t})^{{3}}.$$
 An application of Gronwall's inequality then implies that
\begin{equation}\begin{aligned}\label{0319}
\|\sqrt{\bar{\rho}}\mathbf{u}(t)\|^2_{L^2}+\mu\int_0^t\|\nabla\mathbf{u}(s)\|^2_{L^2}\mm{d}s
\leq C \int_0^t(\delta e^{\Lambda^* s})^{{3}}e^{2\Lambda(t-s)} \mm{d}s\leq  C(\delta e^{\Lambda^* t})^{{3}}
\end{aligned}\end{equation}
for any $t\geq 0$, where we have used the fact that $3\Lambda^*-2\Lambda>0$.
Thus, making use of \eqref{0314}, \eqref{0316} and \eqref{0319}, we deduce that
\begin{equation}\label{nablaestime1357}\begin{aligned}
\frac{1}{\Lambda}\|\sqrt{\bar{\rho}} \mathbf{u}_t(t)\|^2_{L^2}+
\mu\|\nabla\mathbf{u}(t)\|^2_{L^2}  \leq  {\Lambda}\|\sqrt{\bar{\rho}}\mathbf{u}(t)\|^2_{L^2}+2
{\Lambda}\mu\int_0^t\|\nabla\mathbf{u}(s)\|^2_{L^2}\mm{d}s  \leq C(\delta e^{\Lambda^* t})^{{3}}.
\end{aligned}\end{equation}

Now we use \eqref{0319}, \eqref{nablaestime1357} and \eqref{newfor1232}$_1$ to arrive at
\begin{equation}\begin{aligned}\label{0323}
\|\varrho(t)\|_{X}\leq
\int_0^t\|\varrho_s(s)\|_{X}\mm{d}s
\leq \|\bar{\rho}'\|_{L^\infty}\int_0^t\|u_3(s)\|_{X}\mm{d}s\leq  C (\delta e^{\Lambda^* t})^{\frac{3}{2}},
\end{aligned}\end{equation}
where $X=L^2$ or $H^1$. Similarly to the derivation of \eqref{0323}, one obtains
\begin{equation*}\begin{aligned}
\|\mf{N}(t)\|_{L^2}\leq&
\int_0^t\|\mf{N}_s(s)\|_{L^2}\mm{d}s
\leq  \|\bar{\mf{M}}\|_{L^\infty}\int_0^t\|\nabla \mf{u}(s)\|_{L^2}\mm{d}s\leq C(\delta e^{\Lambda^* t})^{\frac{3}{2}}.
\end{aligned}\end{equation*}
Thus, putting the above four estimates together, we immediately get \eqref{0508}. This completes the proof of the lemma.
\hfill $\Box$
\end{pf}

\begin{rem}\label{rem0401n} Slightly modifying the above proof,
we can easily show that any linear solution $(\varrho,\mf{u},\mf{N},q)$, enjoying
proper regularity, of \eqref{0108} with $\bar{\mf{M}}={M}e_i$ ($i=1,3$) satisfies
 the following estimates for any $t\geq 0$:
\begin{align}
& \|\varrho(t)\|_{X}^2\leq Ce^{2\Lambda
t}(\|\varrho_0\|_{X}^2+\|\mathbf{u}_0\|_{H^2}^2+\|\nabla\mf{N}_0\|_{L^2}),\quad
X=L^2\mbox{  or }H^1,\nonumber\\[1mm]
& \|\mathbf{u}(t)\|_{H^1 }^2+\|\mathbf{u}_t(t)\|^2_{L^2 }+
\int_0^t\|\nabla\mathbf{u}(s)\|^2_{L^2}\mm{d}s\leq
Ce^{2\Lambda t}(\|\varrho_0\|_{L^2}^2+\|\mathbf{u}_0\|_{H^2}^2+\|\nabla \mf{N}_0\|_{L^2}^2),\nonumber\\[1mm]
&\|\mf{N}(t)\|_{L^2}\leq Ce^{\Lambda
t}(\|\varrho_0\|_{L^2}+\|\mathbf{u}_0\|_{H^2}+\|\mf{N}_0\|_{H^1}),\nonumber
\end{align}
where the constant $C$ depends on $\bar{\rho}$, $\mu$, $\Lambda$ and $\bar{\mf{M}}$.
\end{rem}

Now, we claim that
\begin{equation}\label{n0508}
T^\delta=\min\left\{T^\delta,T^*,T^{**}\right\},
 \end{equation}
provided that small $\varepsilon_0$ is taken to be
 \begin{equation}\label{defined}
\varepsilon_0=\min\left\{\frac{{\delta_0}}
{4\sqrt{2C_2}},\frac{C_1^2}{3C_3^2}, \frac{m_0^2}{8C_3^2} \right\}.
 \end{equation}
 Indeed, if $T^*=\min\{T^{\delta},T^*,T^{**}\}$, then $T^*<\infty$; moreover, by \eqref{0504} and \eqref{times} we have
 \begin{equation*}
 \left\|\left(\varrho^\delta,
{\mathbf{u}}^\delta,{\mathbf{N}}^\delta
 \right)(T^*)\right\|_{\mathcal{E}}\leq \sqrt{2C_2}\delta e^{\Lambda^* T^*}
\leq \sqrt{2C_2}\delta e^{\Lambda^* T^\delta}=2\sqrt{2C_2}\varepsilon_0<{\delta_0},
 \end{equation*}
 which contradicts with \eqref{0502n1}$_1$. On the other hand, if $T^{**}=\min\{T^{\delta},T^*,T^{**}\}$, then $T^{**}<T^{\max}$. Moreover,
in view of \eqref{0501}, \eqref{times} and  \eqref{0508}, we have
 \begin{equation*}\begin{aligned}
 \left\|\left( \varrho^\delta,
{\mathbf{u}}^\delta,{\mathbf{N}}^\delta
 \right)(T^{**})\right\|_{L^2}
\leq  & \left\|\left( \varrho^\mm{a}_{\delta},
{\mathbf{u}}^\mm{a}_{\delta},{\mathbf{N}}^\mm{a}_{\delta}
 \right)(T^{**})\right\|_{L^2} +\left\|\left( \varrho^{\mathrm{d}},
{\mathbf{u}}^{\mathrm{d}},{\mathbf{N}}^{\mathrm{d}}
 \right)(T^{**})\right\|_{L^2} \\
\leq  &\delta \left\|\left( \varrho^\mm{l},
{\mathbf{u}}^{\mm{l}},{\mathbf{N}}^\mm{l}
\right)(T^{**})\right\|_{L^2}+C_3(\delta e^{\Lambda^* T^{**}})^{\frac{3}{2}} \\
\leq & \delta C_1e^{\Lambda^* T^{**}}+C_3(\delta e^{\Lambda^* T^{**}})^{\frac{3}{2}}
\leq \delta e^{\Lambda^* T^{**}}(C_1+C_3\sqrt{2\varepsilon_0})\\
<&2\delta C_1  e^{\Lambda^* T^{**}},
 \end{aligned} \end{equation*}
which also contradicts with \eqref{0502n1}$_2$. Therefore, \eqref{n0508} holds.

 Finally, we again use \eqref{defined}, \eqref{0508} and \eqref{n0502} to deduce that
 \begin{equation}\begin{aligned}
 \|u_3^{\delta}(T^\delta)\|_{L^2}\geq &
\|u^{\mathrm{a}}_{\delta 3}(T^{\delta})\|_{L^2}-\|u_3^{\mm{d}}(T^{\delta})\|_{L^2}
= \delta\|u_3^{\mathrm{l}}(T^{\delta})\|_{L^2}-\|u_3^{\mm{d}}(T^{\delta})\|_{L^2} \\
 \geq & m_0\delta e^{\Lambda^* T^\delta} -
C_3(\delta e^{\Lambda^* T^{\delta}})^{\frac{3}{2}}
 \geq  2m_0\varepsilon_0-C_3(2\varepsilon_0)^{\frac{3}{2}}
 \geq  m_0{\varepsilon_0},
 \end{aligned}      \end{equation}
 where $u^{\mathrm{a}}_{\delta 3}$ and $u_3^{\mm{d}}(T^{\delta})$ denote the third component of
$\mf{u}^{\mathrm{a}}_{\delta}$ and $\mf{u}^{\mm{d}}(T^{\delta})$, respectively.
Similarly, we have
$$\|\varrho(T^\delta)\|_{L^2(\Omega)},\ \|({u}_1,u_2)(T^\delta)\|_{L^2(\Omega)} \geq m_0{\varepsilon_0}.
$$   Moreover,
\begin{equation*} \|N_3(T^\delta)\|_{L^2(\Omega)}\geq m_0\varepsilon_0\quad\mbox{ if }\; \bar{\bf{M}}=M\mf{e}_3.
\end{equation*}
Finally, we take $\varepsilon=m_0\varepsilon_0$ and thus obtain Theorem \ref{thm:0101} immediately.

\section{Proof of the local well-posedness}\label{appendix}
In this section, we  adapt the standard method to briefly show the local well-posedness for the problem \eqref{0105}--\eqref{0107}
by three steps, in which we shall exploit some mathematical techniques to deal with the horizontally periodic domain.
Firstly we solve the following linearized problem for given $\mathbf{v}$:
 \begin{equation}\label{0501problem}\left\{\begin{array}{l}
\varrho_t+\mf{v}\cdot\nabla
\varrho=-\mf{v}\cdot\nabla
\bar{\rho}, \\[1mm]
(\varrho+\bar{\rho}){\bf u}_t+\nabla\tilde{q}-\mu \Delta \mathbf{u}= (\mf{N}+\bar{\mf{M}})\cdot\nabla \mf{N}
-g \varrho \mf{e}_3-(\varrho+\bar{\rho}){\bf v}\cdot\nabla \mf{v},\\[1mm]
\mf{N}_t+\mf{v}\cdot\nabla \mf{N}-\nabla \mf{v}\mf{N}=\bar{\mf{M}}\cdot\nabla \mf{v},\\[1mm]
\mathrm{div}\mathbf{u}=0,\ \mathrm{div}\mathbf{N}=0\end{array}\right.  \end{equation}
with initial and boundary conditions:
\begin{equation}\label{bounndar0502}
\begin{aligned}
& (\varrho,\mathbf u,\mathbf{N} )|_{t=0}=(\varrho_0,{\mathbf u}_0,\mathbf{N}_0)\in H^3\times H^4\times H^3,\\
& \lim_{|{x}_3|\rightarrow +\infty}(\varrho,{\bf u},\mf{N})(t,\mathbf{x}',x_3)={\bf 0}\;\mbox{ for any }t>0,
\end{aligned}\end{equation}
where $(\mf{u}_0,\mf{N}_0)$ satisfies the compatibility conditions
$\mathrm{div}\mathbf{u}_0=0$ and $\mathrm{div}\mathbf{N}_0=0$.
 Secondly, we use the technique of iteration                      
to construe a sequence of solutions $\{(\varrho^k,\mf{u}^k,\mf{N}^k)\}_{k=1}^\infty$, in which
$(\varrho^k,\mf{u}^k,\mf{N}^k)$ solves the above linearized problem with
$(\varrho^k,\mf{u}^k,\mf{N}^k)$ in place of $(\varrho,\mf{u},\mf{N})$ and $\mf{u}^{k-1}$ in place of $\mf{v}$,
and show the uniform boundedness of the solution sequence in some function space.
Finally, we further prove that the sequence of solutions is a Cauchy sequence, and thus the limit
function is the unique solution of the original problem.

Throughout the rest of this article we shall repeatedly use the abbreviations:
\begin{equation*}
\begin{aligned}
&  H^1_\sigma:=\{\mf{u}\in H^1~|~\mm{div}\mf{u}=0\},\ Q_T:=I_T\times \Omega,\\
& {V}_T:=\{\mf{v}\in C^0(\bar{I}_T,H^4)~|~
\mf{v}_t\in C^0(\bar{I}_T,H^2),\ \mf{v}_{tt}\in C^0({I}_T,L^2),\\
 &\qquad \qquad \qquad \qquad \qquad \quad \mf{v}_{t}\in L^2({I}_T,H^3),\ \mf{v}_{tt}\in L^2({I}_T,H^1),\ \mm{div}\mf{v}=0\},\\
& {H}_T:=\{(\varrho,\mf{N})\in C^0(\bar{I}_T,H^3)~|~(\varrho_t,\mf{N}_t)\in C^0(\bar{I}_T,H^2),\ (\varrho_{tt},\mf{N}_{tt})\in C^0(\bar{I}_T,H^1)\},\\
& {F}_T:=\{\mf{f}\in C^0(\bar{I}_T,H^2)~|~\mf{f}_{t}\in
C^0(\bar{I}_T,L^2)\cap L^2({I}_T,H^1),\ \mf{f}_{tt}\in
L^2({I}_T,L^2)  \}.
\end{aligned}\end{equation*}

\subsection{Unique solvability  of linearized problems}
To show the solvability of the linearized problem \eqref{0501problem}--\eqref{bounndar0502}, it suffices to solve the following two problems
 \begin{equation}\label{hypobolic}\left\{\begin{array}{l}
\varrho_t+\mf{v}\cdot\nabla
\varrho=-\mf{v}\cdot\nabla
\bar{\rho}, \\[1mm]
\mf{N}_t+\mf{v}\cdot\nabla \mf{N}-\nabla \mf{v}\mf{N}=\bar{\mf{M}}\cdot\nabla \mf{v},\\
\mathrm{div}\mathbf{N}=0\end{array}\right.
\end{equation}
with initial and boundary conditions
\begin{equation}\label{intbouadnry51}
\begin{aligned}
& (\varrho,\mathbf{N} )|_{t=0}=(\varrho_0, \mathbf{N}_0)\in H^3\times H^3 ,\\
& \lim_{|{x}_3|\rightarrow +\infty}(\varrho,\mf{N})(t,\mathbf{x}',x_3)={\bf 0}\;\;\mbox{ for any }t>0,
\end{aligned}\end{equation}
and
\begin{equation}\label{loc0503jjw}
\left\{    \begin{array}{ll}
(\varrho+\bar{\rho}){\bf u}_t+\nabla \tilde{q}-\mu \Delta \mathbf{u}=\mf{f}, \\
\mathrm{div}\mathbf{u}=0
  \end{array}  \right.   \end{equation}
  with initial and boundary conditions
\begin{equation}\label{intabioudr502}
\mathbf{u}|_{t=0}={\mathbf u}_0\in H^4,\ \lim_{|{x}_3|\rightarrow +\infty}\mf{u}(t,\mathbf{x}',x_3)={\bf 0}\
\mbox{ for any }t>0,
\end{equation}
where  $\mf{f}:=(\mf{N}+\bar{\mf{M}})\cdot\nabla \mf{N}-g \varrho \mf{e}_3- (\varrho+\bar{\rho}){\bf v}\cdot\nabla\mf{v}$.

By the standard hyperbolic theory, we have the following result on the existence of a unique solution to the first problem \eqref{hypobolic}--\eqref{intbouadnry51}:
\begin{lem}\label{lem:0501}Let $\mf{v}\in {V}_T$, $(\varrho_0, \mf{N}_0)\in \mf{H}^3\times
\mf{H}^3$ and $\mm{div}\mf{N}_0=0$, then, for any $T>0$, the problem \eqref{hypobolic}--\eqref{intbouadnry51} possesses a unique classical solution $(\varrho,\mf{N})\in {H}_T$ emanating from the initial data $(\varrho_0, \mf{N}_0)$. Moreover, $(\varrho,\mf{N})$
satisfies
the following identities:
\begin{align}
&\label{fracrho}\frac{1}{2}\frac{d}{dt}\|\partial^{\alpha}\varrho\|_{L^2}^2
=-\int_\Omega[\partial^\alpha(\mf{v}\cdot\nabla \bar{\rho})+\partial^{\alpha}( \mf{v}\cdot\nabla \varrho)]\partial^\alpha\varrho\mm{d}\mf{x},\\
&\label{fracrho12}\frac{1}{2}\frac{d}{dt}\|\partial^{\alpha}\mf{N}\|_{L^2}^2=\int_\Omega[\partial^\alpha(\bar{\mf{M}}\cdot
\nabla\mf{v}+\nabla\mf{v}\mf{N})-\partial^{\alpha}(\mf{v}\cdot \nabla\mf{N}) ]\partial^\alpha\mf{N}\mm{d}\mf{x},
\end{align}
where we have defined that
$$\mf{v}\cdot \partial^{\alpha}\nabla\varrho=0\mbox{ and }\mf{v}\cdot \partial^{\alpha}\nabla\mf{N}=0\;\;\mbox{ for any }|\alpha|=3.$$
\end{lem}
\begin{pf}
Following the proof of \cite[Theorem 2.16]{NASII04}, we can easily check that
\cite[Theorem 2.16]{NASII04} still holds with the horizontally periodic domain $\Omega$ in place of $\mathbb{R}^3$.
Thus, we immediately get the unique solvability of the problem \eqref{hypobolic}--\eqref{intbouadnry51}.

The strong continuity \eqref{fracrho} and \eqref{fracrho12} can be shown by adapting the idea in the proof of \cite[Lemm 6.9]{NASII04}.
Here we give the proof \eqref{fracrho} for the reader's convenience. Without a generality, we assume that
$\alpha=(3,0,0)$.
Similarly in \cite[Lemm 6.9]{NASII04}, we can deduce from \eqref{hypobolic} that
\begin{equation}\begin{aligned}&\partial_t S_\varepsilon^2(\partial^3_{x_1}\varrho)+\mm{div}(S_\varepsilon^2(\partial^3_{x_1}\varrho)\mf{v})-2r_\varepsilon S_\varepsilon(\partial^3_{x_1}\varrho)\\
&= -2 S_\varepsilon(\partial^3_{x_1}(\mf{v}\cdot\nabla \bar{\rho})+\partial_{x_1}^2(\partial_{x_1}\mf{v}\cdot\nabla \varrho)+\partial_{x_1}^2\mf{v}\cdot \partial_{x_1}\nabla \varrho+2\partial_{x_1}\mf{v}\cdot \partial_{x_1}^2\nabla\varrho )S_\varepsilon(\partial^3_{x_1}\varrho)
\end{aligned}
\end{equation}
a.e. in $Q_T$, where $S_\varepsilon$ is a standard mollifier,
 and $r_\varepsilon =\mm{div}(S_\varepsilon(\partial^3_{x_1}\varrho)\mf{v})-\mm{div}(S_\varepsilon(\partial^3_{x_1}\varrho\mf{v}))\to 0$
 in $L^2(I_T,K)$ for any bounded domain $K\subset \Omega$. Note that $S_\varepsilon(f)\in H^1\cap C^\infty$ if $f\in H^1$.
 Let $\psi\in C_0^\infty(I_T)$ and
  \begin{equation}\label{phconstctin}\begin{aligned}
&\phi_n(x_3)=1\mbox{ for }|x_3|<n,\ \phi_n(x_3)=0\mbox{ for }|x_3|>2n,\\
&0\leq \phi_n(x_3)\leq 1,\ |\phi_n'(x_3)|\leq c/n\mbox{ for
any }x_3\in (-2n,-n)\cup(n,2n),
\end{aligned}  \end{equation}
where the constant $c$ is independent of $n$.
Multiplying \eqref{hypobolic}$_1$ by $\phi_n\psi\in C_0^\infty(I_T)$, and integrating the resulting equality, we get
$$\begin{aligned}&\int_0^T\psi_t\int_\Omega S_\varepsilon^2(\partial^3_{x_1}\varrho)\phi_n\mm{d}\mf{x}\mm{d}t
+\int_0^T\psi\int_\Omega S_\varepsilon^2(\partial^3_{x_1}\varrho)v_3\phi_n'\mm{d}\mf{x}\mm{d}t
+2\int_0^T\psi\int_\Omega r_\varepsilon S_\varepsilon(\partial^3_{x_1}\varrho)\phi_n\mm{d}\mf{x}\mm{d}t\\
&=2\int_0^T\psi\int_\Omega S_\varepsilon(\partial^3_{x_1}(\mf{v}\cdot\nabla \bar{\rho})+\partial_{x_1}^2(\partial_{x_1}\mf{v}\cdot\nabla \varrho)
+\partial_{x_1}^2\mf{v}\cdot \partial_{x_1}\nabla \varrho+2\partial_{x_1}\mf{v}\cdot \partial_{x_1}^2\nabla\varrho )S_\varepsilon(\partial^3_{x_1}\varrho)\phi_n\mm{d}\mf{x}\mm{d}t.
\end{aligned}$$
Thus, taking $\varepsilon\to 0$, and then $n\rightarrow +\infty$, we immediately get \eqref{fracrho}.
\hfill$\Box$
\end{pf}

We turn to the existence proof of a unique solution to the second problem \eqref{loc0503jjw}--\eqref{intabioudr502} by employing the Galerkin method, the
domain expansion technique and a technique of improving regularity.
\begin{lem}\label{lem:0502}
 Let $T>0$, $\mf{u}_0\in H^4$, $\mf{v}\in  {V}_T$, $\mf{f}\in  {F}_T$, $\inf_{x\in \Omega}\{(\varrho_0+\bar{\rho})(\mf{x})\}>0$ and
 $\varrho$ be constructed in Lemma \ref{lem:0501}, then there exists a unique classical solution $\mathbf{u} \in
 {V}_T$ to the problem \eqref{loc0503jjw}--\eqref{intabioudr502} with an associated pressure $\tilde{q}$ satisfying
\begin{equation}\label{pressureesimta}
\tilde{q}\in L^\infty(I_T,L^6),\;\;\nabla\tilde{q}\in C^0(\bar{I}_T,H^2),\;\; (\nabla \tilde{q})_t\in C^0(\bar{I}_T,L^2).
\end{equation}
\end{lem}
\begin{pf}
 We divide the proof of Lemma \ref{lem:0502} into three steps.

(1) First we solve \eqref{loc0503jjw}--\eqref{intabioudr502} without pressure by the Galerkin method, and the existence of solutions to the Galerkin approximate problem is established
by adapting the basic idea from the proof of \cite[Theorem 4.3]{GYTIPDE}.
Since $H^2_{\sigma}:=\{\mf{v}\in H^2~|~\mm{div}\mf{v}=0\}$ is separable, it possesses a countable basis $\{\mf{w}^j\}_{j=1}^\infty$. For each $m\geq 1$ we define a finite dimensional space $\mathcal{U}_m:=\mm{span}\{\mf{w}^1,\ldots,\mf{w}^m\}\subset H^2_\sigma$ and an approximate solution
$$\mf{u}^m(t)=a_j^m(t)\mf{w}_j\;\mbox{ with }\; a_j^m:\ I_T\rightarrow\mathbb{R}\;\mbox{ for }\; j=1,\ldots, m,$$
where as usual we have used the Einstein convention of summation of the repeated index $j$. Clearly, for each $\mf{v}\in H^2_\sigma$,
we have $\mathcal{P}_m\mf{v}\rightarrow \mf{v}$ as $m\rightarrow \infty$ and $\|\mathcal{P}_m\mf{u}\|_{H^2}\leq \|\mf{u}\|_{H^2}$,
where $\mathcal{P}_m:~H_\sigma^2\rightarrow \mathcal{U}_m$ is the $H^2_\sigma$ the orthogonal projection onto $\mathcal{U}_m$.
In addition, by the method of characteristics, we can deduce from the mass equation $\eqref{hypobolic}_1$ that
$$\tilde{\rho}:=\sup_{\mf{x}\in \Omega}\{(\varrho_0+\bar{\rho})(\mf{x})\}\geq (\varrho+\bar{\rho})(t,\mf{x})\geq \underline{\rho}:=
\inf_{\mf{x}\in \Omega}\{(\varrho_0+\bar{\rho})(\mf{x})\}>0\;\; \mbox{ for any }(t,\mf{x})\in Q_T.$$

By the theory of ODEs, we can find the coefficients $a_j^m\in C^{2}(\bar{I}_T)$, such that
\begin{equation}\label{appweaksolu}
\int_{\Omega}(\varrho+\bar{\rho})\mf{u}^m_t\cdot\mf{w}\mm{d}\mf{x}+\mu\int_{\Omega}
\nabla \mf{u}^m:\nabla \mf{w}\mm{d}\mf{x}=\int_\Omega \mf{f}\cdot\mf{w}\mm{d}\mf{x}
\end{equation}
with initial data
$$ \mf{u}^m(0)=\mathcal{P}^m\mf{u}_0\in \mathcal{U}_m $$
for each $\mf{w}\in \mathcal{U}_m$.
Next, we derive uniform estimates for $\mf{u}_m$. In what follows, we denote by $C(\cdots)$ a generic positive constant depending only on its variables.
The notation $a\lesssim b$ means that $a\leq \tilde{C}b$ for a universal constant $\tilde{C}>0$,
which may depend on $T$, $\bar{\mf{M}}$, $\mu$, $g$, $\underline{\rho}$, $\tilde{\rho}$, $\bar{\rho}$,
and the norms $\|(\varrho_0,\mf{N}_0)\|_{H^3}$, $\|\mf{u}_0\|_{H^4}$, $\|\varrho_t(0)\|_{L^\infty}$, $\|\mf{f}(0)\|_{H^2}$,
and $\|\mf{f}_t(0)\|_{L^2}$. Moreover, $\tilde{C}$ is increasing with $T$, and the norms $\|\varrho_t(0)\|_{L^\infty}$, $\|\mf{f}(0)\|_{H^2}$ and $\|\mf{f}_t(0)\|_{L^2}$,
if $\tilde{C}$ depends on them.

Taking $\mf{w}=\mf{u}^m$ in \eqref{appweaksolu}, we see that
\begin{equation}\label{baseestimate}
\frac{1}{2}\frac{\mm{d}}{\mm{d}t}\int_{\Omega}(\varrho+\bar{\rho})|\mf{u}^m|^2\mm{d}\mf{x}
+\mu\int_{\Omega} |\nabla \mf{u}^m|^2\mm{d}\mf{x}=\int_\Omega \mf{f}\cdot\mf{u}^m\mm{d}\mf{x}
+\frac{1}{2}\int_{\Omega}\varrho_t|\mf{u}^m|^2\mm{d}\mf{x},
\end{equation}
which yields \begin{equation*}
\frac{\mm{d}}{\mm{d}t}\|\sqrt{\varrho+\bar{\rho}}\mf{u}^m\|_{L^2}\leq2 \left\|\frac{\mf{f}}{\sqrt{\varrho+\bar{\rho}}}\right\|_{L^2}
+\left\|\frac{\varrho_t}{\varrho+\bar{\rho}}\right\|_{L^\infty}\|\sqrt{\varrho+\bar{\rho}}\mf{u}^m\|_{L^2}.
\end{equation*}
Applying Gronwall's inequality, we get
$$ \|\mf{u}^m \|_{L^2}^2\lesssim e^{\frac{2T}{\underline{\rho}}\|\varrho_t\|_{L^\infty(Q_T)}}
\left(1 + {T}\|\mf{f}\|_{L^2(Q_T)}^2\right),$$
which, together with \eqref{baseestimate}, implies that
\begin{equation} \label{lowestimesv}
\|\mf{u}^m(t)\|_{L^2}^2+\int_0^t\|\nabla \mf{u}^m(s)\|_{L^2}^2\mm{d}s
\lesssim e^{\frac{3T}{\underline{\rho}}\|\varrho_t\|_{L^\infty(Q_T)}} \left( 1 + T\|\mf{f}\|_{L^2(Q_T)}^2\right).
\end{equation}
We differentiate \eqref{appweaksolu} in time and take $\mf{w}=\mf{u}^m_t$ to deduce
\begin{equation*}
\begin{aligned}
\frac{1}{2}\frac{\mm{d}}{\mm{d}t}\int_{\Omega}(\varrho+\bar{\rho})|\mf{u}^m_t|^2\mm{d}\mf{x}
+\mu\int_{\Omega}
|\nabla \mf{u}^m_t|^2\mm{d}\mf{x}=\int_\Omega \mf{f}_t\cdot\mf{u}^m_t\mm{d}\mf{x}
+\frac{1}{2}\int_{\Omega}\varrho_t|\mf{u}^m_t|^2\mm{d}\mf{x},
\end{aligned}
\end{equation*}
which yields
\begin{equation*}
\begin{aligned}
\|\mf{u}^m_t(t)\|^2_{L^2} +\int_{0}^t \|\nabla \mf{u}^m_t(s)\|^2_{L^2}\mm{d}s
\lesssim e^{\frac{3T}{\underline{\rho}}\|\varrho_t\|_{L^\infty(Q_T)}}\left(1+T\|\mf{f}_t \|_{L^2(Q_T)}^2\right).
\end{aligned}\end{equation*}
Finally, taking $\mf{w}=\mf{u}_t^m$ in \eqref{appweaksolu}, we find that
\begin{equation}\label{lowesti2}
\begin{aligned}
\| \nabla \mf{u}^m(t)\|_{L^2}^2 +\int_0^t\|\mf{u}_t^m(s)\|_{L^2}^2\mm{d}s\lesssim 1+\|\mf{f}\|_{L^2(Q_T)}^2.
\end{aligned}
\end{equation}

Summing up the previous estimates and making use of the estimate
 $$\|\mf{f}(t)\|_{L^2}=\left\|\mf{f}(0)+\int_0^t\mf{f}_{s}(s)\mm{d}t\right\|_{L^2}\lesssim
1+\sqrt{T}\|\mf{f}_t\|_{L^2(Q_T)},$$ we arrive at
\begin{equation}\label{sumestimeat}
\begin{aligned}
&\|   \mf{u}^m(t)\|_{H^1}^2+\|   \mf{u}^m_t(t)\|_{L^2}^2+\int_0^t
(\|\nabla \mf{u}^m(s)\|_{L^2}^2+\|\mf{u}_t^m(s)\|_{H^1}^2)\mm{d}s\\
&\lesssim e^{\frac{3T}{\underline{\rho}}\|\varrho_t\|_{L^\infty(Q_T)}}\left(1+
T\|\mf{f}_t\|_{L^2(Q_T)}^2\right).
\end{aligned}\end{equation}

(2) Now we can obtain a strong solution from the above uniform estimates for the approximate solutions.
In view of \eqref{sumestimeat}, up to the extraction of a subsequence, one has
$$\begin{aligned}
&\mf{u}^m\rightarrow \mf{u}\mbox{ weakly-* in }L^\infty(I_T,H^1_\sigma),\ \mf{u}^m_t\rightarrow \mf{u}_t\mbox{ weakly-* in }L^\infty(I_T,L^2),\\
&\nabla \mf{u}^m\rightarrow \nabla\mf{u}\mbox{ weakly  in }L^2(I_T,L^2),\ \mf{u}^m_t\rightarrow \mf{u}_t\mbox{ weakly  in }L^2(I_T,H^1),\\
&\mf{u}^m\rightarrow \mf{u}\mbox{ strongly in }C^0(\bar{I_T},L^2(K))\mbox{ for any bounded domain }K\subset \Omega,\
\mf{u}(0)=\mf{u}_0,
\end{aligned}$$
and $\mm{div}\mf{u}_t=0$.
Thus, if we take limit in \eqref{appweaksolu}, we get
\begin{equation}\label{0507n1}
\int_{\Omega}(\varrho+\bar{\rho})\mf{u}_t\cdot\mf{w}\mm{d}\mf{x}+\mu\int_{\Omega}
\nabla\mf{u}:\nabla\mf{w}\mm{d}\mf{x} =\int_\Omega\mf{f}\cdot\mf{w}\mm{d}\mf{x}\;\;\mbox{ for any }\mf{w}\in H^2_\sigma,
\end{equation}
which can be regarded as a weak solution of the following Stokes equations
 \begin{equation*}
-\mu \Delta \mathbf{u}+ \nabla \tilde{q}=\mf{F}:=\mf{f}-(\varrho+\bar{\rho}){\bf u}_t\mbox{ in }\Omega, \end{equation*}
where $\mf{F}\in L^\infty(I_T,L^2)\cap L^2(I_T,H^1)$.

Next we proceed to derive more estimates for $\mf{u}$ in the Stokes equations by the domain expansion technique
(i.e., the classical regularity theory on the Stokes equations).
Let $L_n=(-n,n)$ and $\Omega_n=(2\pi L\mathbb{T})^2\times L_n$. Similarly to \cite[Propositions 2.9 and 3.7]{GYTIPDE}, we can show that the Stokes problem
 \begin{equation}\label{appstokes}
\left\{
  \begin{array}{ll}
-\mu \Delta \mathbf{u}_n+\nabla \tilde{q}_n=\mf{F} \mbox{ in }\Omega_n, \\
\mm{div}\mf{u}_n=0, \\
   \mf{u}|_{ x_3=-2\pi nL}=\mf{u}|_{x_3=2\pi nL}=0
  \end{array}    \right.
 \end{equation}
 admits a unique strong solution $\mf{u}_n $ with a unique associated pressure $\tilde{q}_n$, such that
 \begin{equation*}\begin{aligned}
& \|\nabla^k \mf{u}_n\|_{H^2(\Omega_n)}^2 + \|\nabla^{k+1}\tilde{q}_n\|_{L^2(\Omega_n)}^2
\leq C(\Omega_n,\mu) \|\nabla^k \mf{F}\|_{L^2(\Omega_n)}^2\;\;\mbox{ for }k=0,1,\nonumber\\
& \|\tilde{q}_n\|_{L^2(\Omega_n)}^2 \leq C(\Omega_n) \|\mf{F}\|_{L^2(\Omega_n)}^2\;\mbox{ and }\;
\int_{\Omega_n}\tilde{q}_n\mm{d}\mf{x}=0.
\end{aligned}\end{equation*}
By scaling the spatial variables on a cuboid domain $(0,2\pi nL)^2\times L_n$ and using the horizontally periodic property
(i.e., $\mf{u}_n(x_1,x_2,x_3)=\mf{u}_n(x_1+2\pi L,x_2+2\pi L,x_3)$), we can deduce the following uniform estimates
\begin{align*}&\| \nabla^{2+k} \mf{u}_n\|_{L^2(\Omega_n)}^2+\|  \nabla^{1+k} \tilde{q}_n\|_{L^2(\Omega_{n_0})}^2\leq C(\mu) \|\nabla^k \mf{F}\|_{L^{2}}^2
\quad\mbox{ for }k=0,1,\nonumber\\
&  \|\tilde{q}_n\|_{L^2(\Omega_{n})}^2\leq  C(\mu)n^2\|\mf{F}\|_{L^2}^2\;\mbox{ and }\;\|\tilde{q}_n\|_{L^6(\Omega_{n})}^2\leq  C(\mu)\|\mf{F}\|_{L^2}^2.
\end{align*}
Moreover, using H\"older's inequality, we get
$$\|\tilde{q}_n\|_{L^2(\Omega_{n_0})}^2\leq  C(\mu)n_0^{2/3}\|\mf{F}\|_{L^2}^2\quad\mbox{ for any }n\geq n_0.$$

We make use of \eqref{0507n1} and \eqref{appstokes} to deduce
$$ \mu\int_{\Omega_n}
|\nabla \mf{u}_n|^2\mm{d}\mf{x}=\int_{\Omega_n} \mf{F}\cdot\mf{u}_n\mm{d}\mf{x} =\mu\int_{\Omega_n}
\nabla \mf{u}:\nabla \mf{u}_n\mm{d}\mf{x},  $$
whence,
$$ \|\nabla \mf{u}_n\|_{L^2(\Omega_n)}\leq \|\nabla \mf{u}\|_{L^2}. $$
 Noting that
 $$\| \mf{u}_n\|_{L^{6}(\Omega_n)}\leq C(\Omega_n)\| \nabla \mf{u}_n\|_{L^{2}(\Omega_n)},$$
we may scale the spatial variables on a cuboid domain $(0,2\pi nL)^2\times L_n$  and utilize the horizontally periodic property to get
$$ \| \mf{u}_n\|_{L^{6}(\Omega_n)}\leq c\|\nabla \mf{u}\|_{L^2}\quad\mbox{ for some constant }c. $$

Now, recalling
$\mf{F}\in L^\infty(I_T,L^2)\cap L^2(I_T,H^1)$ and using the established uniform estimates for $(\mf{u}_n,\tilde{q}_n)$,
up to the extraction of a subsequence, one has
$$\begin{aligned}&\nabla^k\mf{u}_n\rightarrow\nabla^k \bar{\mf{u}}\mbox{ weakly-* in }L^\infty(I_T,L^2)\mbox{ for }k=1,2,\;\;
\nabla^3\mf{u}_n\rightarrow \nabla^3 \bar{\mf{u}}\mbox{ weakly in }L^2(I_T,L^2),\\
& \nabla\tilde{q}_n\rightarrow\nabla\tilde{q}\mbox{ weakly-* in }L^\infty(I_T,L^2),\;\; \nabla^2  \tilde{q}_n\rightarrow \nabla^2  \tilde{q} \mbox{ weakly in }L^2(I_T,L^2),\\
&\mf{u}_n\rightarrow\bar{\mf{u}}\mbox{ weakly-* in }L^\infty(I_T, L^6),\;\;\tilde{q}_n\rightarrow \tilde{q}\mbox{ weakly-* in } L^\infty(I_T,L^6),
\end{aligned} $$
where we have extended $\nabla^k \mf{u}_n$ and $\nabla^m \tilde{q}_n$ by zero extension outside $\Omega_n$.
Thus, we have the following classical regularity estimates on the Stokes problem:
 \begin{align}&\label{stokestimates}
\| \nabla^{2} \mf{u}\|_{L^\infty(I_T,L^2)}^2+\|  \nabla \tilde{q}\|_{L^\infty(I_T,L^2)}^2+
\|\tilde{q}\|_{L^\infty(I_T,L^6)}^2\leq C(\mu) \|\mf{F}\|_{L^\infty(I_T,L^{2})}^2,\\
& \| \nabla^{3} \mf{u}\|_{L^2(I_T,L^2)}^2+\|  \nabla^2 \tilde{q}\|_{L^2(I_T,L^2)}^2\leq C(\mu) \|\nabla \mf{F}\|_{L^2(Q_T)}^2,\nonumber\\
&\label{qomesgq}\|  \tilde{q}\|_{L^\infty(I_T,L^2(\Omega_{n_0}))} \leq C(\mu) n_0^{\frac{1}{3}} \|\mf{F}\|_{L^\infty(I_T,L^2)}\quad\mbox{ for any }n_0>0.
\end{align}
Moreover,\begin{align}
&\label{stoequeiona}
-\mu \Delta \bar{\mathbf{u}}+ \nabla \tilde{q}=\mf{F}, \quad \mm{div}\bar{\mf{u}}=0, \\
&\label{0514n}\mu\int_\Omega \nabla \bar{\mathbf{u}}:\nabla \mf{w}\mm{d}\mf{x} =\int_\Omega\mf{F}\cdot
\mf{w}\mm{d}\mf{x}
 \end{align}
 for any $$\begin{aligned}
\mf{w}\in C_{0,\sigma}^\infty:=\{\mf{w}\in C^\infty~|~&\mm{div}\mf{w}=0\mbox{ and }\exists~R>0,\mbox{ such that }\\
&\mf{w}(x',x_3)=0\mbox{ for any }|x_3|>R\}.
\end{aligned}$$ It should be noted that $C_{0,\sigma}^\infty$
 is dense in $H_\sigma^1$, which will be proved in Lemma \ref{lem:0503}. By a density argument, \eqref{0514n} holds for any $\mf{w}\in H^1_\sigma$,
 which, together with \eqref{0507n1}, implies that $ \mf{u}=\bar{\mf{u}}$.

Summing up the above estimates, we obtain the following inequality:
\begin{equation*}
\begin{aligned}
&\|   \mf{u}(t)\|_{H^2}^2+\|   \mf{u}_t(t)\|_{L^2}^2+\|  \nabla \tilde{q}(t)\|_{L^2}^2+
\int_0^t
\left(\| \mf{u}(s)\|_{H^3}^2+\|\mf{u}_t(s)\|_{H^1}^2+\|\nabla^2\tilde{q}\|_{L^2}^2\right)\mm{d}s\\
&\lesssim e^{\frac{3T}{\underline{\rho}}\|\varrho_t\|_{L^\infty(Q_T)}}\left(1+
\|\varrho\|_{H^2}^2\right)
\left(1+T\|\mf{f}_t\|_{L^2(Q_T)}^2+\|\nabla \mf{f}
\|_{L^2(Q_T)}^2\right).
\end{aligned}\end{equation*}

(3) Before improving the regularity of $\mf{u}$, we show the continuity of $(\mf{u},\nabla\tilde{q})$.
First, similarly to \cite[Remark 6]{YCHKimUnique}, we can prove that for a.e. $t\in I_T$,
\begin{equation}\label{eq0510}
\begin{aligned}
\frac{1}{2}\frac{\mathrm{d}}{\mathrm{d}t}\int_{\Omega}(\varrho+\bar{\rho})|\mf{u}_t|^2
\mm{d}\mf{x}=&\int_{\Omega}\mf{f}_t\cdot{\mf{u}}_t\mm{d}\mf{x}+\frac{1}{2}\int_\Omega \varrho_t|\mf{u}_t|^2\mm{d}\mf{x}-\mu\int_{\Omega}|\nabla \mathbf{u}_t
 |^2\mm{d}\mf{x}.\end{aligned}\end{equation}
It follows from \eqref{0507n1} that for any $\varphi\in H_\sigma^1$,
\begin{equation}\label{absolution}
\begin{aligned}
 \frac{\mathrm{d}}{\mathrm{d}t}\int_{\Omega }(\varrho+\bar{\rho})\mf{u}_t\cdot\varphi\mm{d}\mf{x}
=\int_{\Omega}\mf{f}_t\cdot\varphi\mm{d}\mf{x}-\mu \int_{\Omega} \nabla \mathbf{u}_t:\nabla \varphi\mm{d}\mf{x}\mbox{ holds for a.e. } t\in I_T.
\end{aligned} \end{equation}
From the regularity of $(\varrho,\mf{u})$ we get
\begin{equation}\label{fracg-1}((\varrho+\bar{\rho})\mf{u}_t)_t\in L^2((0,T),H_\sigma^{-1})\mbox{ and }
\frac{\mm{d}}{\mm{d}t}\int_{\Omega }(\varrho+\bar{\rho})\mf{u}_t\cdot\varphi
\mm{d}\mf{x}=<((\varrho+\bar{\rho})\mf{u}_t)_t,\varphi>,
\end{equation}
where $H^{-1}_\sigma$ denotes the dual space of
$H_\sigma^1$ and $<\cdot,\cdot>$ the corresponding dual product.
Therefore, the identity \eqref{eq0510} follows immediately from \eqref{absolution},
\eqref{fracg-1} and the following identity
$$\frac{\mm{d}}{\mm{d}t}\int_{\Omega }(\varrho+\bar{\rho})|\mf{u}_t|^2
\mm{d}\mf{x}=2<((\varrho+\bar{\rho})\mf{u}_t)_t,\mf{u}_t>-\int_\Omega \varrho_t|\mf{u}_t|^2\mm{d}\mf{x},$$
which can be easily established by means of a classical regularization method.

With \eqref{eq0510} in hand, we can infer that
\begin{equation}\label{continutiofu}\mf{u}_t\in C^0(\bar{I}_T,L^2).
\end{equation}
In fact, it follows from \eqref{absolution} and \eqref{eq0510} that $(\varrho+\bar{\rho})\mf{u}_t\in C^0(\bar{I}_T, H^{-1}_\sigma)$ and $\|\sqrt{\varrho+\bar{\rho}}\mf{u}_t\|_{L^2}^2\in C^0(\bar{I}_T)$.
Since $(\varrho+\bar{\rho})\mf{u}_t\in L^\infty(I_T, L^2)$ and $L^2\hookrightarrow H^{-1}_\sigma$, it follows from the embedding theorem
that $(\varrho+\bar{\rho})\mf{u}_t\in C(\bar{I}_T, L^2_{\mm{weak}})$. 
Thus, recalling $\varrho+\bar{\rho}\in C^0(\bar{I}_T,H^{2})$, we have $\sqrt{(\varrho+\bar{\rho})}\mf{u}_t\in C^0(\bar{I}_T,L^2)$
and \eqref{continutiofu}. 
Noting that $\mf{u}\in C^0(\bar{I}_T,H^2)$ by virtue of the embedding theorem, we see that $\nabla \tilde{q}\in  C^0(\bar{I}_T,L^2)$.
In view of \eqref{stoequeiona}, we find that
 \begin{equation}\Delta\tilde{q}=\frac{\nabla (\varrho+\bar{\rho})\cdot\nabla \tilde{q}}{\varrho+\bar{\rho}}-\frac{\mu\nabla (\varrho+\bar{\rho})\cdot\Delta\mf{u}}{\varrho+\bar{\rho}}+(\varrho+\bar{\rho})\mm{div}\left(\frac{\mf{f}}{\varrho+\bar{\rho}}\right)
\in C^0(\bar{I}_T,L^2).  \end{equation}
Hence, by virtue of \eqref{qomesgq}, the elliptic estimates, the periodic property and a cut-off technique,
there exists a constant $c$, such that
\begin{equation*}  \begin{aligned}
\|\nabla^2\tilde{q}(t)\|_{L^2}\leq c\left(\|\Delta\tilde{q}(t)\|_{L^2}+\|\nabla\tilde{q}(t)\|_{L^2}+
\|\mf{F}(t)\|_{L^2}\right)\quad \mbox{ for a.e. }t\in I_T,
\end{aligned}\end{equation*}
which implies that
\begin{equation}\label{eqiatonqegH}    \nabla^2 \tilde{q}\in L^\infty({I},L^2),\end{equation}
whence, $\nabla\tilde{q}|_{t=0}\in H^1$.

 Now, we improve the regularity of $\mf{u}$. Consider the following problem:
 \begin{equation}\label{loc0503jjwtimes}
\left\{
  \begin{array}{ll}
(\varrho+\bar{\rho}){\bf v}_t + \nabla \tilde{p}-\mu \Delta \mathbf{v}= \mf{g}\quad\mbox{ in }\Omega, \\
\mm{div}\mf{v}=0,\\
   \mf{v}|_{t=0}= \mf{v}_0:= [(\mf{f}-\nabla \tilde{q}+\mu \Delta \mathbf{u})/(\varrho +\bar{\rho})]|_{t=0},
  \end{array}  \right.
 \end{equation}
where $\mf{g}:=\mf{f}_t-\varrho_t\mf{u}_t$. Keeping in mind that
$\mf{g}\in C^0(\bar{I}_T,{L^2})$ and $\mf{v}_0\in H^1_\sigma$, referring to the construction
 of the solution $\mf{u}$ and the derivation of regularity for $\mf{u}$, we can easily obtain a strong solution
$\mf{v}$ to the problem \eqref{loc0503jjwtimes} with an associated pressure $\tilde{p}$,
where $\mf{v}\in C^0(\bar{I}_T,H^1)\cap L^2(I_T,H^2),\ \mf{v}_t\in L^2(I_T,L^2)$, $\nabla \tilde{p}\in L^2(I_T,L^2)$ and
\begin{equation*}
\begin{aligned}
&\|   \mf{v}(t)\|_{H^1}^2+\int_0^t
(\|\nabla \mf{v}(s)\|_{H^1}^2+\|\mf{v}_t(s)\|_{L^2}^2)\mm{d}s\\
&\lesssim e^{\frac{6T}{\underline{\rho}}\|\varrho_t\|_{L^\infty(Q_T)}} \left(1+T\|\varrho_t\|^2_{L^\infty(Q_T)}\right) \left(1+\|\mf{f}_t
\|_{L^2(Q_T)}^2\right).
\end{aligned}\end{equation*}
Moreover, for any $\varphi\in H_\sigma^1$, \begin{equation}\label{equan0520}\frac{d}{dt}\int_{\Omega}(\varrho+\bar{\rho})\mf{v}\cdot\varphi\mm{d}\mf{x}+\mu\int_{\Omega}
\nabla \mf{v}:\nabla \varphi\mm{d}\mf{x}=\int_\Omega (\mf{g}+\varrho_t\mf{v})\cdot\varphi\mm{d}\mf{x}.
\end{equation}
Denoting $\tilde{\mf{v}}:=\mf{v}-\mf{u}_t$, and recalling
$\tilde{\mf{v}}|_{t=0}=0$ and $\sqrt{(\varrho+\bar{\rho})}\tilde{\mf{v}}\in C(\bar{I}_T,L^2)$, we
subtract \eqref{equan0520} from \eqref{absolution} to infer that
$$\frac{d}{dt} \int_{\Omega}(\varrho+\bar{\rho})\tilde{\mf{v}}\cdot\varphi\mm{d}\mf{x} + \mu\int_{\Omega}
\nabla \tilde{\mf{v}}:\nabla \varphi\mm{d}\mf{x}\mm{d}t=\int_\Omega \varrho_t\tilde{\mf{v}}\cdot\varphi\mm{d}\mf{x}.$$
 Similarly to the derivation of \eqref{eq0510}, we can establish
$$\frac{d}{dt}\int_{\Omega}(\varrho+\bar{\rho})|\tilde{\mf{v}}|^2\mm{d}\mf{x}\mm{d}t+\mu \int_{\Omega}
|\nabla \tilde{\mf{v}}|^2 \mm{d}\mf{x} =\int_\Omega \varrho_t|\tilde{\mf{v}}|^2\mm{d}\mf{x},$$
which implies $\mf{v}=\mf{u}_t$. Therefore, $\mf{u}_t\in C^0(\bar{I}_T,H^1)$, $\mf{u}_{tt}\in L^2(I_T,L^2)$ and
$\nabla \tilde{p}=(\nabla \tilde{q})_t$.

Using the regularity estimate for the Stokes equations, we further have $\mf{u}\in L^\infty(\bar{I}_T,H^3)$.
Hence, $\Delta\tilde{q}\in L^\infty(\bar{I}_T,H^1)$. In view of the elliptic interior regularity, we see that
$\nabla \tilde{q}\in L^\infty({I},H^2_{\mm{loc}}(\mathbb{R}^3))$. Thus, similarly to the derivation of
\eqref{eqiatonqegH}, we have $\nabla^3 \tilde{q}\in L^\infty({I},L^2)$, which yields
$\nabla\tilde{q}|_{t=0}\in H^2$. Consequently, $\mf{v}_0\in H^2$.
Recalling that $\mf{g}_t\in L^2(I_T,{L^2})$, we easily see that $\mf{u}_t\in C^0(\bar{I}_T, H^2)$ solves
the problem \eqref{loc0503jjwtimes} and enjoys the following estimate
\begin{equation*}
\begin{aligned}
&\|   \mf{u}_t(t)\|_{H^2}^2+\|   \mf{u}_{tt}(t)\|_{L^2}^2+\int_0^t
\left(\| \mf{u}_t(s)\|_{H^3}^2+\|\mf{u}_{tt}(s)\|_{H^1}^2\right)\mm{d}s\\
&\lesssim e^{\frac{9T}{\underline{\rho}}\|\varrho_t\|_{L^\infty(Q_T)}} (1+\|\varrho\|_{H^2}^2)
\left(1+T\|\varrho_{t}\|_{L^\infty(I_T,H^2)}^2+T\|\varrho_{tt}\|^2_{L^\infty(I_T,H^1)}\right)\\
&\qquad\times \left(1+T\|\varrho_t\|^2_{L^\infty(Q_T)}
\right) \left(1+T\|\mf{f}_{tt}\|_{L^2(Q_T)}^2+\|\nabla \mf{f}_t\|_{L^2(Q_T)}^2\right);
\end{aligned}\end{equation*}
and moreover, $\nabla\tilde{p}\in C^0(\bar{I}_T,L^2)\cap L^2(I_T, H^1)$.
Employing the regularity estimate on the Stokes equations again, we arrive at
\begin{equation}\label{estimatesofu}
\begin{aligned}
&\|   \mf{u}(t)\|_{H^4}^2+\|   \mf{u}_t(t)\|_{H^2}^2+\|   \mf{u}_{tt}(t)\|_{L^2}^2+\int_0^t
\left(\| \mf{u}_t(s)\|_{H^3}^2+\|\mf{u}_{tt}(s)\|_{H^1}^2\right)\mm{d}s\\
&\leq \tilde{C} e^{\frac{9T}{\underline{\rho}}\|\varrho_t\|_{L^\infty(Q_T)}} (1+\|\varrho\|_{H^2}^2)^2
\left(1+T\|\varrho_{t}\|_{L^\infty(I_T,H^2)}^2+T\|\varrho_{tt}\|^2_{L^\infty(I_T,H^1)}\right)\\
&\quad\times \left(1+T\|\varrho_t\|^2_{L^\infty(Q_T)}
\right) \left(1+T\|\mf{f}_{tt}\|_{L^2(Q_T)}^2+\|\nabla \mf{f}_t\|_{L^2(Q_T)}^2+\|\mf{f}\|_{L^\infty(I_T,H^2)}^2\right);
\end{aligned}\end{equation}
and moreover,
$\mf{u}\in C^0(\bar{I}_T, H^4)$ and $\nabla \tilde{q}\in C^0(\bar{I}_T,H^2)$.
Summing up the above estimates, we obtain the desired conclusions immediately.\hfill$\Box$
\end{pf}

Finally, we briefly prove that $C_{0,\sigma}^\infty$ is dense in $H_\sigma^1$.
\begin{lem}\label{lem:0503}
 $C_{0,\sigma}^\infty$ is dense in $H_\sigma^1$.
\end{lem}
 \begin{pf} Let $\mf{u}\in H_\sigma^1$ and $n>0$. We define $\Omega_n:=(2\pi L\mathbb{T}^2 )\times
 (-n,n)$, $\Omega_{n,2n}:=\Omega_{2n}/\Omega_{n}$,
 $\tilde{\Omega}_n:=(0,2\pi L)^2\times (-2n,2n)$ and
$\phi_n(x_3)\in C_0^\infty(\mathbb{R})$ be constructed as in \eqref{phconstctin}.
In view of \cite[Theorem III.3.4]{GPGAI}, for each given $n$, there exists a function
$\mf{w}_n\in H^{1}(\Omega_{2n})\cap H_0^1(\tilde{\Omega}_n)$, satisfying
$$\begin{aligned}&\mm{div}\mf{w}_n=-u_3 \phi_n',\\
& \|\nabla \mf{w}_n\|_{L^2(\Omega_{2n})}\leq c\|u_3 \phi_n'\|_{L^2(\Omega_{2n})},
\end{aligned}$$
where and throughout this proof, $c$ denotes a constant independent of the domain. Thus,
$$\|\mf{w}_n\|_{L^2(\Omega_{2n})}\leq cn^{\frac{1}{3}}\|\mf{w}_n\|_{L^6}
\leq cn^{\frac{1}{3}}\|\nabla \mf{w}_n\|_{L^2(\Omega_{2n})}\leq c\|{u}_3\|_{L^2(\Omega_{n,2n})}.$$
We define $\mf{w}_n=0$ on $\mathbb{R}^3\setminus{\Omega}_{2n}$. Let $\mf{v}_n=\phi_n\mf{u}+\mf{w}_n$ and $S_\varepsilon$
be a standard mollifier, then $ S_\varepsilon(\mf{v}_{n})\in C_{0,\sigma}^\infty$. Moreover,
\begin{equation*}
\begin{aligned}\|\mf{u}-S_\varepsilon(\mf{v}_n)\|_{L^2}\leq & \|\mf{u}-\mf{v}_n\|_{L^2}+\|\mf{v}_n-S_\varepsilon(\mf{v}_n)\|_{L^2} \\
\leq& \|(1-\phi_n)\mf{u}\|_{L^2}+\|\mf{w}_n\|_{L^2}+\|\mf{v}_n-S_\varepsilon(\mf{v}_n)\|_{L^2},
\end{aligned}\end{equation*}
which implies that
$$ S_\varepsilon(\mf{v}_n)\rightarrow \mf{u}\mbox{ as }\varepsilon\rightarrow 0\mbox{ and }n\rightarrow +\infty.$$
This shows that $C_{0,\sigma}^\infty$ is dense in $H_\sigma^1$.
\hfill$\Box$\end{pf}

\subsection{Uniform estimates}

Let $(\varrho_0,\mf{u}_0,\mf{N}_0)$ satisfy the assumption in Proposition \ref{pro:0401new}, then, in view of Lemmas \ref{lem:0501} and \ref{lem:0502},
for any given $T>0$ and any given $\mf{v}\in V_T$ satisfying $\mf{v}(0)=\mf{u}_0$,
the linearized problem \eqref{0501problem}--\eqref{bounndar0502} possesses a unique classical solution $((\varrho,\mf{N}),\mf{u})\in H_T\times V_T$.
To emphasize such relation between $\mf{u}$ and $\mf{v}$, we can define $\mf{u}:=S(\mf{v})$.
In addition, we can deduce from \eqref{0501problem}$_2$ that
\begin{equation}\label{partialut}
\|{\bf u}_t(0)\|_{H^1}\leq C_0, \end{equation}
where the positive constant $C_0$ depends on the initial data $(\varrho_0,\mf{u}_0,\mf{N}_0)$
 and other physical parameters in the perturbed equations, and is independent of $\mf{v}$.

Now, let us denote $ {D}_{\kappa}^T=\{\mf{v}\in V_T~|~\|\mf{v}\|_{V_T}\leq \kappa,\ \|\mf{v}_t(0)\|_{H^1}\leq C_0$, $\mf{v}(0)=\mf{u}_0\}$ and
$$\|\mf{v}\|_{V_T}:=\sqrt{\|\mf{v}\|_{L^\infty( {I}_T,H^4)}^2+\|\mf{v}_t\|_{L^\infty( {I}_T,H^2)}^2+\|\mf{v}_{tt}\|_{L^\infty( {I}_T,L^2)}^2+\|\mf{v}_{t}\|_{L^2( {I}_T,H^3)}^2
+\|\mf{v}_{tt}\|_{L^2({I}_T,H^1)}^2}.$$
We show next that there exists two positive constants $T\in (0,1)$ and $\kappa\geq 1$ depending on
$(\varrho_0,\mf{u}_0,\mf{N}_0)$ and other physical parameters in the perturbed equations, such that
\begin{equation}\label{uniforestimate}
 \kappa^8 T= 1\mbox{ and }\|S(\mf{v})\|_{V_{T}}\leq \kappa \quad \mbox{ for any } \mf{v}\in D^T_{\kappa}.
\end{equation}
Throughout this and next subsections, the notation $a\lesssim b$ means that $a\leq Cb$ for some constant $C>0$ which may depend on
$(\varrho_0,\mf{u}_0,\mf{N}_0)$ and other physical parameters in the perturbed equations.

First we temporarily suppose that $\kappa^8 T=1$, $\kappa\geq 1$ and $T\in (0, 1)$. Similarly to the derivation of \eqref{density},
we can deduce from \eqref{0501problem}$_1$ and \eqref{0501problem}$_2$ that
 \begin{equation*}\label{} \begin{aligned}
 \frac{d}{dt}\|(\varrho,\mf{N}) \|^2_{H^3}\lesssim
\|\mf{v}\|_{H^4}(1+\|(\varrho ,\mf{N})\|_{H^3})\|(\varrho, \mf{N}) \|_{H^3},
\end{aligned} \end{equation*}
which yields
\begin{equation}\label{engyrofdensiy}\|(\varrho ,\mf{N})\|_{L^\infty(I_T,H^3)}\lesssim e^{\kappa T}(1+\kappa T)\leq C.
\end{equation}
Furthermore, we can get from \eqref{0501problem}$_1$ and \eqref{0501problem}$_2$ that
\begin{align}
& \|\varrho_t(0)\|_{L^\infty}\leq C,\\
& \label{Nutestimae}\|(\varrho,\mf{N})_t\|_{L^\infty(Q_T)}  \lesssim \|(\varrho,\mf{N})_t\|_{L^\infty(I_T,H^2)}  \lesssim \kappa,\\
&\|(\varrho,\mf{N})_{tt}\|_{L^\infty(I_T,H^1)}\lesssim \kappa+\kappa^2.
\end{align}
Letting $\mf{f}:=(\mf{N}+\bar{\mf{M}})\cdot\nabla \mf{N}-g \varrho \mf{e}_3-(\varrho+\bar{\rho}){\bf v}\cdot\nabla\mf{v}$,
we see that
\begin{align}
&\|\mf{f}(0)\|_{H^2}^2+\|\mf{f}_t(0)\|_{L^2}^2\leq C,\\
& \|\mf{f}_{tt}\|_{L^2(Q_T)} \lesssim \kappa+\kappa^2+\kappa^3+\kappa^4,\\
&\label{} \|\nabla \mf{f}_{t}\|_{L^\infty(I_T,L^2)} \lesssim \kappa+\kappa^2+\kappa^3,\\
&\label{sevestemsnew} \|\mf{f}\|_{L^\infty(I_T,H^2)} \lesssim 1+\kappa^2 T.
\end{align}

Recalling that $\mf{u}$ satisfies \eqref{estimatesofu} and $\tilde{C}$ is increasing in $\|\varrho_t(0)\|_{L^\infty}$, $\|\mf{f}(0)\|_{H^2}$ and $\|\mf{f}_t(0)\|_{L^2}$,
we substitute \eqref{engyrofdensiy}--\eqref{sevestemsnew} into \eqref{estimatesofu} to arrive at
\begin{equation}\label{finestim524}
\begin{aligned}\|   S(\mf{v})\|_{V_T}\leq  C.
\end{aligned}\end{equation}

Finally, if we take $\kappa=\max\{C,1\}$ and $T:=\kappa^{-6}$ where the positive constant $C$ is from
\eqref{finestim524}, then we obtain \eqref{uniforestimate}.
\subsection{Taking the limit}
Let $(\varrho_0,\mf{u}_0,\mf{N}_0)$ satisfy the assumption in  Proposition \ref{pro:0401new}
and $\mf{u}^0\equiv\mf{u}_0$.
In view of Lemmas \ref{lem:0501} and \ref{lem:0502}, we can construct a function sequence $\{(\varrho^k,\mf{u}^k,\mf{N}^k,\tilde{q}^{k})\}_{k=1}^\infty$ satisfying
$((\varrho^k,\mf{N}^k),\mf{u}^k)\in H_T\times V_T$ and
 \begin{equation}\label{iteratiequat}\left\{\begin{array}{l}
\varrho_t^k+\mf{u}^{k-1}\cdot\nabla
\varrho^k=-\mf{u}^{k-1}\cdot\nabla
\bar{\rho}, \\[1mm]
(\varrho^{k}+\bar{\rho}){\bf u}_t^{k}+\nabla\tilde{q}^{k}-\mu \Delta \mathbf{u}^{k}= (\mf{N}^{k}+\bar{\mf{M}})\cdot\nabla \mf{N}^{k}-g \varrho^{k} \mf{e}_3-(\varrho^{k}+\bar{\rho}){\mf{ u}^{k-1}}\cdot\nabla
\mf{u}^{k-1},\\[1mm]
\mf{N}_t^{k}+\mf{u}^{k-1}\cdot\nabla \mf{N}^{k}-\nabla \mf{u}^{k-1}\mf{N}^{k}=
\bar{\mf{M}}\cdot\nabla \mf{u}^{k-1},\\[1mm]
\mathrm{div}\mathbf{u}^{k}=0,\ \mathrm{div}\mathbf{N}^{k}=0\end{array}\right.  \end{equation}
with initial data:
\begin{equation*}
(\varrho^k,\mathbf{u}^k,\mathbf{N}^k )|_{t=0}=(\varrho_0,{\mathbf u}_0,\mathbf{N}_0).
\end{equation*}
Moreover, by virtue of the uniform estimates \eqref{partialut}--\eqref{engyrofdensiy}
and \eqref{Nutestimae}, there exits $T\in (0,1)$ such that the solution sequence $\{(\varrho^k,\mf{u}^k,\mf{N}^k)\}_{k=1}^\infty$ satisfies the following uniform estimate
\begin{equation}\label{n0531} \|(\varrho^k,\mf{N}^k)\|_{L^\infty(I_T,H^3)}+
\|(\varrho^k,\mf{N}^k)_t\|_{L^\infty(I_T,H^2)}+\|\mf{u}^k\|_{{V}_T}\leq C\quad\mbox{for any }k\geq 1,\end{equation}
which, together with the  regularity estimates on the Stokes equations, implies that
\begin{equation}\label{n0531presure}
\|\nabla{\tilde{q}}^k\|_{L^\infty(I_T,H^2)}+\|\tilde{q}^k\|_{L^\infty(I_T,L^6)} \leq C\quad\mbox{for any }k\geq 1.\end{equation}
In addition,  by the mass equation \eqref{iteratiequat}$_1$,
\begin{equation*}
\label{}\inf_{\mathbf{x}\in\Omega}\{(\varrho_0+\bar{\rho})(\mathbf{x})\}
\leq \varrho^k(t,\mf{x})+\bar{\rho}(\mf{x})\leq \sup_{\mathbf{x}\in\Omega}\{(\varrho_0+\bar{\rho})({\mathbf{x}})\}\quad\mbox{for any }k\geq 1.
\end{equation*}

In order to take the limit in \eqref{iteratiequat}, we shall further show $\{\mf{u}_k\}_{k=1}^\infty$
is a Cauchy sequence. To this end, we define
$$(\bar{\varrho}^{k+1},\bar{\mf{u}}^{k+1},\bar{\mf{N}}^{k+1},\nabla \bar{q}^{k+1})=
({\varrho}^{k+1}-{\varrho}^{k},\mf{u}^{k+1}-\mf{u}^k,{\mf{N}}^{k+1}-{\mf{N}}^{k},\nabla ( \tilde{q}^{k+1}-\tilde{q}^{k})),$$
which satisfies
 \begin{equation}\label{difeequion}\left\{\begin{array}{l}
\bar{\varrho}_t^{k+1}+\mf{u}^{k}\cdot\nabla
\bar{\varrho}^{k+1}=-\bar{\mf{u}}^{k}\cdot\nabla
\bar{\rho}-\bar{\mf{u}}^k\cdot\nabla \varrho^k, \\[1mm]
(\varrho^{k+1}+\bar{\rho})\bar{\mf{ u}}_t^{k+1}+\nabla\bar{q}^{k+1}-\mu \Delta \bar{\mathbf{u}}^{k+1}= (\mf{N}^{k+1}
+\bar{\mf{M}})\cdot\nabla \bar{\mf{N}}^{k+1}-g \bar{\varrho}^{k+1} \mf{e}_3-\bar{\varrho}^{k+1}{\mf{ u}}_t^{k}\\
-(\varrho^{k+1}+\bar{\rho}){\mf{ u}^{k}}\cdot\nabla
\bar{\mf{u}}^{k}-\bar{\mf{N}}^{k+1}\cdot\nabla {\mf{N}}^{k}-(\varrho^{k}+\bar{\rho}){\bar{\mf{ u}}^{k}}\cdot\nabla
\mf{u}^{k-1}-\bar{\varrho}^{k+1}{{\mf{ u}}^{k-1}}\cdot\nabla
\mf{u}^{k-1}:=\mf{f}_k,\\[1mm]
\bar{\mf{N}}_t^{k+1}+\mf{u}^{k}\cdot\nabla \bar{\mf{N}}^{k+1}-\nabla \mf{u}^{k}\bar{\mf{N}}^{k+1}=
\bar{\mf{M}}\cdot\nabla \bar{\mf{u}}^{k}-\bar{\mf{u}}^{k}\cdot\nabla {\mf{N}}^{k}+\nabla \bar{\mf{u}}^{k}{\mf{N}}^{k},\\[1mm]
\mathrm{div}\mathbf{u}^{k+1}=0,\ \mathrm{div}\mathbf{N}^{k+1}=0\end{array}\right.  \end{equation}
with initial data
 \begin{equation*}\label{}(\bar{\varrho}^{k+1},\bar{\mf{u}}^{k+1}, \bar{\mf{N}}^{k+1})=(0,\mf{0},\mf{0})\  \mbox{in}\ \Omega.
\end{equation*}

Recalling that $\kappa^8 T = 1$ and $T\leq 1$,  we easily find by \eqref{difeequion}$_1$ and \eqref{difeequion}$_2$ that
\begin{align}
&\label{rhoNestimes}\sup_{t\in I_T}\|(\bar{\varrho}^{k+1},\bar{\mf{N}}^{k+1})\|_{H^1}\lesssim    T\sup_{t\in I_T}\|\bar{\mf{u}}^k\|_{H^2}, \\
&\sup_{t\in I_T}\|(\bar{\varrho}^{k+1},\bar{\mf{N}}^{k+1})_t\|_{L^2}\lesssim \sup_{t\in I_T}\|\bar{\mf{u}}^k\|_{H^2}.\nonumber \end{align}
Similarly to the derivation of \eqref{lowestimesv} and \eqref{lowesti2}, we get from \eqref{difeequion}$_2$ that
\begin{equation}
\|\bar{\mf{u}}^{k+1}(t)
\|_{L^2}^2+\int_0^t\|\nabla \bar{\mf{u}}^{k+1}(s)\|_{L^2}^2\mm{d}s
\lesssim\|\mf{f}_k\|_{L^2(Q_T)}^2
\end{equation}
and
\begin{equation}   \begin{aligned}
\| \nabla \bar{\mf{u}}^{k+1}(t)\|_{L^2}^2+\int_0^t \|\bar{\mf{u}}_s^{k+1}(s)\|_{L^2}^2\mm{d}s\lesssim \|\mf{f}_k\|_{L^2(Q_T)}^2\mbox{ for any }t\in {I}_t.
\end{aligned}\end{equation}

Next we continue to derive a bound for $\nabla^2 \bar{\mf{u}}^{k+1}$. To this end, we differentiate \eqref{difeequion}$_2$ in $t$,
multiply the resulting equations by $\bar{\mf{u}}_{t}^{k+1}$ in $L^2$ to obtain
\begin{equation*}  \begin{aligned}
&\frac{1}{2}\frac{\mm{d}}{\mm{d}t}\int_{\Omega}(\varrho^{k+1}+\bar{\rho})|\bar{\mf{u}}^{k+1}_t|^2\mm{d}\mf{x}
+\mu\int_{\Omega} |\nabla \bar{\mf{u}}^{k+1}_t|^2\mm{d}\mf{x}\\
&=\frac{1}{2}\int_{\Omega}\varrho_t^{k+1}|\bar{\mf{u}}^{k+1}_t|^2\mm{d}\mf{x}
+\int_\Omega [\mf{N}^{k+1}_t\cdot\nabla \bar{\mf{N}}^{k+1}-\varrho^{k+1}_t{\mf{ u}^{k}}\cdot\nabla
\bar{\mf{u}}^{k}-(\varrho^{k+1}+\bar{\rho})\mf{ u}^{k}_t\cdot\nabla
\bar{\mf{u}}^{k}]\cdot{\bar{\mf{u}}^{k+1}_t}\mm{d}\mf{x}   \\
&\quad -\int_\Omega \left\{(\mf{N}^{k+1}+\bar{\mf{M}})\cdot\nabla \bar{\mf{u}}^{k+1}_{t}\cdot\bar{\mf{N}}^{k+1}_t
-[(\varrho^{k+1}+\bar{\rho}){\mf{ u}^{k}}\cdot\nabla\bar{\mf{u}}^{k+1}_{t}+{\mf{ u}^{k}}\cdot\nabla (\varrho^{k+1}+\bar{\rho})
\bar{\mf{u}}^{k+1}_{t}]\cdot\bar{\mf{u}}^{k}_t\right\}\mm{d}\mf{x}  \\
&\quad -\int_\Omega [g \bar{\varrho}^{k+1} \mf{e}_3+\bar{\varrho}^{k+1}{\mf{ u}}_t^{k}
+\bar{\mf{N}}^{k+1}\cdot\nabla {\mf{N}}^{k}+(\varrho^{k}+\bar{\rho}){\bar{\mf{ u}}^{k}}\cdot\nabla
\mf{u}^{k-1} +\bar{\varrho}^{k+1}{{\mf{ u}}^{k-1}}\cdot\nabla \mf{u}^{k-1}_t]_t\cdot\bar{\mf{u}}^{k+1}_t\mm{d}\mf{x},
\end{aligned}
\end{equation*}
where the terms on the right hand side can be bounded from above by
$$C[\|\sqrt{(\varrho^{k+1}+\bar{\rho})}\bar{\mf{u}}_t^{k+1}\|_{L^2}^2
+\sup_{t\in I_T}(\|\bar{\mf{u}}^k\|_{H^2}^2+\|\bar{\mf{u}}^{k}_t\|_{L^2}^2)]+\frac{\mu}{2}\|\nabla
\bar{\mf{u}}_t^{k+1}\|_{L^2}^2.$$
Therefore, we have
$$
\frac{\mm{d}}{\mm{d}t}\|\sqrt{(\varrho^{k+1}+\bar{\rho})}\bar{\mf{u}}^{k+1}_t\|^2_{L^2}+\|\nabla \bar{\mf{u}}_t^{k+1}\|_{L^2}^2
\lesssim \|\sqrt{(\varrho^{k+1}+\bar{\rho})}\bar{\mf{u}}_t^{k+1}\|_{L^2}^2
+\sup_{t\in I_T}(\|\bar{\mf{u}}^k\|_{H^2}^2+\|\bar{\mf{u}}^{k}_t\|_{L^2}^2).    $$
Applying Gronwall's lemma to the above inequality and recalling $\bar{\mf{u}}_t^{k+1}(0)=0$, we conclude
 $$ \|\bar{\mf{u}}^{k+1}_t(t)\|_{L^2}^2 + \int_0^t\|\nabla \bar{\mf{u}}_s^{k+1}(s)\|_{L^2}^2\mm{d}s
 \lesssim  T\sup_{t\in I_T}(\|\bar{\mf{u}}^k\|_{H^2}^2 + \|\bar{\mf{u}}^{k}_t\|_{L^2}^2)\mbox{ for any }t\in I_T,$$
With the help of the classical regularity estimate on the Stokes problem, we obtain
\begin{align}&\begin{aligned}
\|\nabla^2 \bar{\mf{u}}^{k+1}\|_{L^2}^2 +\|\nabla\bar{q}^{k+1}\|_{L^2}^2\lesssim \|\bar{\mf{u}}^{k+1}_t\|_{L^2}^2+\|\mf{f}_k\|_{L^2}^2,
\end{aligned}
\end{align}
where $\|\mf{f}_k\|_{L^2}$ can be bounded as follows.
\begin{equation}\label{fkestimates}\|\mf{f}_k\|_{L^2}^2\lesssim T\sup_{t\in I_T}(\|\bar{\mf{u}}^k\|_{H^2}^2+
\|\bar{\mf{u}}^{k}_t\|_{L^2}^2)+T\|\nabla \bar{\mf{u}}_t^{k}\|_{L^2(Q_T)}^2.
\end{equation}

Summing up the estimates \eqref{rhoNestimes}--\eqref{fkestimates}, we arrive at
\begin{equation*}\label{}
\begin{aligned}
&\sup_{t\in I_T}(\|(\bar{\varrho}^{k+1},\bar{\mf{N}}^{k+1})\|_{H^1}^2+\|\nabla\bar{q}^{k+1}\|_{L^2}^2+\|\bar{\mf{u}}^{k+1} \|_{H^2}^2+\|\bar{\mf{u}}^{k+1}_t
\|_{L^2}^2)+\|\nabla \bar{\mf{u}}^{k+1}_t\|_{L^2(Q_T)}^2  \\
& \quad \leq \bar{C}T\left[\sup_{t\in I_T}(\|\bar{\mf{u}}^k\|_{H^2}^2 + \|\bar{\mf{u}}^{k}_t\|_{L^2}^2)+\|\nabla \bar{\mf{u}}_t^{k}\|_{L^2(Q_T)}^2\right]\\
&\quad \leq \frac{1}{2}\left[\sup_{t\in I_T}(\|\bar{\mf{u}}^k\|_{H^2}^2+\|\bar{\mf{u}}^k_t\|_{L^2}^2)+\|\nabla \bar{\mf{u}}_t^{k}\|_{L^2(Q_T)}^2\right]\qquad\mbox{ for any }k\geq 1,
\end{aligned}     \end{equation*}
provided $T$ is chosen so small that $T\leq ({2\bar{C}})^{-1}$, where the constant $\bar{C}$ may depend on $(\varrho_0,\mf{u}_0,\mf{N}_0)$
and other physical parameters in the perturbed equations.
The above inequality implies that
$$\sum_{k=1}^\infty\left(\|(\bar{\varrho}^{k+1},\bar{\mf{N}}^{k+1})\|_{L^\infty(I_T,H^1)}^2+
\|\bar{\mf{u}}^{k+1}\|_{L^\infty(I_T,H^2)}^2+\|\bar{\mf{u}}^{k+1}_t\|_{L^\infty(I_T,L^2)}^2\right)<\infty. $$
Hence, $({\varrho}^{k},{\mf{u}}^{k},{\mf{N}}^{k},\mf{u}_t^{k})$ is a Cauchy sequence in
$C(\bar{I}_T,H^1)\times C(\bar{I}_T,H^2)\times C(\bar{I}_T,H^1)\times C(\bar{I}_T,L^2)$, and
\begin{equation}\label{strongconvegneuN}
({\rho}^{k},{\mf{u}}^{k},{\mf{N}}^{k},\mf{u}_t^{k})\rightarrow   ( {\rho}, {\mf{u}},{\mf{N}},\mf{u}_t)
\end{equation}
strongly in $C(\bar{I}_T,H^1)\times C(\bar{I}_T,H^2)\times C(\bar{I}_T,H^1)\times C(\bar{I}_T,L^2)$.

Consequently, by \eqref{n0531}, \eqref{n0531presure} and \eqref{strongconvegneuN},
we easily verify that $((\varrho,\mf{N}),\mf{u})\in H_T\times V_T$ is a unique solution to the following problem with an associated
pressure $\tilde{q}$ enjoying the regularity \eqref{pressureesimta}:
\begin{equation*}\left\{\begin{array}{l}
\varrho_t+{\bf u}\cdot\nabla
 \varrho=-{\bf u}\cdot\nabla\bar{\rho} , \\[1mm]
(\varrho+\bar{\rho}){\bf u}_t+(\varrho+\bar{\rho}){\bf u}\cdot\nabla
{\bf u}+\nabla\tilde{q}-\mu \Delta \mathbf{u}= (\mf{N}+\bar{\mf{M}})\cdot\nabla \mf{N}-g \varrho \mf{e}_3,\\[1mm]
\mf{N}_t+\mf{u}\cdot\nabla \mf{N}=
(\mf{N}+\bar{\mf{M}})\cdot\nabla \mf{u},\\[1mm]
\mathrm{div}\mathbf{u}=0,\ \mathrm{div}\mathbf{N}=0 \end{array}\right.  \end{equation*}
with initial data
\begin{equation*}
(\varrho,\mathbf{u},\mathbf{N} )|_{t=0}=(\varrho_0,{\mathbf u}_0,\mathbf{N}_0).
\end{equation*}

Finally, if we define
$${q}:=\tilde{q}-|\mf{N}+\bar{\mf{M}}|^2/2, $$
and use the facts that
$$(\nabla \times(\mf{N}+\bar{\mf{M}}) )\times (\mf{N}+\bar{\mf{M}})=(\mf{N}+\bar{\mf{M}})\cdot\nabla \mf{N}-\nabla \frac{|\mf{N}+\bar{\mf{M}}|^2}{2},$$
and $$\nabla \times (\mf{u}\times (\mf{N}+\bar{\mf{M}}))= (\mf{N}+\bar{\mf{M}})\cdot\nabla \mf{u}-\mf{u}\cdot\nabla \mf{N},$$
we obtain Proposition \ref{pro:0401new}.
\vspace{4mm}

 \noindent\textbf{Acknowledgements.}  
The research of Fei Jiang was supported by the NSFC (Grant Nos. 11301083 and 11471134) and the NSF of Fujian Province of China (Grant No. 2014J01011),
the research of Song Jiang by the National Basic Research Program
under the Grant 2011CB309705, NSFC (Grant Nos. 11229101, 11371065) and the Beijing Center for Mathematics and Information Interdisciplinary Sciences.

\renewcommand\refname{References}
\renewenvironment{thebibliography}[1]{%
\section*{\refname}
\list{{\arabic{enumi}}}{\def\makelabel##1{\hss{##1}}\topsep=0mm
\parsep=0mm
\partopsep=0mm\itemsep=0mm
\labelsep=1ex\itemindent=0mm
\settowidth\labelwidth{\small[#1]}%
\leftmargin\labelwidth \advance\leftmargin\labelsep
\advance\leftmargin -\itemindent
\usecounter{enumi}}\small
\def\newblock{\ }
\sloppy\clubpenalty4000\widowpenalty4000
\sfcode`\.=1000\relax}{\endlist}

\end{document}